\documentstyle[amssymb,amsfonts,12pt]{amsart}
\newtheorem{theorem}{Theorem}[section]
\newtheorem{lemma}[theorem]{Lemma}
\newtheorem{corollary}[theorem]{Corollary}
\newtheorem{proposition}[theorem]{Proposition}
\newtheorem{definition}[theorem]{Definition}

\newtheorem{question}[theorem]{Question}

\theoremstyle{remark}
\newtheorem{remark}[theorem]{Remark}

\def\QSet{\mbox{\rm\kern.24em
\vrule width.03em height1.48ex depth-.051ex \kern-.26em Q}}

\def\E{{\mbox{\rm I\kern-.22em E}}}
\def\P{{\bf P}}
\def\D{{\bf D}}
\def\T{{\bf T}}
\def\S{{\bf S}}
\def\Z{{\bf Z}}
\def\R{{\bf R}}
\def\N{{\bf N}}
\def\C{{\bf C}}
\def\supp{{\operatorname{supp}}}
\def\size{{\operatorname{size}}}
\def\diam{{\operatorname{diam}}}

\def\th{{\operatorname{th}}}

\def\H{{\mathcal H}}
\def\P{{\mathcal P}}

\def\M{{\operatorname{M}}}
\def\BMO{{\operatorname{BMO}}}
\def\F{{\mathcal F}}
\def\L{{\mathcal L}}

\def\G{{\mathcal G}}
\def\O{{\mathcal S}}
\def\U{{\mathcal C}} 
\def\X{{\mathcal A}}
\def\B{{\mathcal B}}
\def\I{{\mathcal I}}

\def\dist{{\operatorname{dist}}}

\def\univ{{\operatorname{univ}}}

\def\bas{\begin{align*}}
\def\eas{\end{align*}}
\def\bi{\begin{itemize}}
\def\ei{\end{itemize}}
\newenvironment{proof}{\noindent {\bf Proof} }{\endprf\par}
\def \endprf{\hfill  {\vrule height6pt width6pt depth0pt}\medskip}
\def\emph#1{{\it #1}}

\begin{document}
\title{Breaking the duality in the return times theorem}

\author{Ciprian Demeter}
\address{Department of Mathematics, UCLA, Los Angeles CA 90095-1555}
\email{demeter@@math.ucla.edu}

\author{Michael Lacey}
\address{ School of Mathematics,
Georgia Institute of Technology, Atlanta,  GA 30332 USA}
\email{lacey@@math.gatech.edu}

\author{Terence Tao}
\address{Department of Mathematics, UCLA, Los Angeles CA 90095-1555}
\email{tao@@math.ucla.edu}

\author{Christoph Thiele}
\address{Department of Mathematics, UCLA, Los Angeles CA 90095-1555}
\email{thiele@@math.ucla.edu}
\thanks{The first author was supported by NSF Grant DMS-0556389}
\thanks{The second author was supported by an NSF Grant}
\thanks{The third author was supported by a grant from the Macarthur Foundation}
\thanks{The fourth author was supported by NSF Grant DMS-0400879}
\begin{abstract}
We prove Bourgain's Return Times Theorem for a range of exponents $p$ and $q$ that are outside the duality range. An oscillation result is used to
prove hitherto unknown almost everywhere convergence for the signed average analog of Bourgain's averages. 
\end{abstract}

\maketitle
\section{Introduction}

Almost everywhere convergence results for ergodic weighted averages of
various kinds typically are proved in two steps: first one proves
convergence for a small class of functions (typically $L^\infty$)
and then one proves a priori bounds for maximal operators which allow
to extend the almost everywhere convergence result to
larger classes of functions (typically $L^p$). In many instances,
both steps can require rather sophisticated analysis and offer their
own challenges. As the exponent $p$ is lowered, it gets increasingly
harder to prove $L^p$ bounds for maximal operators and there may be
several thresholds at which certain methods break down. It has been
recognized in \cite{La} that time frequency methods as
pioneered in  \cite{Car}, \cite{F}, \cite{LTBilH}, \cite{LTCar} give the strongest maximal theorems known
to date for the operators that they apply to. The purpose of
the current paper is to make time frequency methods available for a
much wider class of ergodic averages that have enjoyed some prominence
in ergodic theory in recent history. In particular we are able to break the
threshold of exponents in duality in Bourgain's Return Times Theorem
\cite{Bo5}, \cite{Bo15}. The methods in this paper are rather robust and typically
apply not only to standard averages but for example to signed and weighted
averages with Hilbert kernels as weights. Moreover, the method typically
provides a priori estimates for oscillation norms along with a priori
estimates for maximal operators, and thus abandons the need to prove
convergence for dense subclass ($L^\infty$) along with bounds for
maximal operators. In this paper an oscillation result is used to
prove hitherto unknown convergence for the signed average analog of
Bourgain's Return Times Theorem, and to provide a separate proof of Bourgain's theorem.
As in earlier works such as \cite{D}, \cite{DTT} and \cite{LTer}, our methods are almost entirely analytic in nature, however the results have independent interest from both an ergodic theoretic and harmonic analytic point of view.

Let ${\bf X}=(X,\Sigma,\mu, \tau)$ be a dynamical system, that is a Lebesgue space $(X,\Sigma,\mu)$ equipped with an invertible bimeasurable  measure preserving transformation $\tau:X\to X$. We recall that a complete probability space $(X,\Sigma,\mu)$ is called a \emph{Lebesgue space} if it is isomorphic with the ordinary Lebesgue measure space $([0,1),\L,m)$, where $\L$ and $m$ denote the usual Lebesgue algebra and measure (see \cite{Ha} for more on this topic). In particular, the $\sigma$-algebra $\Sigma$ (and hence all the spaces $L^p(X)$) will be separable, a property that will be used later to argue that a certain class of operators act measurably. The system ${\bf X}$ is called \emph{ergodic} if $A\in\Sigma$ and $\mu(A\vartriangle\tau^{-1}A)=0$ imply $\mu(A)\in\{0,1\}.$

In \cite{Bo5}, \cite{Bo15} (see also \cite{BFKO}) Bourgain proved the following result.

\begin{theorem}[Return times theorem, \cite{Bo5}, \cite{Bo15}]
\label{Bretthm}
For each function $f\in L^{\infty}(X)$ there is a universal set $X_0\subseteq X$ with $\mu(X_0)=1$, such that for each second dynamical system  ${\bf Y}=(Y,\F,\nu,\sigma)$, each $g\in L^{\infty}(Y)$ and each $x\in X_0$, the averages
$$\lim_{N \to \infty} \frac1{N}\sum_{n=0}^{N-1}f(\tau^nx)g(\sigma^ny)$$
converge $\nu$-almost everywhere.
\end{theorem}

If in the above theorem  $f$ is taken to be a constant function, one recovers the classical Birkhoff's pointwise ergodic theorem, see \cite{Bi}. However, Theorem ~\ref{Bretthm} is much stronger, in that it shows that given $f$, almost every sampling sequence $(f(\tau^nx))_{n\in\N}$ forms a system of universal weights for the  pointwise ergodic theorem.

Interest in results like Theorem ~\ref{Bretthm} can be traced back to the result of Wiener and  Wintner \cite{WW}, whose equivalent formulation is that for each integrable function $f$, almost every sampling sequence $(f(\tau^nx))_{n\in\N}$ is a universal system of weights for the mean ergodic theorem:

\begin{theorem}[Wiener-Wintner theorem, \cite{WW}]
\label{WWthm}
For each function $f\in L^{1}(X)$  there is a universal set $X_0\subseteq X$ with $\mu(X_0)=1$, such that for each $\theta\in[0,1)$ and each $x\in X_0$ the following averages converge
$$\lim_{N \to \infty} \frac1{N}\sum_{n=0}^{N-1}f(\tau^nx)e^{2\pi in\theta}.$$
\end{theorem}

This result also is an immediate consequence of Theorem ~\ref{Bretthm}. Indeed, for each $\theta\in[0,1)$ we can apply the theorem to the system ${\bf Y}$ consisting of  the interval $[0,1)$ equipped with the Lebesgue algebra and measure, together  with the transformation $\sigma y:=y+\theta\pmod 1,$  and to the function $g(y):=e^{2\pi iy}.$

An alternative proof of Theorem ~\ref{Bretthm}, based on the machinery of joinings, is due to Rudolph \cite{Ru}. The same author refines his techniques in \cite{Ru1} to prove a deep multiple return times theorem. H\"older's inequality and an elementary density argument show that Bourgain's theorem holds for $f\in L^{p}(X)$ and $g\in L^{q}(Y)$, whenever $1\le p,q\le \infty$ and $\frac1p+\frac1q\le 1$, see \cite{Ru} and also Section ~\ref{sec:approx} here. On the other hand, it has been recently proved by Assani, Buczolich and   Mauldin \cite{ABM} that this result fails when $p=q=1$:

\begin{theorem}\cite{ABM}
Let ${\bf X}=(X,\Sigma,\mu, \tau)$ 
 be an ergodic dynamical system. There exist a function $f\in L^1(X)$ and a subset $X_0\subseteq X$ of full  measure with the following property: for each $x_0\in X_0$ and for each ergodic dynamical system ${\bf Y}=(Y,\F,\nu,\sigma)$, there exists $g\in L^{1}(Y)$ such that the averages 
$$\lim_{N \to \infty} \frac1N\sum_{n=0}^{N-1}f(\tau^nx_0)g(\sigma^ny)$$
diverge for almost every y. 
\end{theorem}
The need for ergodicity in the above theorem is apparent from the observation that if either $\tau$ or $\sigma$ is the (nonergodic) identity transformation, then a positive result is easily seen to hold instead, for all integrable functions $f$ and $g$.

An interesting question arises on whether Theorem ~\ref{Bretthm} holds outside the duality range:

\begin{question}
\label{WWavethmqqqjy}
Do there exist indices $1< p,q<\infty$ with $\frac1p+\frac1q>1$ such that for each dynamical system ${\bf X}=(X,\Sigma,\mu, \tau)$ and each  $f\in L^{p}(X)$ there is a universal set $X_0\subseteq X$ with $\mu(X_0)=1$, such that for each second dynamical system  ${\bf Y}=(Y,\F,\nu,\sigma)$, each $g\in L^{q}(Y)$ and each $x\in X_0$, the averages
\begin{equation}
\frac1N\sum_{n=0}^{N-1}f(\tau^nx)g(\sigma^ny)
\end{equation}
converge $\nu$-almost everywhere?
\end{question}

 Supporting evidence for a positive result in this direction comes from the fact that the duality is indeed broken if either the weights or the test process is replaced with a sequence of i.i.d. random variables:

\begin{theorem}[I. Assani 2003, \cite{A2}, \cite{A3}]
\label{thm:hihihi0}

Let $(X_n)$ be a sequence of i.i.d. random variables with finite $p^{\th}$ moment for some $1<p\le \infty$, defined  on the probability space $(X,\Sigma, \mu)$. Then there exists a subset $X^{*}\subseteq X$ of full measure such that for each $x\in X^{*}$ the following holds:  for any  dynamical system  ${\bf Y}=(Y,\F,\nu,\sigma)$ and  $g\in L^q(Y)$, $1< q\le\infty$,  we have
$$\lim_{N\to\infty}\frac1N\sum_{n=0}^{N-1}X_n(x)g(\sigma^ny)=\E(X_0)\int g d\nu$$ for $\nu$-almost every $y$.
\end{theorem}

\begin{theorem}[I. Assani 1997, \cite{A1}; J. Baxter, R. Jones, M. Lin, J. Olsen 2003, \cite{BJLO}]
\label{thm:asdg}
Assume that either $p>1$ and $q=1$,  or $p=1$ and $q>1$. 
For each dynamical system ${\bf X}=(X,\Sigma, \mu,\tau)$ and each $f\in L^p(X)$ there is a set $X^{*}\subseteq X$ of full measure, such that for each sequence of $L^p$ i.i.d. random variables $Y_n$ defined on the probability space $(Y, \F,\nu)$ and each $x\in X^{*}$,
$$\lim_{N\to\infty}\frac1N\sum_{n=0}^{N-1}f(\tau^nx)Y_n(y)$$
exists for $\nu$-almost every $y$. 
\end{theorem}

Similar questions arise in the case of summation operators. We recall that the almost everywhere convergence of the ergodic truncated Hilbert transform 
\begin{equation}
\label{HilbserCot}
\lim_{N \to \infty} \sideset{}{'}\sum_{n=-N}^{N}\frac{f(\tau^nx)}{n}
\end{equation}
was proved by Cotlar \cite{Cot}.
The return times results for series are harder; the spectral theory and dynamics methods seem to be inapplicable to address the following question

\begin{question}
\label{WWserthmqqq}
Given  $1< p$ and  $1\le q\le \infty$, is it true that for each  dynamical system ${\bf X}=(X,\Sigma,\mu, \tau)$ and  each function $f\in L^{p}(X)$, there is a universal set $X_0\subseteq X$ with $\mu(X_0)=1$, such that for each second dynamical system  ${\bf Y}=(Y,\F,\nu,\sigma)$, each $g\in L^{q}(Y)$ and each $x\in X_0$, the series
\begin{equation}
\lim_{N \to \infty} \sideset{}{'}\sum_{n=-N}^{N}\frac{f(\tau^nx)g(\sigma^ny)}{n}
\end{equation}
converges $\nu$-almost everywhere?
\end{question} 

It has been shown in \cite{A5} that Question ~\ref{WWserthmqqq} has a negative answer when $p=1$, for arbitrary  $q$. Positive results are again known outside the duality range, in the special case when either the weights or the test process consist of  i.i.d. random variables, see \cite{A4}. However, no positive results were known for  Question ~\ref{WWserthmqqq} prior to this work, not even when $p=q=\infty$. We note that unlike the case of the averages, H\"older's inequality is of no use here due to the lack of summability of the sequence $(\frac1n)_{n\in\N}$.

We close this discussion with a parallel between return times results for averages and series. Spectral theory is an important component of all the four known proofs of Theorem ~\ref{Bretthm}. Three of them use purely dynamical (in particular non-Fourier-analytical) methods  and rely on the spectral decomposition according to which each  function can be decomposed  into a component with a purely discrete spectral measure plus a component with  continuous spectral measure. If $f_1$ and $g_1$ represent the continuous components of $f$ and $g$ while $f_2$ and $g_2$ are the discrete components, then the proof shows that the limit of the averages 
$$\lim_{N \to \infty} \frac{1}{N}\sum_{n=0}^{N-1}f_i(\tau^nx)g_j(\sigma^ny)$$ 
is $0$, as long as $1\in\{i,j\}$.
That is to say, the Kronecker factor (i.e the sub $\sigma$ algebra spanned by the eigenfunctions of the transformation) is characteristic for the (almost everywhere and norm) convergence of these averages.

This type of spectral analysis has not proven successful so far in proving convergence results for Hilbert series like the ones in Question ~\ref{WWserthmqqq}. The Kronecker factor is not expected to play the same role as in the case of averages. In particular, not even the series in ~\eqref{HilbserCot} will converge to zero for all functions with continuous spectrum. This suggests that, perhaps, the answer to these questions does not lie in dynamics, but rather in analytic methods.

In this paper we will answer affirmatively  Questions ~\ref{WWavethmqqqjy}  and will prove a similar result for  Question ~\ref{WWserthmqqq}. These theorems are described in detail in the next section.

\section{Notation and terminology}
\label{sec:4}
If $I\subseteq \R$ is an interval then  $c(I)$ denotes the center of $I$, $|I|$ denotes the length, and $CI$ is the interval with the same center and length $C$ times the length of $I$.
By $ 1_{A}$ we denote the characteristic function of the set $A\subseteq \R$, while for any  interval $I$, we define the weight function $$\chi_I(x):=\left(1+\frac{|x-c(I)|}{|I|}\right)^{-1}.$$
A tile $s$ is a rectangle $s=I_s\times\omega_{s}$ with $I_s$ some dyadic interval and $\omega_s$ some interval satisfying $|I_s|\cdot|\omega_s|=1$. 

The notation $a\lesssim b$ means that $a\le cb$ for some universal constant $c$, and $a\sim b$ means that  $a\lesssim b$ and  $b\lesssim a$.
These constants are allowed to depend on the exponents $p$ and $q$.
Sometimes we will  write $|f(x)|\lesssim \chi_I^M(x)$ with unspecified $M$ to indicate that this inequality holds for all $M\ge 1$, with implicit constant depending only on  $M$. Also, for each $1\le p< \infty$ we use the \emph{$p^{\th}$-power Hardy-Littlewood maximal operator} 
$$\M_p(f)(x):=(\sup_{r>0}\frac{1}{r}\int_{|t|\le r}|f(x+t)|^pdt)^{1/p}$$
and the $\BMO$ norm
$$ \|f\|_{\BMO(\R)} := \sup_I \frac{1}{|I|} \int_I \left|f - \frac{1}{|I|} \int_I f\right|$$
where $I$ ranges over all intervals.

The \emph{Fourier transform} of a function $f:\R\to\R$ is
$$\widehat{f}(\xi):=\F(f)(\xi)=\int f(x)e^{-2\pi i\xi x}dx,$$
while the inverse Fourier transform is
$$\check{f}(\xi):={\F}^{-1}(f)(\xi)=\int f(x)e^{2\pi i\xi x}dx.$$
Define the dilation, translation, and modulation operators 
\begin{equation*}  
\label{e.dilate}
\text{Dil}_s^p f (x) := s^{-1/p}f(x/s), 
\end{equation*}
\begin{equation*}  
\label{e.trans}
\text{Tr}_y f(x):=f(x-y),
\end{equation*}
\begin{equation*}
\text{Mod}_\theta f(x):= e^{2\pi i \theta x} f(x). 
\end{equation*}

\begin{definition}
For each $M\ge 0$, let $A(M)$ be some big universal constants, that will stay fixed throughout this paper.
A function $\phi_I$ is said to be \emph{$C$-adapted} to the interval $I$ if for each such\footnote{Actually, our proof will only require these decay bounds for a finite number of $M$, though the number of such $M$ can depend on exponents such as $p$.} $M\ge 0$
\begin{equation*}
|\phi_I(x)|\le A(M)C \frac1{|I|^{1/2}}\chi_I^{M}(x)
\end{equation*}
\begin{equation*}
|\frac{d}{dx}\phi_I(x)|\le A(M)C \frac1{|I|^{3/2}}\chi_I^{M}(x).
\end{equation*}
\end{definition}
The constant $C$ will vary throughout this paper and will always be specified  explicitly.
\begin{definition}
A function $\phi_I$ is said to have the mean zero property with respect to a frequency $c$ if 
$$\int\phi_I(x)e^{-ixc}=0.$$ 
\end{definition}

\section{Main results and high-level overview of the proof}
\label{sec:2}

Our first result here gives an affirmative answer to Question ~\ref{WWavethmqqqjy}, by extending  Bourgain's Return Times theorem to the range $1<p\le \infty$ and $q\ge 2$.

\begin{theorem}
\label{Bretthmour}
Let $1<p\le \infty$ and $q\ge 2$ be some arbitrary indices.
For each function $f\in L^{p}(X)$  there is a universal set $X_0\subseteq X$ with $\mu(X_0)=1$, such that for each second dynamical system  ${\bf Y}=(Y,\F,\nu,\sigma)$, each $g\in L^{q}(Y)$  and each $x\in X_0$, the averages
$$\frac1{N}\sum_{n=0}^{N}f(\tau^nx)g(\sigma^ny)$$
converge $\nu$-almost everywhere.
\end{theorem}

Given the convergence for $L^{\infty}$ functions $f$ and $g$, an approximation argument like in Theorem ~\ref{thm:approx6y} will immediately prove the above, once we establish the following maximal inequality:

\begin{theorem}
\label{otnumber1}
For each dynamical system ${\bf X}=(X,\Sigma,\mu, \tau)$, each $1<p\le \infty$ and each $f\in L^p(X)$
\begin{equation}
\label{hgdhgdhgh769878}
\|\sup_{(Y,\F,\nu,\sigma)}\sup_{\|g\|_{L^2(Y)=1}}\|\sup_{N}|\frac1{N}\sum_{n=0}^{N}f(\tau^nx)g(\sigma^ny)|\|_{L^2_y(Y)}\|_{L^p_x(X)}\lesssim \|f\|_{L^p(X)},
\end{equation}
where the first supremum in the inequality above is taken over all dynamical systems ${\bf Y}=(Y,\F,\nu,\sigma)$.  Here we have subscripted some of our $L^p$ norms to clarify the variable being integrated over.
\end{theorem}
\begin{remark}
The measurability in both inequality ~\eqref{hgdhgdhgh769878} and in inequality ~\eqref{hgdhgdhgh769878a} from below is proved by an application of Conze's principle (Theorem ~\ref{thm:Conze}) and the separability of each $L^2(Y)$, and the reader is referred to the proof of Theorem ~\ref{thm:approx6y} for details.
\end{remark}
Inequality ~\eqref{hgdhgdhgh769878} is only new for $1< p\le 2$. When $p>2$  it is an immediate consequence of H\"older's inequality and the boundedness of the ergodic maximal function in every $L^p$, $p>1$. 

The analog of Theorem ~\ref{otnumber1} for series also holds:
\begin{theorem}
\label{otnumber2}
For each dynamical system ${\bf X}=(X,\Sigma,\mu, \tau)$, each $1<p< \infty$ and each $f\in L^p(X)$
\begin{equation}
\label{hgdhgdhgh769878a}
\left\|\sup_{(Y,\F,\nu,\sigma)}\sup_{\|g\|_{L^2(Y)=1}}\|\sup_{N}|\sideset{}{'}\sum_{n=-N}^{N}\frac{f(\tau^nx)g(\sigma^ny)}{n}|\|_{L^2_y(Y)}\right\|_{L^p_x(X)}\lesssim \|f\|_{L^p(X)},
\end{equation}
where the first supremum in the inequality above is taken over all dynamical systems ${\bf Y}=(Y,\F,\nu,\sigma)$.
\end{theorem}
Note that no particular case of the maximal inequality  ~\eqref{hgdhgdhgh769878a} was previously known. It is also worth observing  the lack of applicability of H\"older's inequality in this context. The inequalities ~\eqref{hgdhgdhgh769878} and ~\eqref{hgdhgdhgh769878a} are obtained via standard transfer methods from the following general result, as explained in the Section ~\ref{sec:transfer}.

\begin{theorem} 
\label{t.returntime}  
Let $K:\R\to\R$ be an $L^2$ kernel satisfying the requirements:
\begin{align}
\label{reqkernel1}
\widehat{K}&\in C^{\infty}(\R\setminus\{0\}) \\
\label{reqkernel2}
|\widehat{K}(\xi)|&\lesssim \min\{1,\frac{1}{|\xi|}\} \quad \forall \xi\not=0\\
\label{reqkernel3}
|\frac{d^n}{d\xi^n}\widehat{K}(\xi)|&\lesssim \frac{1}{|\xi|^n}\min\{|\xi|,\frac1{|\xi|}\} \quad \forall \xi\not=0,\;n\ge 1.
\end{align}
Then the following inequality holds for each $1< p<\infty$ 
\begin{equation}  
\label{e.returntime}
\left\|\sup_{\| g\|_{L^2(\R)}=1}\left\|\sup_{k\in \Z}|\frac{1}{2^k}\int f(x+y)g(z+y)K(\frac{y}{2^k}) dy|
\right\|_{L^2_z(\R)} \right\|_{L^p_x(\R)}\lesssim \|f\|_{L^p(\R)}.
\end{equation}
\end{theorem}
\begin{remark}
Due to the fact that $K\in L^2$, the quantity 
$$\frac{1}{2^k}\int f(x+y)g(z+y)K(\frac{y}{2^k}) dy$$
is defined for each $g\in L^2$ and every $x$ and $z$, assuming $f$ is an $L^{\infty}$ function with  bounded support. Inequality ~\eqref{e.returntime} will be  proved with this extra requirement about  $f$, then density arguments will provide it with  a meaning for all $f\in L^{p}.$ It further follows that for each $x\in\R$
the quantity 
$$\sup_{\| g\|_{L^2(\R)}=1}\left\|\sup_{k\in \Z}|\frac{1}{2^k}\int f(x+y)g(z+y)K(\frac{y}{2^k}) dy|\right\|_{L^2_z(\R)}$$ is well defined and  gives rise to a measurable function of $x$.
\end{remark}
\begin{remark} 
A somewhat similar, yet distinct operator is the following  bilinear maximal function for which bounds are proved in \cite{La}:
$$B^{*}(f,g)(x)=\sup_{k\in \Z}\left|\frac{1}{2^k}\int f(x+y)g(x-y)K(\frac{y}{2^k}) dy\right|.$$
While the functions $f$ and $g$ play a symmetric role in the above, their contribution to the return times operator in inequality ~\eqref{e.returntime} is significantly different.

Also, unlike in the case of the  bilinear maximal function, the signs of $y$ in the innermost expression in the left hand side of ~\eqref{e.returntime} have no deep significance at all. More generally, a simple scaling-dilation argument  shows that  inequalities ~\eqref{e.returntime}  with $f(x+ay)g(z+by)$ are all equivalent, for each choice of $a,b\not=0$.
\end{remark}

One immediate consequence of the above result is the following. 
\begin{corollary} 
\label{t.returntimeave}  
For each $f\in L^p(\R)$ we have   
\begin{equation}  
\label{e.returntimeave}
\left\|\sup_{\| g\|_{L^2(\R)}=1}\|\sup_{t>0}\frac{1}{2t} \int_{-t}^t|f(x+y)g(z+y)|dy\|_{L^2_z(\R)}\right\|_{L^p_x(\R)}\lesssim \|f\|_{L^p(\R)},\;\;1<p<\infty.
\end{equation}
\end{corollary} 

The  corollary is trivial for $p>2$, by H\"older's inequality. To see how the result for general $p$ follows from that of Theorem  ~\ref{t.returntime}, choose $K$ to be some positive Schwartz function and note that it suffices to assume that $f$ and $g$ are  positive and also to restrict the supremum in ~\eqref{e.returntimeave} to dyadic values of $t$.

The second corollary is the analog of the first one for singular integrals.
\begin{corollary} 
\label{t.returntimeser}  
For each $f\in L^{\infty}(\R)$ with finite support we have  
\begin{equation}  
\label{e.returntimeser}
\|\sup_{\| g\|_{L^2(\R)}=1}\|\sup_{t>0}|\int_{|y|>t}f(x+y)g(z+y)\frac{dy}{y}|\|_{L^2_z(\R)}\|_{L^p_x(\R)}\lesssim \|f\|_{L^p(\R)},\;\;1<p<\infty.
\end{equation}
\end{corollary}
Note again that the integral above is defined for each $g\in L^2$ and each $x$ and $z$, due to the kernel $\frac1y1_{\{|y|>1\}}$ being in $L^2$. 
Consider a $C^{\infty}(\R)$ kernel such that $K(y)=\frac1{y}$ for $|y|\ge 1$. The proof of the above corollary follows from the following two observations. On the one hand, by using Corollary ~\ref{t.returntimeave} it suffices to prove  Corollary ~\ref{t.returntimeser} with $K$ replacing the rough kernel $\frac{1}{y}1_{\{|y|>1\}}$ and with the supremum restricted to dyadic values of $t$. On the other hand, it is an easy exercise to prove  that $K$ satisfies the requirements of Theorem ~\ref{t.returntime}. 

As far as Question ~\ref{WWserthmqqq} is concerned, we remark that Theorem ~\ref{otnumber2} can not provide any answer to it. The reason is that a dense class result is missing. It turns out however that the techniques used in Theorem ~\ref{otnumber2} can be refined  to prove the following analog for series of Bourgain's Return Times theorem. 

\begin{theorem}
\label{Bretthmseriesk=0}
For each function $f\in L^{\infty}(X)$  there is a universal set $X_0\subseteq X$ with $\mu(X_0)=1$, such that for each second dynamical system  ${\bf Y}=(Y,\F,\nu,\sigma)$, each $g\in L^{\infty}(Y)$  and each $x\in X_0$, the series
$$\sideset{}{'}\sum_{n=-N}^{N}\frac{f(\tau^nx)g(\sigma^ny)}{
n}$$
converges $\nu$-almost everywhere.
\end{theorem}

Now Theorems ~\ref{otnumber2} and  ~\ref{Bretthmseriesk=0} together with an approximation argument as in Theorem ~\ref{thm:approx6y} lead to  the following corollary.

\begin{corollary}
\label{Bretthmserieskcor=0}
Let $1<p\le \infty$ and $q\ge 2$ be some arbitrary indices.
For each function $f\in L^{p}(X)$  there is a universal set $X_0\subseteq X$ with $\mu(X_0)=1$, such that for each second dynamical system  ${\bf Y}=(Y,\F,\nu,\sigma)$, each $g\in L^{q}(Y)$  and each $x\in X_0$, the series
$$\sideset{}{'}\sum_{n=-N}^{N}\frac{f(\tau^nx)g(\sigma^ny)}{n}$$
converge $\nu$-almost everywhere.
\end{corollary} 
\begin{remark}
It actually turns out that the same methods can be used to give yet another proof\footnote{The proofs in \cite{BFKO}, \cite{Ru} and \cite{Ru1} use dynamics. Bourgain's original argument \cite{Bo5}, \cite{Bo15}, uses classical Fourier analysis geared towards getting entropy estimates for multipliers.  The proof along the techniques developed in our paper, while inspired by more recent developements in  time-frequency harmonic analysis, shares similarities with Bourgain's argument; in particular, a special case of Theorem \ref{thm:gh098cfhh} here also played a crucial role in Bourgain's original argument.} of  Bourgain's Return Times Theorem ~\ref{Bretthm}, see Section ~\ref{reproveBourgain}. 
\end{remark}

Choose  ${\bf Y}$  to be the interval $[0,1)$ equipped with the Lebesgue algebra and measure together  with the transformation $\sigma(y):=y+\theta\pmod 1,$ while $g(y):=e^{2\pi iy}.$ The above corollary applied to the dynamical system ${\bf Y}$ provides  the following Wiener-Wintner result for series
\begin{corollary}
\label{Cor:ghhnnir45}
Given $1< p\le \infty$,  for each dynamical system ${\bf X}=(X,\Sigma,\mu, \tau)$ and  each function $f\in L^{p}(X)$ there is a universal set $X_0\subseteq X$ with $\mu(X_0)=1$, such that for each $\theta\in[0,1)$ and each $x\in X_0$ the following series converges
$$\sideset{}{'}\sum_{n=-N}^{N}\frac{f(\tau^nx)}{n}e^{2\pi in\theta}.$$
\end{corollary}

A separate proof of the above result appears also in \cite{LTer}. The methods used there are not strong enough to address the rest of the results obtained in this paper. 

Since in general only quantitative inequalities transfer from harmonic analysis to ergodic theory, in order to prove Theorem ~\ref{Bretthmseriesk=0} via a transfer argument, the almost everywhere convergence needs to be quantified in some way. Our approach relies on proving an oscillation inequality, which will be shown to imply\footnote{It will become clear in Section ~\ref{sub:helconv} that the result of Theorem ~\ref{t.returntimeo} for any particular $p$ suffices to imply Theorem ~\ref{Bretthmseriesk=0}.} almost everywhere convergence in Section ~\ref{sec:transfer}.\footnote{This type of approach has been used before in ergodic theory, see for example \cite{Bo1}.}

\begin{theorem} 
\label{t.returntimeo}  
Let $K:\R\to\R$ be an $L^2$ kernel satisfying ~\eqref{reqkernel1}, ~\eqref{reqkernel2} and ~\eqref{reqkernel3}.
Then for each $1<p<\infty$ there is $0<\epsilon(p)<\frac12$ such that the following holds: for each $d=2^{1/n}$, $n\in \N$, and for each finite sequence of integers $k_1<k_2<\ldots<k_J$ 
\begin{equation*}  
\label{e.returntimeo}
\left\|\sup_{\| g\|_{L^2(\R)}=1}\left\|(\sum_{j=1}^{J-1}\sup_{k\in\Z\atop{k_j\le k<k_{j+1}}}|\int f(x+y)g(z+y)(\operatorname{Dil}_{d^k}^1K-\operatorname{Dil}_{d^{k_{j+1}}}^1K)(y)\ dy|^2)^{1/2}\right\|_{L^2_z(\R)}\right\|_{L^p_x(\R)}
\end{equation*}
\begin{equation*}
\lesssim J^{\frac12-\epsilon(p)}\|f\|_{L^p(\R)},
\end{equation*}

with the implicit constants depending only on $n$ and $p$.
\end{theorem}

This theorem is a consequence of two distinct results of dyadic analysis. The first one, Theorem ~\ref{t.returntimeo1}, is the particular case $d=2$ of the above and captures the main difficulty of the problem. The second one, Theorem ~\ref{t.returntimeo2}, is a square function estimate and is meant to control error terms.

To understand better the connection between Theorems ~\ref{t.returntimeo}, ~\ref{t.returntimeo1} and ~\ref{t.returntimeo2} we introduce some notation.
Let $h:(0,\infty)\to\C$. Let also $k_1<\ldots<k_{J}$ be as in Theorem ~\ref{t.returntimeo}  and define integers $a_1\le\ldots\le a_{J}$ such that $a_jn\le k_j<(a_j+1)n$. Then observe that 
\begin{align*}
(\sum_{j=1}^{J-1}\sup_{k\in\Z\atop{k_j\le k<k_{j+1}}}|h(\frac{k}{n})-h(\frac{k_{j+1}}{n})|^2)^{1/2}
& \lesssim \sum_{i=0}^{n-1}(\sum_{j=1}^{J-1}\sup_{k\in\Z\atop{a_j\le k<a_{j+1}}}|h(k+\frac{i}{n})-h(a_{j+1}+\frac{i}{n})|^2)^{1/2} \\
&\quad +\sum_{i,j=0\atop{i\not=j}}^{n-1}(\sum_{k\in\Z}|h(k+\frac{i}{n})-h(k+\frac{j}{n})|^2)^{1/2}.
\end{align*}

Using this inequality and a dilation argument, Theorem ~\ref{t.returntimeo} will follow immediately from the following two results.

\begin{theorem} 
\label{t.returntimeo1}  
Let $K:\R\to\R$ be an $L^2$ kernel satisfying ~\eqref{reqkernel1}, ~\eqref{reqkernel2} and ~\eqref{reqkernel3}.
Then for each $1<p<\infty$ there is $0<\epsilon(p)<\frac12$ such that  for each finite sequence of integers $k_1<k_2<\ldots<k_J$ 
\begin{equation*}  
\label{e.returntimeo1}
\left\|\sup_{\| g\|_{L^2(\R)}=1}\|\left(\sum_{j=1}^{J-1}\sup_{k\in\Z\atop{k_j\le k<k_{j+1}}}|\int f(x+y)g(z+y)(\operatorname{Dil}_{2^k}^1K-\operatorname{Dil}_{2^{k_{j+1}}}^1K)({y})\ dy|^2\right)^{1/2}\|_{L^2_z(\R)}\right\|_{L^p_x(\R)}
\end{equation*}
\begin{equation*}
\lesssim J^{\frac12-\epsilon(p)}\|f\|_{L^p(\R)},
\end{equation*}
with the implicit constants depending only on  $p$.
\end{theorem}

\begin{theorem} 
\label{t.returntimeo2}  
Let $K:\R\to\R$ be an $L^2$ kernel satisfying ~\eqref{reqkernel1}, ~\eqref{reqkernel2} and ~\eqref{reqkernel3} and the extra requirement
\begin{equation}
\label{reqkernel4}
|\widehat{K}(\xi)|\lesssim |\xi|.
\end{equation}
Then for each $1<p<\infty$ the following inequality holds
\begin{equation*}  
\label{e.returntimeo2}
\left\|\sup_{\| g\|_{L^2(\R)}=1}\|\left(\sum_{k\in \Z}|\int f(x+y)g(z+y)(\operatorname{Dil}_{2^{k}}^1 K) ({y})\ dy|^2\right)^{1/2}\|_{L^2_z(\R)}\right\|_{L^p_x(\R)}
\lesssim \|f\|_{L^p(\R)},
\end{equation*}
with the implicit constants depending only on  $p$.
\end{theorem}

Our approach to theorems ~\ref{t.returntime},  ~\ref{t.returntimeo1} and ~\ref{t.returntimeo2} relies on using time-frequency techniques to bound discrete model operators. This amounts to decomposing the time-frequency plane into dyadic rectangles $s=I_s\times \omega_s$ (also called tiles), associated with highly localized  wave packets $\phi_s(x,\theta),\varphi_s(x)$. The decomposition is guided by the nature of the operator under investigation, and the goal is to reduce the proof of its boundedness to that of the discrete model sums $$\sum_s\langle f,\varphi_s\rangle\phi_s(x,\theta)$$ in some appropriate norm. Our proof of Theorem ~\ref{t.returntime} has emerged from the discovery of striking connections  between the model operator for the return times operator  and the \emph{Carleson-Hunt's operator}
\begin{equation*}  
\label{e.Car}
Cf(x,\theta):=p.v.\int_{\R}\frac{f(x-y)}{y}e^{i\theta y}dy
\end{equation*}
which controls the convergence of the Fourier series.
To  clarify this connection  we introduce some notation. For each $1\le p\le \infty$, the \emph{$M_p$ multiplier norm} of a function  $m:\R\to\R$ is defined as
$$\|m\|_{M_p(\R)}= \|m(\theta)\|_{M_{p,\theta}(\R)}:=\sup_{\|h\|_p=1}\|\int m(\theta)\widehat{h}(\theta)e^{2\pi i\theta x}d\theta\|_{L^p_x(\R)}.$$ 
Of course the $M_2(\R)$ norm is just the $L^\infty(\R)$ norm, $\|m\|_{M_2(\R)}=\|m\|_{L^{\infty}(\R)}.$
Similarly, the \emph{$M_p^{*}$ norm} of a sequence of multipliers $m_k:\R\to\R$ is defined as
$$\|(m_k)_{k \in \Z}\|_{M_p^{*}(\R)} = \|(m_k(\theta))_{k \in \Z}\|_{M_{p,\theta}^{*}(\R)} 
 :=\sup_{\|h\|_p=1}\left\|\sup_k|\int m_k(\theta)\widehat{h}(\theta)e^{2\pi i\theta x}d\theta|\right\|_{L^p_x(\R)}.$$ 
The celebrated theorem of Carleson-Hunt asserts the following:
\begin{theorem}
\label{Carleson-Hunt}
For each $1<p<\infty$ and each $f\in L^p(\R)$,
$$\left\|\|Cf(x,\theta)\|_{L^{\infty}_\theta(\R)} \right \|_{L^p_x(\R)}\lesssim \|f\|_{L^p(\R)}$$
or equivalently
$$\left\|\|Cf(x,\theta)\|_{M_{2,\theta}(\R)} \right \|_{L^p_x(\R)}\lesssim \|f\|_{L^p(\R)}$$
\end{theorem}

It turns out that there is an appropriate choice of wave packets $\phi_s$ and $\varphi_s$ such that  Theorem ~\ref{Carleson-Hunt} can be reduced to showing that 
\begin{equation}
\label{discterecarleson}
\left\|\|\sum_{s}\langle f,\varphi_s\rangle\phi_s(x,\theta)\|_{M_{2,\theta}(\R)}\right\|_{L^p_x(\R)}\lesssim \|f\|_{L^p(\R)},
\end{equation}
 while Theorem ~\ref{t.returntime} can be reduced to showing that 
\begin{equation}
\label{discterereturntimes}
\|\|(\sum_{s:|I_s|\ge 2^k}\langle f,\varphi_s\rangle\phi_s(x,\theta))_{k \in \Z}\|_{M_{2,\theta}^{*}(\R)}\|_{L^p_x(\R)}\lesssim \|f\|_p.
\end{equation}

The proof of Theorem ~\ref{t.returntimeo1} relies on the same techniques as the ones utilized in Theorem ~\ref{t.returntime}, with an extra twist created by the  oscillations of the operators in question. The main new ingredient here is Theorem ~\ref{Bourgeneralo}, whose estimates incorporate both the maximal and the oscillatory behavior of the multiplier. 

In contrast, Theorem  ~\ref{t.returntimeo2} does not encode any maximal or oscillatory behavior. Its proof does not need any new ingredients, other than the ones we use to produce an (implicit) proof of the Carleson-Hunt theorem. 
  
The main novelty of our approach in this paper resides in getting \emph{local} type of estimates for the model operator, as opposed to proving \emph{global} estimates via dualization. This latter strategy was successful in dealing with maximal operators of similar complexity, as those in \cite{D}, \cite{DTT}, \cite{La}. Our search for this new type of approach was guided by the the nature of the $M_2^*$ norm, which makes the dualization of ~\eqref{discterereturntimes} extremely  hard to handle. We thus had to develop a set of techniques that do not involve the dual of the $M_2^*$ norm. We note that the $M_2$ norm is much more amenable to dualization. This fact was observed in  \cite{LTCar} in the context of the Carleson-Hunt operator, where dualization of the $M_2$ norm was used to create an interplay between energy and mass.

Here is an overview of our proof of inequality ~\eqref{discterereturntimes}. In Section ~\ref{sec:5} we indicate how to reduce  theorems ~\ref{t.returntime},  ~\ref{t.returntimeo1} and  ~\ref{t.returntimeo2} to  similar statements about discrete model operators. The details for our main result, Theorem  ~\ref{t.returntime}, are as follows. 
For each scale $k\in \Z$ we further decompose the model  operator $\sum_{s:|I_s|\ge 2^k}\langle f,\varphi_s\rangle\phi_s(x,\theta)$ into the sum of two distinct operators with good frequency localization. 

The first one is controlled by a weighted version of the aforementioned maximal multiplier result of Bourgain, in which the multiplier assumes different values depending on $k$ and on the frequency base point. The proof of this result is presented in Section ~\ref{sec:7} and its later application depends on variational estimates proved in Section ~\ref{sec:8}.
The second operator is essentially a composition of the original Bourgain's maximal operator  and Carleson's operator, and as a consequence its boundedness depends on the boundedness of these two fundamental operators.  

Our analysis of the return times operator is then guided by time localization, in  that for each $x$ on the time axis we split the contribution coming from various trees in terms of their spatial localization with respect to $x$. We then get pointwise -rather than global $L^p$ norm- estimates for the model operator at each $x$ outside an appropriately chosen exceptional set. One immediate advantage of this type of localization is that it reduces substantially  the combinatorial difficulty of organizing the trees  into structured subcollections. Indeed, the contribution to a given $x$ on the time axis will essentially come from just one  stack of trees. The fact that $x$ is chosen outside the exceptional set will guarantee  control both over the number of trees in the stack (which makes the weighted Bourgain's multiplier result effective, see Section ~\ref{sec:7}) and over the size of the weights (via $\BMO$ estimates, see Section ~\ref{sec:8}). The remaining details of the proof are then presented in the last two sections of the paper.  

This new method of estimating the model operator locally, as opposed to the previously employed global approach, has first led us to a new proof of the Carleson-Hunt theorem. One which is in the spirit of Carleson's original argument in that it uses energy but not mass, however it uses a completely different mathematical language and set of tools. This proof is incorporated in the main argument, and is used to control the second operator mentioned above.   

An approach to the  Return Times theorem in the case $1<q<2$ along the lines  of Theorem ~\ref{t.returntime}  would involve estimates both 
on the $M_q^{*}$ norm of the weighted  Bourgain's maximal multilinear operator in Theorem ~\ref{thm:gh098cfhh} and on the $M_q$  norm of the model sums associated with Carleson's operator. 
Crucial to our proof of the case $q=2$ in Theorem ~\ref{t.returntime} is the fact that the $M_2^{*}$ norm of the first operator is small as a function of the number $L$ of frequency basepoints\footnote{The bound obtained in Theorem ~\ref{thm:gh098cfhh} is of the order $L^{\epsilon}$, for arbitrarily small $\epsilon>0$. Any improvement over the trivial bound of $L^{1/2}$ produces  positive results for some range of $p<2$ and the fact that the bound is actually $L^{\epsilon}$ extends the result to the full range $1<p<\infty$. While the $L^{\epsilon}$ bound suffices for our applications here, it would be interesting to know its correct order of magnitude.}. The  $M_q^{*}$ norm is significantly larger when $q\not=2$. More precisely, it is shown in Section ~\ref{sec:7} that  this norm is  at least of the order of $L^{|1/2-1/q|}$ for each $q\in (1,2)\cup(2,\infty)$.

On the other hand, the $L$ dependency of the $M_q^{*}$ norm  for $1\le q<2$ is of at most $L^{2/q-1}$, which is what one gets by interpolating with the $M_1^{*}$ norm. Even with this large bounds our methods still seem to produce partial results in Theorem ~\ref{t.returntime} for other values of $q$, assuming  good control over the the $M_q$  norm of the model sums associated with Carleson's operator\footnote{This is currently investigated by the last two authors here together with other authors.}. This will appear elsewhere.

\section{The approximation argument}
\label{sec:approx}

Let $(Y_i,\F_i,\nu_i)$, $i=1,2$, be some arbitrary Lebesgue spaces. Denote by $\C(Y_i)$ the family of all the $\nu_i$-measure preserving transformations on $Y_i$. Equip $\C(Y_i)$ with the topology of weak convergence, in which $\tau_{s}\to \tau$ if and only if $\nu_i(\tau_sA\Delta\tau A)\to 0$ for each $A\in \F_i$. We will also denote by $\C(Y_1,Y_2)$ the set of all invertible, bimeasurable transformations $\beta:Y_1\to Y_2$ which take the measure $\nu_1$ to the measure $\nu_2$. The following result is due to Halmos \cite{Ha}.
\begin{lemma}
\label{lem:Halmos}
If ${\bf Y_1}=(Y_1,\F_1,\nu_1,\sigma_1)$ is ergodic then the set
$$\{\beta\sigma_1\beta^{-1}, \beta\in \C(Y_1,Y_2)\}$$
is dense in $\C(Y_2)$ in the weak topology.
\end{lemma}

Consider now  a sequence $S_N$ of weighted operators acting on the measurable functions in each system, according to the formula
$$S_Ng(y) :=\sum_{n=r_N}^{p_N}w_{N,n}g(\sigma_i^ny),$$
where the weights $\{w_{N,n}\}$ are arbitrary complex numbers. Denote by $S^{*}$ the maximal operator $S^{*}g=\sup_{N}|S_{N}g|$.
The following version of the so called Conze's principle is a consequence of the above lemma (see \cite{Co} for a similar version of this result).

\begin{theorem}[Conze's Principle]
\label{thm:Conze}

Let ${\bf Y}_i=(Y_i,\F_i, \nu_i, \sigma_i)$, $i=1,2$,  be two dynamical systems, with  ${\bf Y_1}$  ergodic. Then for each $1\le q<\infty$
$$\sup_{\|g\|_{L^q(Y_1)}=1}\|S^{*}g\|_{L^q(Y_1)}\ge \sup_{\|g\|_{L^q(Y_2)}=1}\|S^{*}g\|_{L^q(Y_2)}.$$
In particular, if both systems are ergodic then the left and right hand sides  are equal.
\end{theorem}

We use this to prove the following general approximation result.

\begin{theorem}[The approximation argument]
\label{thm:approx6y}
Fix some $1\le p,q< \infty$ and consider the dynamical systems  ${\bf X}=(X,\Sigma,\mu, \tau)$ and ${\bf Y}=(Y,\F,\nu,\sigma)$, where the second one is assumed to be ergodic.  Consider a sequence of bilinear operators defined as $$T_N(f,g)(x,z):=\sum_{n=-N}^{N}w_{N,n}f(\tau^nx)g(\sigma^nz)$$ for each dynamical system ${\bf Z}=(Z,\Upsilon, m, \rho)$, each $f\in L^p(X)$ and $g\in L^q(Z)$. Assume that 
\begin{equation}
\label{eq:approx6y}
\left\|\sup_{\|g\|_{L^q(Y)=1}}\|\sup_{N}|T_N(f,g)(x,y)|\|_{L^q_y(Y)}\right\|_{L^p_x(X)}\lesssim \|f\|_{L^p(X)},
\end{equation}
Assume also that for each function $f\in L^{\infty}(X)$ there is a universal set $X_0\subseteq X$ with $\mu(X_0)=1$, such that for each ${\bf Z}=(Z,\Upsilon, m, \rho)$, each $g\in L^{\infty}(Z)$ and each $x\in X_0$, the sequence
$$T_N(f,g)(x,z)$$
converge for $m$-almost every $z$.
Then the last statement above  also holds for each $f\in L^p(X)$ and each $g\in L^q(Z)$.
\end{theorem}
\begin{proof}
For each $f\in L^p(X)$ and each $x\in X$ define
$$R^{*}f(x):=\sup_{{\bf Z}}\sup_{\|g\|_{L^q(Z)=1}}\|\sup_{N}|T_N(f,g)(x,z)|\|_{L^q_z(Z)},$$
where the first supremum above is taken over all dynamical systems ${\bf Z}$. Note first that Theorem \ref{thm:Conze} implies that 
$$R^{*}f(x)=\sup_{\|g\|_{L^q(Y)=1}}\|\sup_{N}|T_N(f,g)(x,y)|\|_{L^q_y(Y)}.$$
Second, for each  $g\in L^q(Y)$ the quantity
$$\|\sup_{N}|T_N(f,g)(x,y)|\|_{L^q_y(Y)}$$
gives rise to a measurable function of $x$, by Fubini's theorem.
Third, since $L^q(Y)$ is separable it follows that for each $x$ the latter supremum can  be taken over a fixed countable family of functions $g_n$ which is dense in $L^2(Y)$. With these  observations, the fact that $R^{*}f(x)$ is a measurable function of $x$ follows immediately. Moreover, \eqref{eq:approx6y}  implies that
$$\|R^{*}f(x)\|_{L^p(X)}\lesssim \|f\|_{L^p(X)}.$$
Fix $f\in L^p(X)$. Let $f_i\in L^{\infty}(X)$ be such that $\|f-f_i\|_{L^p(X)}\to 0$. For each $i$ denote by $X_{0,i}$ the universal set corresponding to $f_i$. Define $X_0^{0}:=\cap X_{0,i}$ and note that it has full measure. For each dynamical system ${\bf Z}$ as above and for each $g\in L^q(Z)$ let $g_i\in L^{\infty}(Z)$ be such that $\|g_i\|_{L^q(Z)}\le 2$ and  $\|g-g_i\|_{L^q(Z)}\to 0.$
Now for each $x\in X_0^{0}$ and each $g\in L^q(Z)$ with $\|g\|_{L^q(Z)}=1$ we have
\begin{align*}
\|\limsup_{N,M\to\infty}&|T_N(f,g)(x,z)-T_M(f,g)(x,z)|\|_{L^q_z(Z)}\\&\le 2\inf_{i}\|\sup_{N}|T_N(f-f_i,g_i)(x,z)|\|_{L^q_z(Z)}+ 2\inf_{i}\|\sup_{N}|T_N(f,g-g_i)(x,z)|\|_{L^q_z(Z)}\\&\le 4\inf_{i}R^{*}(f-f_i)(x)+2\inf_{i}R^{*}f(x)\|g-g_i\|_{L^p(Z)}\\&= 4\inf_{i}R^{*}(f-f_i)(x)
\end{align*}

We deduce that
\begin{align*}
\|\sup_{{\bf Z}}\sup_{\|g\|_{L^q(Z)=1}}\|\limsup_{N,M\to\infty}&|T_N(f,g)(x,z)-T_M(f,g)(x,z)|\|_{L^q_z(Z)}\|_{L^p_x(X)}\\&\le 4\inf_{i}\|R^{*}(f-f_i)(x)\|_{L^p_x(X)}\\&\lesssim \inf_{i}\|f-f_i\|_{L^p(X)}=0
\end{align*}
The universal set $X_0$ associated with $f$ is obtained as the intersection between the set $X_0^{0}$ and the set of those $x\in X$ for which 
$$\sup_{{\bf Z}}\sup_{\|g\|_{L^q(Z)=1}}\|\limsup_{N,M\to\infty}|T_N(f,g)(x,z)-T_M(f,g)(x,z)|\|_{L^q_z(Z)}=0.$$
\end{proof}

\section{Transfer to ergodic theory}
\label{sec:transfer}
 
We first sketch the argument on how inequalities  ~\eqref{e.returntimeave} and ~\eqref{e.returntimeser} imply their counterparts in ergodic theory, that is  ~\eqref{hgdhgdhgh769878} and ~\eqref{hgdhgdhgh769878a}, respectively. In the end of the section we prove that Theorem ~\ref{t.returntimeo} implies Theorem ~\ref{Bretthmseriesk=0}, and indicate how a similar argument and Theorem ~\ref{t.returntimeo} imply yet another proof of Bourgain's Return Times theorem. 

\subsection{Transfer for maximal averages (~\eqref{e.returntimeave} $\Rightarrow$ ~\eqref{hgdhgdhgh769878} ).}
 
Fix some $\phi:\Z\to\Z_{+}$ with finite support. For each $a\in \Z$, denote with $C(\phi)(a)$ the best constant which makes the following inequality true for each finitely supported $\psi:\Z\to\Z_{+}$
$$\|\sup_{N}\frac1N\sum_{b=0}^{N-1}\phi(a+b)\psi(c+b)\|_{l^2_c(\Z)}\le C(\phi)(a)\|\psi\|_{l^2(\Z)}.$$
We claim that for each $1<p<\infty$ we have
\begin{equation}
\label{uru2461}
 \|C(\phi)\|_{l^p_a(\Z)}\lesssim\|\phi\|_{l^p(\Z)},
\end{equation}
with bounds independent of $\phi.$ To see this, for each $\phi$ and $\psi$ as above define  $f:\R\to\R$ with $f(x):=\phi([x])$ and  $g:\R\to\R$ with $g(x):=\psi([x])$.  Note that for each $a\le x<a+\frac12$ and each $c\le z<c+\frac12$ we have that 
$$\frac1N\sum_{b=0}^{N-1}\phi(a+b)\psi(c+b)\lesssim \frac1{N}|\int_{0}^{N}f(x+y)g(z+y)dy|,$$ 
uniformly in $x,z,N$. Note also that $\|\phi\|_{l^p(\R)}\sim\|f\|_{L^p(\R)}$, $\|\psi\|_{l^2(\Z)}\sim\|g\|_{L^2(\R)}$.
It turns out that
$$\|\sup_{N}\frac1N\sum_{b=0}^{N-1}\phi(a+b)\psi(c+b)\|_{l^2_c(\Z)}\lesssim  \inf_{a\le x<a+\frac12}\|\sup_{t>0}\frac1t|\int_{0}^{t}f(x+y)g(z+y)dy|\|_{L^2_z(\R)},
$$
and so 
$$C(\phi)(a)\lesssim  \inf_{a\le x<a+\frac12}\sup_{\|g\|_{L^2(\R)}=1}\|\sup_{t>0}\frac1t|\int_{0}^{t}f(x+y)g(z+y)dy|\|_{L^2_z(\R)},$$
which upon using ~\eqref{t.returntimeave} finalizes with
\begin{align*}
\|C(\phi)\|_{l^p_a(\Z)}&\lesssim \|\sup_{\|g\|_{L^2(\R)}=1}\|\sup_{t>0}\frac1t|\int_{0}^{t}f(x+y)g(z+y)dy|\|_{L^2_z(\R)}\|_{L^p_x(\R)}\\&\lesssim \|f\|_{L^p(\R)}\\&\lesssim \|\phi\|_{l^p(\Z)}.
\end{align*}

Consider next two dynamical systems ${\bf X}=(X,\Sigma,\mu, \tau)$ and ${\bf Y}=(Y,\F,\nu,\sigma)$, where the second one is assumed to be ergodic. Fix some large $K>0$, a positive function $f\in L^p(X)$, and the point  $x\in X$. For each $0\le a\le \frac{K}{2}$ and each  $y\in Y$ define $C(a,x,y)$ to be the smallest constant for which
\begin{equation}
\label{bhhhapp90}
\sum_{0\le c\le K/2}(\sup_{N\le K/2}\frac1N\sum_{b=0}^{N-1}f(\tau^{a+b}x)g(\sigma^{c+b}y))^2\le C^2(a,x,y)\sum_{0\le n\le K}g^2(\sigma^ny),
\end{equation}
for each positive function $g\in L^2(Y).$ 
It is an immediate consequence of ~\eqref{uru2461} that
\begin{equation}
\label{bhhhapp901}
\sum_{0\le a\le K/2}\sup_{y\in Y}C(a,x,y)^p\lesssim \sum_{0\le n\le K}f^p(\tau^nx).
\end{equation}
 To see this it suffices to apply ~\eqref{uru2461}
to the functions $\phi,\psi:\Z\to\Z$ defined by
$$\phi(n):=\begin{cases}f(\tau^nx)&:\quad 0\le n \le K \\ \hfill  0&:\quad \text{otherwise}\end{cases}\;\;\;\;  \psi(n):=\begin{cases}g(\sigma^ny)&:\quad 0\le n \le K \\ \hfill  0&:\quad \text{otherwise}\end{cases}$$
By integrating with respect to $y$ in ~\eqref{bhhhapp90} we get for each $x$, $0\le a\le K/2$ and each $g\in L^2(Y)$
$$[\frac K2+1]\int(\sup_{N\le K/2}\frac1N\sum_{b=0}^{N-1}f(\tau^{a+b}x)g(\sigma^{b}y))^2dy\le [K+1]\sup_{y\in Y}C(a,x,y)^2\int g^2(y)dy.$$
Given the universality of $C(a,x,y)$ we get 
$$\sup_{\|g\|_{L^2(Y)}=1\atop{g\ge 0}}\int(\sup_{N\le K/2}\frac1N\sum_{b=0}^{N-1}f(\tau^{a+b}x)g(\sigma^{b}y))^2dy\lesssim \sup_{y\in Y}C(a,x,y)^2.$$
Combining this with ~\eqref{bhhhapp901} we get 
$$\sum_{0\le a\le K/2}\sup_{\|g\|_{L^2(Y)}=1\atop{g\ge 0}}(\int(\sup_{N\le K/2}\frac1N\sum_{b=0}^{N-1}f(\tau^{a+b}x)g(\sigma^{b}y))^2dy)^{p/2}\lesssim \sum_{0\le n\le K}f^p(\tau^nx).$$
Integrate the above with respect to $x$ and divide by $K$ to get 
$$\int\sup_{\|g\|_{L^2(Y)}=1\atop{g\ge 0}}(\int(\sup_{N\le K/2}\frac1N\sum_{b=0}^{N-1}f(\tau^{b}x)g(\sigma^{b}y))^2dy)^{p/2}dx\lesssim \int f^p(x)dx.$$ Finally, let $K\to\infty$ and use the Monotone Convergence Theorem to conclude that
$$\left\|\sup_{\|g\|_{L^2(Y)}=1}\|\sup_{N}\frac1N|\sum_{b=0}^{N-1}f(\tau^{b}x)g(\sigma^{b}y)|\|_{L^2_y(\R)}\right\|_{L^p_x(\R)}\lesssim \|f\|_{L^p(\R)}.$$ Note that this together with Theorem \ref{thm:Conze} immediately imply ~\eqref{hgdhgdhgh769878}.

\subsection{Transfer for maximal truncated series (~\eqref{e.returntimeser} $\Rightarrow$ ~\eqref{hgdhgdhgh769878a} )}

The transfer from ~\eqref{e.returntimeser} to  ~\eqref{hgdhgdhgh769878a} involves similar steps. We start by first observing the following immediate consequence of Corollaries ~\ref{t.returntimeave} and ~\ref{t.returntimeser}:

\begin{corollary} 
\label{t.returntimeser1}  
For each $1<p<\infty$ and each $f\in L^p(\R)$ we have:  
$$
\left\|\sup_{\| g\|_{L^2(\R)}=1}\|\sup_{N}|\sum_{n=1}^{N}\frac{1}{n}(\int_{n}^{n+1}f(x+y)g(z+y)dy-\int_{-n}^{-n+1}f(x+y)g(z+y)dy)|\|_{L^2_z(\R)}\right\|_{L^p_x(\R)}$$
$$\lesssim \|f\|_{L^p(\R)}.$$
\end{corollary} 
\begin{proof}
First note that Corollary ~\ref{t.returntimeser} implies that
$$
\left\|\sup_{\| g\|_{L^2(\R)}=1}\|\sup_{n\ge 1}|\int_{1\le |y|\le n}f(x+y)g(z+y)\frac{dy}{y}|\|_{L^2_z(\R)}\right\|_{L^p_x(\R)}\lesssim \|f\|_{L^p(\R)},\;\;1<p<\infty.
$$
Also,
$$\left|\sum_{n=-N}^N\frac1n\int_{n}^{n+1}f(x+y)g(z+y)dy-\int_{1\le |y|\le n}f(x+y)g(z+y)\frac{dy}{y}\right|$$
$$\lesssim \sup_{t>0}t^{-1}\int_{-t}^t|f(x+y)g(z+y)|dy,$$
and thus  Corollary ~\ref{t.returntimeave} finishes the proof.
\end{proof}

Fix again some $\phi:\Z\to\Z$ with finite support. For each $a\in \Z$, denote with $C(\phi)(a)$ the best constant which makes the following inequality true for each finitely supported $\psi:\Z\to\Z$
$$\left\|\sup_{N}|\frac1N\sideset{}{'}\sum_{b=-N}^{N}\frac{\phi(a+b)\psi(c+b)}{b}|\right\|_{l^2_c(\Z)}\le C(\phi)(a)\|\psi\|_{l^2(\Z)}.$$
We claim that for each $1<p<\infty$ we have
\begin{equation}
\label{uru2461hhhwq}
 \|C(\phi)(a)\|_{l^p_a(\Z)}\lesssim\|\phi\|_{l^p(\Z)},
\end{equation}
with bounds independent of $\phi.$ To see this, for each $\phi$ and $\psi$ as above define  $f,g:\R\to\R$ with 
$$f(x)=\begin{cases}\phi([x])&:\quad [x]+\frac14\le x\le [x]+\frac12 \\ \hfill  0&:\quad \text{otherwise}\end{cases}\;\;\;\; g(x)=\begin{cases}\psi([x])&:\quad [x]+\frac14\le x\le [x]+\frac12 \\ \hfill  0&:\quad \text{otherwise}\end{cases}$$  Note that for each $a\le x<a+\frac1{10}$ and each $c\le z<c+\frac1{10}$ we have that 

$$|\sideset{}{'}\sum_{b=-N}^{N}\frac{\phi(a+b)\psi(c+b)}{b}|\lesssim |\sum_{n=1}^{N}\frac1n(\int_{n}^{n+1}f(x+y)g(z+y)dy-\int_{-n}^{-n+1}f(x+y)g(z+y)dy)|,$$ 
uniformly in $x,z,N$. Note also  that $\|\phi\|_{l^p(\Z)}\sim\|f\|_{L^p(\R)}$, $\|\psi\|_{l^2(\Z)}\sim\|g\|_{L^2(\R)}$. Inequality ~\eqref{uru2461hhhwq} follows as before, by using Corollary ~\ref{t.returntimeser1}. The transfer from $\Z$ to dynamical systems follows exactly the same path as in the case of averages.

\subsection{Transfer for the pointwise convergence (Theorem ~\ref{t.returntimeo} $\Rightarrow$ Theorem ~\ref{Bretthmseriesk=0} )}
\label{sub:helconv}

We first observe that it suffices to prove the convergence of the series in Theorem ~\ref{Bretthmseriesk=0} along a lacunary subsequence. Indeed, fix some $f\in L^{\infty}(X)$ and assume that for each $d_i=2^{1/i}$, $i\in\N$, we know that there exists a universal set $X_i\subseteq X$ with $\mu(X_i)=1$, such that for each second dynamical system  ${\bf Y}=(Y,\F,\nu,\sigma)$, each $g\in L^{\infty}(Y)$  and each $x\in X_i$, the limit
\begin{equation}
\label{lacunarysubseqhh}
\lim_{N\to\infty}\sideset{}{'}\sum_{-d_i^{N}\le n\le d_i^N}\frac{f(\tau^nx)g(\sigma^ny)}{n}
\end{equation}
exists $\nu$-almost everywhere. Let $\tilde{X}$ be a subset of $X$ of full measure such that $|f(\tau^nx)|\le \|f\|_{L^{\infty}(X)}$ for each $x\in\tilde{X}$ and  each $n\in\Z$.
We then use the boundedness of both the weight and the test function to argue that for each $g\in L^{\infty}$, for each $x\in \tilde{X}$, for each $i\in\N$ and for almost every $y\in Y$ we have
$$\limsup_{N,M\to\infty}(\sideset{}{'}\sum_{n=-{N}}^{N}\frac{f(\tau^nx)g(\sigma^ny)}{n}-\sideset{}{'}\sum_{n=-M}^{M}\frac{f(\tau^nx)g(\sigma^ny)}{n})\le$$
$$\limsup_{N,M\to\infty}(\sideset{}{'}\sum_{-d_i^{N}\le n\le d_i^N}\frac{f(\tau^nx)g(\sigma^ny)}{n}-\sideset{}{'}\sum_{-d_i^{M}\le n\le d_i^M}\frac{f(\tau^nx)g(\sigma^ny)}{n})+$$
$$+C\log d_i\|f\|_{L^{\infty}(X)}\|g\|_{L^{\infty}(X)}.$$
Since $i$ can be chosen arbitrarily large,  for each $x\in X_0:=\bigcap_{i\in \N}X_i\cap\tilde{X}$ we get that
$$\limsup_{N,M\to\infty}\left(\sideset{}{'}\sum_{n=-{N}}^{N}\frac{f(\tau^nx)g(\sigma^ny)}{n}-\sideset{}{'}\sum_{n=-M}^{M}\frac{f(\tau^nx)g(\sigma^ny)}{n}\right)=0,$$
for $\nu$- almost every $y$.

It remains to prove that the convergence of the subsequences in ~\eqref{lacunarysubseqhh} follows from Theorem ~\ref{t.returntimeo}. To ease the exposition we will restrict the attention to the case $d=2$ (that is $i=1$). The argument for general $i$ poses no further difficulties.
 Let $K$ be a $C^{\infty}(\R)$ kernel which satisfies the requirements of Theorem ~\ref{t.returntimeo} and in addition satisfies $K(x)=\frac1x$ for $|x|\ge 1$. Introduce the kernels $H_k:\R\to\R,\; k\ge 1$  (these are rough versions of the kernels $\hbox{Dil}_{2^k}^1K$) defined by the formula
$$H_k(x):=\sum_{-2^k\le i\le 2^{k}-1}1_{[i,i+1)}(x)\frac{1}{2^k}K(\frac{i}{2^k})+\sum_{i\in\Z\setminus [-2^k,2^k-1]}1_{[i,i+1)}(x)\frac{1}{i}.$$
Take  an arbitrary sequence $k_1<k_2<\ldots<k_J$ of positive integers. Let $\epsilon(2)$ be such that Theorem ~\ref{t.returntimeo} holds when $p=2$ and $d=2$. As a consequence of this theorem we get that for each $f\in L^2(\R)$ 
\begin{equation*}  
\left\|\sup_{\| g\|_{L^2(\R)}=1}\|\left(\sum_{j=1}^{J-1}\sup_{k_j\le k<k_{j+1}}|\int f(x+y)g(z+y)(H_k({y})-H_{k_{j+1}}({y}))dy|^2\right)^{1/2}\|_{L^2_z(\R)}\right\|_{L^2_x(\R)}\end{equation*}
\begin{equation}  
\label{e.returntime14o}
\lesssim J^{\epsilon(2)}\|f\|_{L^2(\R)},
\end{equation}
with some universal implicit constant (independent of $J$, in particular).
Indeed, note that 
$$|H_k(y)-\hbox{Dil}_{2^k}^1K(y)|\lesssim \begin{cases}&\frac1{2^{2k}},\;\; |y|\le 2^k\\&\frac{1}{y^2},\;\;|y|\ge 2^k\end{cases},$$ 
with the implicit constant independent of $k$.
From the boundedness of the maximal averages (Corollary ~\ref{t.returntimeave}) we deduce that
\begin{equation*}  
\label{e.returntime15o}
\left\|\sup_{\| g\|_{L^2(\R)}=1}\|\left(\sum_{k\ge 1}|\int f(x+y)g(z+y)(H_k({y})-\hbox{Dil}_{2^k}^1K(y))dy|^2\right)^{1/2}\|_{L^2_z(\R)}\right\|_{L^2_x(\R)}\lesssim \|f\|_2.
\end{equation*}
This together with the inequality in Theorem ~\ref{t.returntimeo} and the fact that the terms $k_j$ are positive proves ~\eqref{e.returntime14o}.

The next step consists of transferring  ~\eqref{e.returntime14o} to integers. By following the same lines like in the previous subsections, that is by considering functions  $f,g:\R\to\R$ with 
$$f(x):=\begin{cases}\phi([x])&:\quad [x]+\frac14\le x\le [x]+\frac12 \\ \hfill  0&:\quad \text{otherwise}\end{cases}\;\;\;\;g(x):=\begin{cases}\psi([x])&:\quad [x]+\frac14\le x\le [x]+\frac12 \\ \hfill  0&:\quad \text{otherwise}\end{cases},$$ we get that for each $\phi:\Z\to\Z$ with finite support 
$$  
\left\|\sup_{\|\psi\|_{L^2(\R)}=1}\|\left(\sum_{j=1}^{J-1}\sup_{k_j\le k<k_{j+1}}|\sum_{b\in\Z} \phi(a+b)\psi(c+b)(H_k(b)-H_{k_{j+1}}(b))|^2\right)^{1/2}\|_{l^2_c(\Z)}\right\|_{l^2_a(\Z)}$$
\begin{equation}\label{e.returntime16o} 
\lesssim J^{\epsilon(2)}\|\phi\|_{l^2(\Z)},
\end{equation}
where the first supremum above is taken over all finitely supported functions $\psi:\Z\to\Z$.

For each $k\ge 1$ introduce the kernels $A_k:\Z\to\Z$ and $S_k:\Z\to\Z$ defined by
$$
A_k(i):=\begin{cases}&H_{k}(i),\;\;-2^k\le i\le 2^k\\&0,\;\;\hbox{otherwise}\end{cases},\;\;\;
S_k(i):=\begin{cases}&\frac1i,\;\;-2^k\le i\le 2^k,\;i\not=0\\&0,\;\;\hbox{otherwise}\end{cases},
$$
and note that for each $k<k'$
$$H_k-H_{k'}=O_k-O_{k'}:=(A_k-S_k)-(A_{k'}-S_{k'}).$$ Thus ~\eqref{e.returntime16o} gives 
\begin{equation*}  
\left\|\sup_{\|\psi\|_{L^2(\R)}=1}\|\left(\sum_{j=1}^{J-1}\sup_{k_j\le k<k_{j+1}}|\sum_{b\in\Z} \phi(a+b)\psi(c+b)(O_k(b)-O_{k_{j+1}}(b))|^2\right)^{1/2}\|_{l^2_c(\Z)}\right\|_{l^2_a(\Z)}
\end{equation*}
\begin{equation*}
\label{e.returntime17o}
\lesssim J^{\epsilon(2)}\|\phi\|_{l^2(\Z)},
\end{equation*}
where the first supremum above is taken over all finitely supported functions $\psi:\Z\to\Z$. Standard transfer to a dynamical system ${\bf X}=(X,\Sigma, \mu,\tau)$, as described earlier, leads to 
\begin{equation*}  
\left\|\sup_{(Y,\F,\nu,\sigma)}\sup_{\|g\|_{L^2(Y)}=1}\|\left(\sum_{j=1}^{J-1}\sup_{k_j\le k<k_{j+1}}|\sum_{n\in\Z} f(\tau^nx)g(\sigma^ny)(O_k(n)-O_{k_{j+1}}(n))|^2\right)^{1/2}\|_{L^2_y(Y)}\right\|_{L^2_x(X)}
\end{equation*}
\begin{equation}
\label{e.returntime18o}
\lesssim J^{\epsilon(2)}\|f\|_{L^2(\R)},
\end{equation}
with some universal implicit constant, where the first supremum is taken over all possible dynamical systems ${\bf Y}={(Y,\F,\nu,\sigma)}$. It is then easy to see that this implies the following statement:

({\bf S}): {\em For each function $f\in L^{\infty}(X)$ there is a universal set $X_0\subseteq X$ with $\mu(X_0)=1$, such that for each second dynamical system  ${\bf Y}=(Y,\F,\nu,\sigma)$, each $g\in L^{\infty}(Y)$ and each $x\in X_0$, the weighted averages
$$\sum_{n\in\Z} f(\tau^nx)g(\sigma^ny)O_k(n)$$
converge $\nu$-almost everywhere as $k\to\infty.$}

To see this, assume for contradiction that the above fails for some $f\in L^{\infty}(X)$. It follows that there is a measurable set $X'\subset X$ of positive $\mu$ measure, such that for each $x\in X'$ there is a system ${\bf Y}_x$, a function  $g_x\in L^{2}(Y_x)$ with $\|g_x(y)\|_{L^{2}_y(Y_x)}=1$ and $\alpha(x), \beta(x)>0$ such that 
$$\limsup_{k\to\infty}\sum_{n\in\Z} f(\tau^nx)g_x(\sigma^ny)O_k(n)-\liminf_{k\to\infty}\sum_{n\in\Z} f(\tau^nx)g_x(\sigma^ny)O_k(n)>\alpha(x)$$
for $y$ in a set of $\nu$ measure $\beta(x)$. An elementary measure theoretic argument shows that one can choose a set $X''\subset X'$ of positive $\mu$ measure such that $\alpha(x)>\alpha$ and $\beta(x)> \beta$ for each $x\in X''$, for some $\alpha,\beta>0$. A similar argument shows the existence of set $X'''\subset X''$ of positive $\mu$ measure and of a sequence of positive integers $(k_j)_{j\in\N}$ such that

$$\sup_{k_j\le k<k_{j+1}}\left|\sum_{n\in\Z} f(\tau^nx)g_x(\sigma^ny)O_k(n)-\sum_{n\in\Z} f(\tau^nx)g_x(\sigma^ny)O_{k_{j+1}}(n)\right|>\frac{\alpha}2,$$
for each $j\in \N$ and for each $(x,y)\in X'''\times Y_x'$, where $\nu(Y_x')>\beta$. We immediately get that for each $J$

\begin{equation*}  
\left\|\sup_{(Y,\F,\nu,\sigma)}\sup_{\|g\|_{L^2(Y)}=1}\left\|\left(\sum_{j=1}^{J}\sup_{k_j\le k<k_{j+1}}|\sum_{n\in\Z} f(\tau^nx)g(\sigma^ny)(O_k(n)-O_{k_{j+1}}(n))|^2\right)^{1/2}\right\|_{L^2_y(Y)}\right\|_{L^2_x(X)},
\end{equation*}
\begin{equation*}
\ge\frac{\alpha}{2}(\beta \mu(X'''))^{1/2} J^{1/2},
\end{equation*}
which together with the fact that $\epsilon(2)<\frac12$ contradicts inequality 
~\eqref{e.returntime18o}. The reader is referred to Section ~\ref{sec:approx} for measurability issues regarding the selections of the various sets in the above argument.

The last portion of the argument is devoted to  proving the  statement ({\bf S}) for the weighted averages where $A_k(n)$ replaces $O_k(n)$. This will follow from Bourgain's result for standard averages, Theorem ~\ref{Bretthm}, by means of a common averaging procedure described below. We analyze the two one-sided sums separately, since the mean zero property is no longer crucial in this case. Note that in particular for each $k\ge 1$
$$\sum_{n\ge 1}f(\tau^nx)g(\sigma^ny)A_k(n)=\sum_{n=1}^{2^k}n(A_{k}(n)-A_{k}(n+1))(\frac1n\sum_{i=1}^{n}f(\tau^ix)g(\sigma^iy)).$$
By using Bourgain's result, the fact that 
$$\lim_{k\to\infty}n(A_{k}(n)-A_{k}(n+1))=0$$
for each $n\ge 1$
and the fact that
$$\sup_{k\ge 0}\sum_{n\ge 1}|n(A_{k}(n)-A_{k}(n+1))|<\infty,$$
it follows that we have the return times result for 
$\sum_{n\ge 1} f(\tau^nx)g(\sigma^ny)A_k(n).$
A similar argument works for 
$\sum_{n\le -1} f(\tau^nx)g(\sigma^ny)A_k(n).$ We also trivially have the same result for 
$$f(\tau^0x)g(\sigma^0y)A_k(0)=\frac{K(0)}{2^k}f(x)g(y).$$
This ends the argument.
\endprf

\subsection{Proof of Bourgain's Return Times theorem (Theorem ~\ref{t.returntimeo} $\Rightarrow$ Theorem ~\ref{Bretthm} ).}
\label{reproveBourgain}

 The argument goes as in the previous subsection. The only difference is that this  time we apply  Theorem ~\ref{t.returntimeo} for each $i\in\N$ to a $C^{\infty}(\R)$ kernel $K_i$ which equals $1$ on $[-1,1]$ and 0 on $\T\setminus [-1-\frac1i,1+\frac1i]$, and which also satisfies $\|K_i\|_{L^{\infty}}\le 1$. The error term  caused by the restriction of $K_i$ to ${1\le |x|\le 1+\frac1i}$ is $O(\frac1i)$, and hence can be eliminated by letting $i\to\infty.$
\endprf

\section{Discretization}
\label{sec:5}
We begin this section with  the definition of a (saturated) grid. 
\begin{definition}
\label{d:defgrid}
A set $\G'$ of  intervals each with length in the set $\{2^k:k\in\Z\}$ is called a \emph{saturated grid} if
\begin{enumerate}
\item  for each  $k\in\Z$ there exists $o(k)\in\R$ such that $[o(k)+n2^{k},o(k)+(n+1)2^{k}]\in\G'$ for each $n\in\Z$
\item for every $I,I'\in \G'$ with $I\cap I'\not=\emptyset$ we  have that either $I\subseteq I'$ or $I'\subseteq I$.
\end{enumerate}
If only the second axiom is satisfied then we call $\G'$ a \emph{grid}.
\end{definition}

The endpoints of the intervals in the grid are called \emph{dyadic points}.
We note that if  $\G'$ is a saturated grid, then for each interval $\omega=[a,b]\in\G'$, the subintervals $\omega_1=[a,b]_1:=[a,\frac{a+b}{2}]$ and  $\omega_2=[a,b]_2:=[\frac{a+b}{2},b]$, called the sons of $\omega$ are also in $\G'$. We define the \emph{descendants} of $\omega$ as the collection of all element of $\G'$ which contains its sons, the sons of its sons and so on.
In general, the intervals on the frequency axis will be referred to by the letter $\omega$ while those on the time axis by the letter $I$.

The \emph{standard saturated grid} $\O$ is defined by
$$\O:=\{\;[2^il,2^i(l+1)]:\;i,l\in\Z\}.$$
We will be interested in the following types of grids on the frequency axis: for each odd integer $N\ge 3$, $0\le j\le N-2$ and $0\le L\le N-1$ the collection
$$\G_{N,j,L}:=\left\{\left[2^{i}\left(l+\frac{L}{N}\right),2^{i}\left(l+\frac{L}{N}+1\right)\right]\;:i\equiv j\pmod {N-1},\;l\in\Z\right\}$$
is a  grid, as it easily follows from the fact that $2^{N-1}\equiv 1\pmod N.$ It is not in general a saturated grid, since the  first requirement in Definition ~\ref{d:defgrid} is only satisfied for $k\equiv j\pmod {N-1}$. However, one can easily  turn $\G_{N,j,L}$ into a saturated grid denoted by $\G_{N,j,L}'$ by adding all the descendants of the intervals already in the grid. 
Another interesting observation concerns the fact that for each fixed $N$ the grids $\G_{N,j,L}$ are pairwise disjoint, for $0\le j\le N-2$ and $0\le L\le N-1$. 

Fix now a kernel $K$ as in Theorem ~\ref{t.returntime}.
For each $f\in L^{\infty}(\R)$ with finite support and each $x$  define the operator
$$T_{f,x,K}g(z):=\sup_{k}\frac{1}{2^k}\left|\int f(x+y)g(z+y)\operatorname{Dil}_{2^k}^{1}K(y)dy\right|.$$
Note that we have to prove
$$\left\|\|T_{f,x,K}\|_{L^2_z(\R) \to L^2_z(\R)}\right\|_{L^p_x(\R)}\lesssim \|f\|_{L^p(\R)}.$$
Choose  $\eta:\R\to\R$ such that $\widehat{\eta}$ is a $C^{\infty}(\R\setminus\{0\})$ function which equals $\lim_{\xi\to 0^{+}}\widehat{K}(\xi)$ on $\left(0,\frac18\right]$, $\lim_{\xi\to 0^{-}}\widehat{K}(\xi)$ on $\left[-\frac18,0\right)$ and $0$ outside  $\left[-\frac38,\frac38\right]$. The two limits exist due to the fact that $|\frac{d}{d\xi}\widehat{K}(\xi)|\lesssim 1$ for $\xi\not=0$.
 It suffices to prove
\begin{equation}
\label{en.1}
\left\|\|T_{f,x,\eta}\|_{L^2_z(\R) \to L^2_z(\R)}\right\|_{L^p_x(\R)}\lesssim \|f\|_{L^p(\R)}
\end{equation}
\begin{equation}
\label{enjne.1}
\left\|\|T_{f,x,K-\eta}\|_{L^2_z(\R) \to L^2_z(\R)}\right\|_{L^p_x(\R)}\lesssim \|f\|_{L^p(\R)}.
\end{equation}
The proofs for the above inequalities will follow from a more general principle, as explained below. The crucial property of the multiplier $\widehat{K-\eta}$ that will be used later is the following 
\begin{equation}
\label{eq:24gbnu}
|\frac{d^n}{d{\xi}^n}\widehat{K-\eta}(\xi)|\lesssim \frac{1}{|\xi|^n}\min\{|\xi|,\frac{1}{|\xi|}\},\;\;n\ge 0.
\end{equation}
Note that the additional inequality $|\widehat{K-\eta}(\xi)|\lesssim |\xi|$ for $\xi\not=0$ is  a consequence of the fact that $|\frac{d}{d\xi}\widehat{K}(\xi)|\lesssim 1$ for $\xi\not=0$.
Write 
\begin{equation}
\label{eq:24}
\widehat{K-\eta}(\xi)=\sum_{j=-\infty}^{\infty}\widehat{K-\eta}(\xi)q(\frac{\xi}{2^j}),
\end{equation}
 where $q$ is some Schwartz function supported in the annulus $\frac18<|\xi|<\frac38$ such that 
$$\sum_{j\in \Z}q(\frac{\xi}{2^j})=1,\;\;\xi\not=0.$$  As a consequence of ~\eqref{eq:24gbnu}, each function $g_j=\widehat{K-\eta}(\xi)q(\frac{\xi}{2^j})$ will satisfy 
$$\|\frac{d^n}{d\xi^n}g_j(\xi)\|_{L^\infty_\xi(\R)}\lesssim \frac{2^{-|j|}}{2^{jn}},\;\;\xi\not=0,$$ 
for all $n\ge 0$, uniformly in $j\in\Z$. It follows that that each function $\operatorname{Dil}_{2^{-j}}^{\infty}g_j$ satisfies 
$$\|\frac{d^n}{d\xi^n}\operatorname{Dil}_{2^{-j}}^{\infty}g_j(\xi)\|_{L^\infty_\xi(\R)}\lesssim 2^{-|j|},\;\;\xi\not=0,$$
for all $n\ge 0$, uniformly in $j\in\Z$. Moreover, it is  supported in the  annulus $\frac18<|\xi|<\frac38$. Since the operators $T_{f,x,\check{g_j}}$ and $T_{f,x,\operatorname{Dil}_{2^{j}}^{1}\check{g_j}}$ coincide, inequality ~\eqref{enjne.1} will immediately follow if we prove that
\begin{equation}
\label{en.2}
\left\|\|T_{f,x,\check{\psi}}\|_{L^2_z(\R) \to L^2_z(\R)}\right\|_{L^p_x(\R)}\lesssim \|f\|_{L^p(\R)},
\end{equation}
uniformly in all Schwartz functions $\psi$  supported as above and satisfying 
\begin{equation}
\label{gdhgfhgtr}
\|\frac{d^n}{d\xi^n}\psi(\xi)\|_{L^{\infty}_\xi(\R)}\lesssim 1
\end{equation}
for all $n\ge 0$. 
 
From now on $\psi$ will be either a  function as above or the function $\widehat{\eta}$. We next focus on proving ~\eqref{en.2}. By a dilation argument we can assume in addition that  $\psi(\xi)-\psi(2\xi)$ is supported in the annulus $\frac{1}{16}\le |\xi|\le \frac3{8}$. Triangle's inequality further  allows us to assume that the support is inside $\left[\frac{1}{16}, \frac3{8}\times\right]$. Note that 
$$\psi(2^k\xi)=\sum_{i\ge k}\psi_i(\xi),
$$ 
with $\psi_i(\xi):=\psi(2^{i}\xi)-\psi(2^{i+1}\xi)$  supported in  $[\frac{1}{16}2^{-i}, \frac3{8}\times2^{-i}]$.
For each $f$ as above  and for each $x$ we have 
\begin{equation}
\|T_{f,x,\psi}\|_{L^2_z(\R) \to L^2_z(\R)}=\sup_{\|g\|_{2}=1}\left\|\sup_{k\in\Z}|\sum_{i\ge k}\int_{\R}f(x+y)g(z+y)\check{\psi_i}(y)dy|\right\|_{L^2_z(\R)}.
\end{equation}

Pick a Schwartz function $\varphi$ such that $\widehat{\varphi}$ is supported in $[0,\frac{2}{41}]$ and satisfies the following property for every $\xi\in\R$:
$$\sum_{l\in \Z}\left|\widehat{\varphi}\left(\xi-\frac{l}{41}\right)\right|^2=1.
$$
For each scale $i$ use the following expansion for  $f$, valid in every $L^p(\R) ,1<p<\infty$ norm

$$f=\sum_{m,l\in \Z}\langle f,\varphi_{i,m,\frac{l}{41}}\rangle\varphi_{i,m,\frac{l}{41}},$$
where $\varphi_{i,m,l}$ is the modulated wave packet (see \cite{LTBilH} for a similar expansion)
$$\varphi_{i,m,l}(x):=2^{-\frac{i}{2}}\varphi(2^{-i}x-m)e^{2\pi i2^{-i}xl}.$$
Now
\begin{align*}
\sup_{k\in\Z}|\sum_{i\ge n}\int_{\R}f(x+y)g(z+y)\check{\psi_i}(y)dy|
&=\sup_{k\in\Z}|\sum_{i\ge n\atop{m,l\in \Z}}\int_{\R}\langle f,\varphi_{i,m,\frac{l}{41}}\rangle\varphi_{i,m,\frac{l}{41}}(x+y)g(z+y)\check{\psi_{i}}(y)dy|\\&=\sup_{k\in\Z}\left|\tilde{F}_{k,x}\ast g(z)\right|,
\end{align*}
where $\tilde{F}$ denotes the reflection $\tilde{F}(y):=F(-y)$ and 
$$F_{k,x}(y):=\sum_{i\ge k}\sum_{m,l\in \Z}\langle f,\varphi_{i,m,\frac{l}{41}}\rangle\varphi_{i,m,\frac{l}{41}}(x+y)\check{\psi_{i}}(y).$$ 
With this notation, the inequality ~\eqref{e.returntime} follows from
\begin{equation}
\label{en.4}
\left\|\|(\F_y(F_{k,x}))_{k\in\Z}\|_{M_{2,\theta}^*(\R)}\right\|_{L^p_x(\R)}\lesssim\|f\|_{L^p(\R)}.
\end{equation}
Here we use $M_{2,\theta}^*(\R)$ to denote the maximal multiplier norm $M_2^*$ in the $\theta$ variable.
Note that the Fourier transform of $F_{k,x}(y)$ in the $y$ variable is 
$$\F_y(F_{k,x})(\theta)=\sum_{i\ge k}\sum_{m,l\in \Z}\langle f,\varphi_{i,m,\frac{l}{41}}\rangle\int_{\R}\psi_0(2^i(\theta-\xi))\widehat{\varphi}_{i,m,\frac{l}{41}}(\xi)e^{2\pi i\xi x}d\xi.$$ 
Define $\phi_{i,m,\frac{l}{41}}(x,\theta):=\int_{\R}\psi_0(2^i(\theta-\xi))\widehat{\varphi}_{i,m,\frac{l}{41}}(\xi)e^{2\pi i\xi x}d\xi$ and note that
\begin{equation}
\label{xbfgr2}
\phi_{i,m,\frac{l}{41}}(x,\theta)=2^{-\frac{i}{2}}\phi_{0,0,0}(2^i\theta-\frac{l}{41},x2^{-i}-m)e^{2\pi i\frac{l}{41}(x2^{-i}-m)}.
\end{equation}
The function $\phi_{0,0,0}$ is in $C^{\infty}(\R\times \R)$, and as a consequence of ~\eqref{gdhgfhgtr} satisfies the following
\begin{equation}
\label{xbfgr1}
\|\frac{\partial^n}{\partial\theta^n}\frac{\partial^m}{\partial x^m}\phi_{0,0,0}(x,\theta)\|_{L^{\infty}(\theta)}\lesssim \frac{1}{(1+|x|)^M},\;\;\forall n,m,M\ge 0.
\end{equation}
The function $\phi_{i,m,\frac{l}{41}}(x,\theta)$ and its  $x$ Fourier transform $\F_x(\phi_{i,m,\frac{l}{41}}(x,\theta))(\xi)=\psi(2^i(\theta-\xi))\widehat{\varphi}_{i,m,\frac{l}{41}}(\xi)$ are localized as follows:
\begin{align}
\nonumber
\supp_{\theta}(\phi_{i,m,\frac{l}{41}}(x,\theta))&\subseteq \left[2^{-i}\frac{l}{41},2^{-i}\frac{l+2}{41}\right]+\left[2^{-i}\frac{1}{16},2^{-i}\frac{3}{8}\right]\\&=
\label{e.mmccffttss}
\left[2^{-i}\left(\frac{l}{41}+\frac{1}{16}\right),2^{-i}\left(\frac{l+2}{41}+\frac{3}{8}\right)\right],\;\hbox{for each\;}x
\end{align}

\begin{equation}
\label{e.mmccffttss1}
\supp_{\xi}(\F_x(\phi_{i,m,\frac{l}{41}}(x,\theta))(\xi))\subseteq \left[2^{-i}\frac{l}{41},2^{-i}\frac{l+2}{41}\right],\;\hbox{for each\;}\theta.
\end{equation}
The crucial property of these supports is that 
$$\left[2^{-i}\left(\frac{l}{41}+\frac{1}{16}\right),2^{-i}\left(\frac{l+2}{41}+\frac{3}{8}\right)\right]\subseteq \left[2^{-i}\frac{l-18}{41},2^{-i}\left(\frac{l-18}{41}+1\right)\right]_2,$$ and
$$\left[2^{-i}\frac{l}{41},2^{-i}\frac{l+2}{41}\right]\subseteq \left[2^{-i}\frac{l-18}{41},2^{-i}\left(\frac{l-18}{41}+1\right)\right]_1,$$
where $\omega_{i,l}:=\left[2^{-i}\frac{l-18}{41},2^{-i}\left(\frac{l-18}{41}+1\right)\right]$ is in some (unique) grid $\G_{41,j,L}$.

To each $m,i,l\in\Z$ we associate the tile $s=[2^{i}m,2^i(m+1)]\times \omega_{i,l}$ and use the notation $\varphi_s:=\varphi_{i,m,\frac{l}{41}}$, $\phi_s:=\phi_{i,m,\frac{l}{41}}$. As a consequence of ~\eqref{xbfgr2}, ~\eqref{xbfgr1}, ~\eqref{e.mmccffttss} and ~\eqref{e.mmccffttss1}, the localization and decay of $\phi_s$ can now be summarized as follows:
\begin{equation}
\label{xbfgr1456wd}
\supp_{\theta}(\phi_{s}(x,\theta))\subseteq \omega_{s,2}\;\hbox{for each\;}x
\end{equation}
\begin{equation}
\label{xbfgr1456jr}
\supp_{\xi}(\F_x(\phi_{s}(x,\theta))(\xi))\subseteq \omega_{s,1}\;\hbox{for each\;}\theta
\end{equation}
\begin{equation}
\label{xbfgr1456bb}
\sup_{c\in\omega_s}\left\|\frac{\partial^n}{\partial\theta^n}\frac{\partial^m}{\partial x^m}\left[\phi_{s}(x,\theta)e^{-2\pi icx}\right]\right\|_{L^{\infty}_\theta(\R)}\lesssim |I_s|^{(n-m-1/2)}\chi_{I_s}^M(x),\;\; \forall n,m,M\ge 0,
\end{equation}
uniformly in $s$.
We also note that 
\begin{equation}
\label{xbfgr1456wd731q}
\supp (\widehat{\varphi_{s}})\subseteq \omega_{s,1}
\end{equation}
and 
\begin{equation}
\label{xbfgr14561o}
\sup_{c\in\omega_s}\left|\frac{\partial^n}{\partial x^n}\left[\varphi_{s}(x)e^{-2\pi icx}\right]\right|\lesssim {|I_s|}^{-n-\frac12}\chi_{I_s}^M(x),
\end{equation}
uniformly in $s$.

For each $j,L$ as above define a collection of tiles 
$$\S_{j,L}:=\{[2^{i}m,2^i(m+1)]\times \omega:\;\omega\in \G_{41,j,L},\;m\in\Z,\;2^{i}|\omega|=1\}$$
and note that ~\eqref{en.4} is equivalent to 
\begin{equation*}
\|\|(\sum_{j,L}\sum_{s\in \S_{j,L}\atop{|I_s|<2^k}}\langle f,\varphi_s\rangle\phi_s(x,\theta))_{k\in\Z}\|_{M_{2,\theta}^*(\R)}\|_{L^p_x(\R)}\lesssim\|f\|_p.
\end{equation*}
In the above we changed the restriction $|I_s|\ge 2^k$ into the more suitable for later purposes $|I_s|< 2^k$. Note that they are equivalent. Theorem ~\ref{t.returntime} will be a consequence of the following more general result:

\begin{theorem}
\label{ceamaiceahh}
Let $\G'$ be a saturated grid and let $\S$ be some arbitrary finite subcollection of the set of all tiles 
$$\S_{\univ}:=\{[2^{i}m,2^i(m+1)]\times \omega:\;\omega\in \G',\;m,i\in\Z,\;2^{i}|\omega|=1\}.$$
Consider also two  collections $\{\phi_s,s\in\S\}$ and $\{\varphi_s,s\in\S\}$ of Schwartz functions. The functions $\phi_s:\R^2\to\R$ satisfy ~\eqref{xbfgr1456wd}, ~\eqref{xbfgr1456jr} and ~\eqref{xbfgr1456bb}, uniformly in $s$. The functions $\varphi_s:\R\to\R$ satisfy ~\eqref{xbfgr1456wd731q} and ~\eqref{xbfgr14561o}, uniformly in $s$.

Then the following inequality holds for each $f\in L^p(\R)$, $1<p<\infty$
\begin{equation}
\label{en.511}
\|\|(\sum_{s\in \S\atop{|I_s|<2^k}}\langle f,\varphi_s\rangle\phi_s(x,\theta))_{k\in\Z}\|_{M_{2,\theta}^*(\R)}\|_{L^p_x(\R)}\lesssim\|f\|_p,
\end{equation} 
with the implicit constant depending only on $p$ and on the implicit constants in ~\eqref{xbfgr1456wd} and ~\eqref{xbfgr1456bb} (in particular it is independent of the choice of the grid).
\end{theorem}

The same discretization techniques immediately show that Theorem ~\ref{t.returntimeo1} will follow from the following:
\begin{theorem}
\label{ceamaiceahho1}
Assume we are in the settings from the above theorem. For each $1<p<\infty$ there is $0<\epsilon(p)<\frac12$ such that  for each finite sequence of integers $u_1<u_2<\ldots<u_J$ 
$$\left\|\sup_{\|g\|_{L^2(\R)}=1}(\sum_{j=1}^{J-1}\|\sup_{\;u_j\le
k<u_{j+1}}|\F^{-1}_{\theta}\{\sum_{s\in \S\atop{2^{u_j}\le
|I_s|<2^k}}\langle
f,\varphi_s\rangle\phi_s(x,\theta)\widehat{g}(\theta)\}(z)|\|_{L^2_z(\R)}^2)^{1/2}\right\|_{L^p_x(\R)}$$
\begin{equation}
\label{en.5o}
\lesssim
J^{1/2-\epsilon(p)}\|f\|_{L^p(\R)},
\end{equation}
with the implicit constant depending only on $p$ and on the implicit constants in ~\eqref{xbfgr1456wd} and ~\eqref{xbfgr1456bb}.
\end{theorem}

By a very similar argument, Theorem ~\ref{t.returntimeo2} will follow from the following:

\begin{theorem}
\label{ceamaiceahho2}
Assume we are in the settings from  Theorem  ~\ref{ceamaiceahh}. For each $1<p<\infty$ the following inequality holds
\begin{equation}
\label{en.5o1}
\|\|(\sum_{k\in\Z}|\sum_{s\in \S\atop{|I_s|=2^k}}\langle
f,\varphi_s\rangle\phi_s(x,\theta)|^2)^{1/2}\|_{M_{2,\theta}(\R)}\|_{L^p_x(\R)}\lesssim \|f\|_{L^p(\R)},
\end{equation}
with the implicit constant depending only on $p$  and on the implicit constants in ~\eqref{xbfgr1456wd} and ~\eqref{xbfgr1456bb}.
\end{theorem}

In the remaining sections we will prove Theorems ~\ref{ceamaiceahh}, ~\ref{ceamaiceahho1} and ~\ref{ceamaiceahho2}. From now on, by a dyadic frequency interval we will understand any interval of the saturated grid $\G',$ while a dyadic time interval will continue to refer to an interval in the standard dyadic grid.

\section{Trees}
\label{sec:trees}

We now recall some standard terminology concerning trees of tiles. (see \cite{LTBilH} and \cite{MTT1} for more details)

\begin{definition}[Tile order]
For two tiles $s$ and $s'$ we write $s\le s'$ if $I_s\subseteq I_{s'}$ and $\omega_{s'}\subseteq\omega_s$.
\end{definition}

\begin{definition}[Trees]
A \emph{tree} with \emph{top} $T\in \S_{\univ}$ is a set of tiles $\T\subseteq\S$ such that
$s\le T$ for each $s\in \T$.  For $i=1,2$, we say that an \emph{$i$-tree} is a tree $\T$ such that $\omega_T\subseteq\omega_{s,i}$ for each $s\in \T\setminus T$, where the intervals $\omega_{s,1}$ and $\omega_{s,2}$ are the left and right halves of $\omega_s.$
\end{definition}

We will also encounter  a more general instance of a tree called ``quasitree''.
 
\begin{definition}
A \emph{quasitree} with \emph{top} $(I_{\T},\xi_{\T})$, where $I_{\T}$ is an arbitrary (not necessarily dyadic) interval and $\xi_{\T}\in\R$ is a (not necessarily dyadic\footnote{In fact, $\xi_{\T}$ may always be taken to be non-dyadic} point, is a set of tiles $\T\subseteq\S$ such that
$I_s\subset I_{\T}$ and $\xi_{\T}\in\omega_s$ for each $s\in \T$.  If $i=1,2$, an \emph{$i$-quasitree} is a quasitree $\T$ such that $\xi_{\T}\in\omega_{s,i}$ for each $s\in \T$, where the intervals $\omega_{s,1}$ and $\omega_{s,2}$ are the left and right halves of $\omega_s.$
\end{definition}

\begin{remark}
Note that each tree $\T$ with top $T$ is a also a quasitree with top $(I,\xi)$, for each interval $I_T\subseteq I$ and each  $\xi\in \omega_T$ which is not a dyadic point. We will adopt the convention that $I_{\T}=I_T$ and $\xi_{\T}\in \omega_{T,1}$, without any further specification on $\xi_{\T}$.
\end{remark}

The standard decomposition of a quasitree $\T$ with top $(I_{\T},\xi_{\T})$ is the splitting of $\T$ into the $1$-quasitree
$$\T^{(1)}:=\{s\in\T:\xi_{\T}\in \omega_{s,1}\}$$
and the $2$-tree
$$\T^{(2)}:=\{s\in\T:\xi_{\T}\in \omega_{s,2}\}.$$
Note that if $\T$ is a tree with top $T$ then this decomposition does not depend on the choice of $\xi_{\T}\in \omega_{T,1}$, and moreover, if $T\in\T$ then $T\in\T^{(1)}$.

\begin{definition}
Fix some $f:\R\to\R$.
For a  finite subset of tiles $\S'\subseteq\S$ define its \emph{size} relative to $f$ as 
$$\size(\S'):=\sup_{\T}\left(\frac{1}{|I_T|}\sum_{s\in \T}|\langle f, \varphi_{s}\rangle|^2\right)^{\frac{1}{2}}$$ where the supremum is taken over all the $2$-trees $\T \subset \S'$. 
\end{definition}

We recall two important results regarding the size.
\begin{proposition}
\label{g11a00uugg55}
For each $1<t<\infty$, each 2-tree $\T$ with top $T$ and  each $f\in L^t(\R)$ we have 
$$\left(\frac1{|I_T|}{\sum_{s\in\T}|\langle f,\varphi_s\rangle|^2}\right)^{1/2}\lesssim \inf_{x\in I_T}M_t f(x).$$
\end{proposition}
\begin{proof}
See for example Lemma 1.8.1 in \cite{Th}.
\end{proof} 

The following Bessel type inequality from \cite{LTCar} will be useful in organizing  collections of tiles into trees. 
\begin{proposition}
\label{Besselsineq}
Let $\S'\subseteq\S$ be a collection  of tiles and define $\Delta:=[-\log_2(\size(\S'))]$, where the size is understood with respect to some function $f\in L^2(\R)$. Then $\S'$ can be written as a disjoint union $\S'=\bigcup_{n\ge \Delta}\P_n,$ where $\size(\P_n)\le 2^{-n}$ and each $\P_n$ consists of a family $\F_{\P_n}$ of pairwise disjoint trees  satisfying 
\begin{equation}
\label{e:intp}
\sum_{\T\in\F_{\P_n}}|I_T|\lesssim 2^{2n}\|f\|_2^2,
\end{equation} 
with bounds independent of $\S'$, $n$ and $f$. 
\end{proposition}

In the following we will use the notation for the counting function associated with a collection $\F$ of quasitrees
$$N_{\F}(x):=\sum_{\T\in\F}1_{I_T}(x).$$

Let $\T$ be a  2-quasitree with top $(I_{\T},\xi_{\T})$. The following decomposition will be useful in the future.
For each $s\in\T$ and scale $l\ge 0$ we split  $\phi_s(x,\theta)$ as $$\phi_s(x,\theta)=\tilde{\phi}_{s,\T}^{(l)}(x,\theta)+\phi_{s,\T}^{(l)}(x,\theta).$$ For convenience, we set $\phi_{s,\T}^{(0)}:=\phi_{s}$ for each $s\in\T$. For $l\ge 1$ we define the first  piece to be localized in time:
 $$\operatorname {supp}\tilde{\phi}_{s,\T}^{(l)}(\cdot,\theta)\subseteq 2^{l-1}I_s,\;\;\hbox{for each\;}\theta\in\R.$$ 
For the second piece we need some degree of frequency localization, but obviously full localization as in the case of $\phi_s$ is impossible. We will content ourselves with preserving the mean zero property with respect to the top of the quasitree. The advantage of $\phi_{s,\T}^{(l)}$ over $\phi_s$ is that it gains extra decay in $x$. More precisely, we have for each $s\in \T$ and each $M\ge 0$

\begin{equation}
\label{e.m_s-supportati}
\phi_{s,\T}^{(l)}(x,\theta)e^{-2\pi i\xi_{\T}x}\;\hbox{ has mean zero},\;\;\theta\in\R,
\end{equation}

\begin{equation}
\label{e.m_s-size}
\phi_{s,\T}^{(l)}(x,\theta)e^{-2\pi i\xi_{\T}x}\hbox{\;is\;} c(M)2^{-Ml}-\hbox{adapted to}\; I_s,\;\;\hbox{for some constant\;}c(M), \;\;\theta\in\R,
\end{equation}

\begin{equation} 
\label{e.m_s-supportgrr}
\operatorname {supp}\phi_{s,\T}^{(l)}(x,\cdot)\subset \omega_{s,2},\hbox{\;for each\;}x\in\R, 
\end{equation} 

 \begin{equation}
 \label{e.m_s-smoothgrr} 
| \tfrac{ d}{d \theta} \phi_{s,\T}^{(l)}(x,\theta) |\lesssim 2^{-Ml} |I_s|^{\frac12} \chi_{I_s}^{M}(x),\;\;\hbox{uniformly in \;}x,\theta\in\R.
\end{equation}

We achieve this decomposition by first choosing a smooth function $\eta$ such that $\supp(\eta)\subset [-1/2,1/2]$ and $\eta=1$ on $[-1/4,1/4]$. We then define 
$$\tilde{\phi}_{s,\T}^{(l)}(\theta;x):=\phi_{s}(\theta;x)\eta\text{Dil}_{2^{l}I_s}^{\infty}\eta(x)-\frac{e^{2\pi i\xi_{\T}x}\text{Dil}_{2^{l}I_s}^{\infty}\eta(x)}{\int_{\mathbb R}\text{Dil}_{2^{l}I_s}^{\infty}\eta(x)dx}\int_{\mathbb R}\phi_{s}(\theta;x)e^{-2\pi i\xi_{\T}x}\text{Dil}_{2^{l}I_s}^{\infty}\eta(x)dx$$
and
$$\phi_{s,\T}^{(l)}(\theta;x):=\frac{e^{2\pi i\xi_{\T}x}\text{Dil}_{2^{l}I_s}^{\infty}\eta(x)}{\int_{\mathbb R}\text{Dil}_{2^{l}I_s}^{\infty}\eta(x)dx}\int_{\mathbb R}\phi_{s}(\theta;x)e^{-2\pi i\xi_{\T}x}\text{Dil}_{2^{l}I_s}^{\infty}\eta(x)dx+\phi_{s}(\theta;x)(1-\text{Dil}_{2^{l}I_s}^{\infty}\eta(x)).$$
Properties ~\eqref{e.m_s-supportati} through ~\eqref{e.m_s-smoothgrr} are now easy consequences of ~\eqref{xbfgr1456wd}, ~\eqref{xbfgr1456jr} and ~\eqref{xbfgr1456bb}.

In the next two sections we prove some general results of independent interest, which will be used later in the main argument.

\section{A weighted Bourgain's Lemma}
\label{sec:7}
For each $1\le r<\infty$ and each sequence $(x_k)_{k\in \Z}$ in a Hilbert space $\H$, define the \emph{$r$-variational norm} of  $(x_k)_{k \in \Z}$  to be
$$\|x_k\|_{V^r_k(\Z)}:=\sup_{k}\|x_k\|_{\H}+ \|x_k\|_{\tilde V^r_k(\Z)}$$
where $\tilde V^r_k(\Z)$ is the homogeneous $r$-variational seminorm
$$ \|x_k\|_{\tilde V^r_k(\Z)}:= \sup_{M,k_0<k_1< \ldots <k_{M}}(\sum_{m=1}^M\|x_{k_m}-x_{k_{m-1}}\|_{\H}^r)^{1/r}.$$ 
We also write $V^r_k(L)$ and $\tilde V^r_k(L)$ for the variational norm of a sequence $(x_k)_{1 \leq k \leq L}$ of $L$ elements.
Define also the \emph{oscillation norm} $\|\cdot\|_{O_{{\bf U}}}$ of a sequence $(x_k)$  with respect to the sequence of integers ${\bf U}=(u_j)_{j=1}^{J}$ to be
$$\|x_k\|_{O_{{\bf U}}}=(\sum_{j=1}^{J-1}\sup_{u_j\le k<u_{j+1}}\|x_k-x_{u_{j}}\|_{\H}^2)^{1/2}.$$
For each $r>2$ define also the oscillation-variational norm 
 
\begin{equation}
\label{defoscvar}
 \|x_k\|_{O_{{\bf U}} \cap V^r_k(\Z)}:=\|x_k\|_{O_{{\bf U}}}+\|x_k\|_{V^r_k(\Z)}.
\end{equation}

For future reference we record the following easily verified lemma.

\begin{lemma}[Product estimates]
\label{lem:auxq3}
For each $k$, let $a_k,b_k$ be some complex numbers and let ${\bf U}:=u_1<u_2<\ldots<u_j$ be an arbitrary finite sequence of integers. Then for each $r\ge 1$
$$\|a_kb_k\|_{V^r_k(\Z)}\lesssim \|a_{k}\|_{V^r_k(\Z)}\|b_{k}\|_{V^r_k(\Z)},$$
$$\|a_kb_k\|_{O_{{\bf U}}}\lesssim \|a_{k}\|_{O_{{\bf U}}}\|b_k\|_{l_k^\infty(\Z)}+\|b_{k}\|_{O_{{\bf U}}}\|a_k\|_{l_k^\infty(\Z)}.$$
\end{lemma}

Consider a finite set  $\Lambda=\{\lambda_1,\ldots,\lambda_L\}$ such that each  dyadic frequency  interval\footnote{that is, intervals in the saturated grid $\G'$} of length 1  contains at most one  element of $\Lambda$. For each $k\ge 0$ define $R_k$ to be the collection of the $L$ dyadic frequency intervals of length $2^{-k}$ which contain  an  element from $\Lambda$, and denote by $\omega_{k,l}$ the one that contains $\lambda_l$. Also, for each $k\ge 0$ and each $1\le l\le L$ consider some multipliers $m_{k,l}:\R\to\C$. Define
$$\Delta_kf(x):=\sum_{l}\int_{\omega_{k,l}}m_{k,l}(\xi)\widehat{f}(\xi)e^{2\pi i\xi x}d\xi .$$ The following theorem  is a particular case of the main result of this section, Theorem ~\ref{thm:gh098cfhh}. 
\begin{theorem}
\label{thm:gh098cf}
For each $r>2$ we have the inequality
$$\|\sup_{k\ge 0}|\Delta_kf|\|_{L^2(\R)} \lesssim L^{1/2-1/r}\sup_l\sup_{\|g\|_{L^2(\R)} =1}\left\|\|(m_{k,l}1_{\omega_{k,l}}\widehat{g})\check{\ }(z)\|_{V^r_k(L)} \right\|_{L^2_z(\R)} \|f\|_{L^2(\R)},$$
with the implicit constant depending only on $r$.
\end{theorem}

The gain in this theorem is over the exponent of $L$, given the fact that the triangle inequality trivially implies the result with $L$ replacing $L^{1/2-1/r}$. The remaining part of the bound is an amorphous quantity, its later estimate will depend on the multipliers in question.

The proof presented below of the above theorem  relies on a couple of lemmas and is heavily inspired by an argument of Bourgain for a particular case (Corollary ~\ref{labellhgt64}), see \cite{Bo1}.  We will denote by $l^2(L)$ the Hilbert space of all the finite sequences $(c_l)_{1\le l\le L}.$ The following result is classical.
\begin{lemma}
\label{lem:auxq1}
For each finite set $A\subseteq l^2(L)$ with cardinality $\sharp A$, we have
$$\|\sup_{a\in A}|\sum_{l=1}^{L}a_{l}e^{2\pi i \lambda_l y}|\|_{L^2_y([0,1))}\lesssim \min(\sqrt{L},\sqrt{\sharp A})\sup_{a\in A}\|a\|_{l^2(L)}.$$
\end{lemma}

\begin{proof}(Sketch) To obtain the bound involving $\sqrt{L}$, take absolute values everywhere and use Cauchy-Schwarz.  To obtain the bound
involving $\sqrt{\sharp A}$, estimate the left-hand side by the square function
$$(\sum_{a \in A} \| \sum_{l=1}^{L} a_l e^{2\pi i \lambda_l y} \chi(y) \|_{L^2_y(\R)}^2)^{1/2}$$
where $\chi$ is a bump function supported on $[-1,2]$ that equals one on $[0,1]$, and then use Plancherel's theorem.
\end{proof}

We  use this lemma to prove the following.
\begin{lemma}
\label{lem:auxq2}
For each  set $C=\{c_k\}\subseteq l^2(L)$ and each $r>2$ we have
$$\|\sup_{k}|\sum_{l=1}^{L}c_{k,l}e^{2\pi i \lambda_l y}|\|_{L^2_y([0,1))}\lesssim L^{1/2-1/r}\|c_k\|_{V^r_k(L)},$$
with the implicit constant depending only on $r$.
\end{lemma}
\begin{proof}
The proof of this lemma relies on a standard metric entropy approach.
It suffices to prove it in the case  $C$ is finite and then to invoke the Monotone Convergence Theorem. For each $\lambda>0$ denote by $M_{\lambda}$ the minimum number of balls with radius $\lambda$ and centered at elements of $C$, whose union covers $C.$ It is an easy exercise to prove that
\begin{equation}
\label{2uanoap}
\sup_{\lambda>0}\lambda M_{\lambda}^{1/r}\lesssim \|c_k\|_{V^r_k(L)},
\end{equation}
with the implicit constant depending only on $r$.
Let $c^{*}$ be an arbitrary element of $C$. For each $n\ge -\log_2(\diam (C))$, let $C_n$ be a  collection of elements of $(C-C)\cup\{{\bf 0}\}$ such that
\begin{align*}
\|c\|_{l^2(L)}&\le 2^{-n+2}\;\;\hbox{for each}\;c\in C_n,\\
\sharp C_n&\le M_{2^{-n}}+1
\end{align*}
and each $c\in C$ can be written as 
\begin{equation}
\label{sumrepresent}
c=c^{*}+\sum_{n\ge -\log_2(\diam (C))}c_n\;\;\hbox{with}\;c_n\in C_n.
\end{equation}
Here is how $C_n$ is constructed. For each $n\ge -\log_2(\diam (C))$ define $B_{n}$ to be a collection of $M_{2^{-n}}$ elements of $C$ such that the balls with centers in $B_{n}$ and radius $2^{-n}$ cover $C$. If $n=[-\log_2(\diam (C))]-1$ define $B_n=\{c^*\}$. For each $n\ge -\log_2(\diam (C))$ and each $c\in B_{n}$, choose an element $c'\in B_{n-1}$ such that the ball centered at $c$ and with radius $2^{-n}$ intersects  the ball centered at $c'$ and with radius $2^{-n+1}$. Define 
$$C_n:=\{c-c':c\in B_n\}\cup\{{\bf 0}\}.$$
Since $C$ is finite, for each $c\in C$ there is $n$ such that $c\in B_n$. To verify the representation ~\eqref{sumrepresent} for an arbitrary $c\in C$,  denote as above by $c'$ the element from $B_{n-1}$ associated with $c$, by $c''$ the element from $B_{n-2}$ associated with $c'$ and so on, and note that this sequence will eventually terminate with $c^*$. Hence we can write
$$c=(c-c')+(c'-c'')+\ldots+c^*.$$ Note also that by construction, each element of $C_n$ has norm at most $2^{-n+2}$.

If for each $c\in l^2(L)$ we define 
$$X_c(y)=\sum_{l=1}^{L}c_l e^{2\pi i \lambda_l y}$$ then we have
$$X_c(y)=X_{c^*}(y)+\sum_{n\ge -\log_2(\diam (C))}X_{c_n}(y)\;\;\hbox{with}\;c_n\in C_n,$$
for each $c\in C.$ This together with  inequality ~\eqref{2uanoap} and Lemma ~\ref{lem:auxq1} further allows us to write
\begin{align*}
\|\sup_{c\in C}|X_c(y)|\|_{L^2_y([0,1))}&\le \|X_{c^*}(y)\|_{L^2_y([0,1))}+\sum_{n\ge -\log_2(\diam (C))}\|\sup_{c\in C_n}|X_{c}(y)|\|_{L^2_y([0,1))}\\&\le \sup_{c\in C}\|c\|_{l^2(L)}+\sum_{n\ge -\log_2(\diam (C))}\min(\sqrt{L},\sqrt{M_{2^{-n}}+1})\sup_{c\in C_n}\|c\|_{l^2(L)}\\&\lesssim \sup_{c\in C}\|c\|_{l^2(L)}+\sum_{n\in\Z}2^{-n}\min(\sqrt{L},\|c_k\|_{V^r_k(L)}^{r/2}2^{nr/2})\\&\lesssim \|c_k\|_{V^r_k(L)}L^{1/2-1/r}.
\end{align*}
\end{proof}

\begin{lemma}
\label{lem:auxq4}
Let $(\omega_k)_{k\in\Z}$ be a sequence of nested  dyadic frequency intervals with  $|\omega_k|=2^{-k}$ and let also ${\bf U}:=u_1<\ldots<u_J$ be a sequence of positive integers. Then for each $r>2$
$$\|\|\int \widehat{f}(\xi)1_{\omega_k}(\xi)e^{2\pi i\xi x}d\xi\|_{O_{{\bf U}} \cap V^r_k(L)}\|_{L^2_x(\R)}\lesssim \|f\|_{L^2(\R)},$$
with the implicit constants depending only on $r$.
\end{lemma}

\begin{proof}
It suffices to assume that $0\in\omega_k$ for each $k\in\Z.$
We will say that an interval $[a,b]$ lies in the interior of the interval $[c,d]$ if $c<a<b<d,$ and refer to this property as strong nestedness.
Define a sequence $\ldots<k_{-2}<k_{-1}<k_{0}<k_1<k_2<\ldots$ such that for each $i$ the interval $\omega_{k_{i+1}}$ lies in the interior of $\omega_{k_{i}}$ and none of the intervals $\omega_k$ with $k_i<k<k_{i+1}$ lies in the interior of $\omega_{k_{i}}$.  Define also $f_i$ by $\widehat{f_i}:=1_{\omega_{k_i}}\widehat{f}$. 

We first estimate the $l^{\infty}_k(\Z)$ component of the variational norm. Choose some Schwartz function $\zeta$ with $1_{[-1,1]}\le \widehat{\zeta}\le 1_{[-2,2]}$ and for each dyadic $\omega$ define $\widehat{\zeta}_{\omega}(\xi)=\widehat{\zeta}(\frac{\xi-c(\omega)}{|\omega|})$. Note that
$$\sup_k|f\ast \check{1}_{\omega_k}(x)|\le \sup_i\sup_{k_i\le k<k_{i+1}}|f_i\ast {\zeta}_{\omega_k}(x)|+(\sum_i\sum_{k_i\le k<k_{i+1}}|f_i\ast \check{1}_{\omega_k}(x)-f_i\ast {\zeta}_{\omega_k}(x)|^2)^{1/2}.$$ By Plancherel's formula, the square function above is bounded in $L^2$ by a constant multiple of $\|f\|_2$. Since for each $k$ and $x\in\R$
$$|\zeta_{\omega_k}(x)|\le \text{Dil}_{2^{k}}^{1}|\zeta|(x),$$
we get 
\begin{align*}
\sup_i\sup_{k_i\le k<k_{i+1}}|f_i \ast {\zeta}_{\omega_k}(x)|&\le (\sup_k\text{Dil}_{2^{k}}^{1}|\zeta|)(x)\ast(\sup_i |f_i|)\\&\le M_1(\sup_i |f_i|)(x).
\end{align*}
Finally, to control $\sup_i |f_i|$ we note that 
$$\sup_i |f_i(x)|\le \sup_i|f\ast {\zeta}_{\omega_{k_i}}(x)|+(\sum_i|f_i(x)-f\ast {\zeta}_{\omega_{k_i}}(x)|^2)^{1/2},$$
and then use an argument as above to conclude that
$$\|\sup_i |f_i(x)|\|_{L^2_x(\R)}\lesssim \|f\|_{L^2(\R)}.$$
This shows that 
$$\|\sup_k|f\ast \check{1}_{\omega_k}(x)|\|_{L^2_x(\R)}\lesssim \|f\|_{L^2(\R)}.$$

The estimates for the oscillation norm are now immediate consequences of the maximal estimates and the orthogonality of 
$$f_j:=(\hat{f}1_{\omega_{u_j}\setminus\omega_{u_{j+1}}})\check{\ }.$$
Indeed
$$
(\sum_{j=1}^{J-1}\|\sup_{u_j\le k<u_{j+1}}|(\widehat{f}1_{\omega_{k}})\check{\ }-(\widehat{f}1_{\omega_{u_{j}}})\check{\ }|\|_{L^2(\R)}^2)^{1/2}
=(\sum_{j=1}^{J-1}\|\sup_{u_j\le k<u_{j+1}}|(\widehat{f_j}1_{\omega_{k}})\check{\ }-(\widehat{f_j}1_{\omega_{u_{j}}})\check{\ }|\|_{L^2(\R)}^2)^{1/2}$$
$$\lesssim (\sum_{j=1}^{J-1}\|\sup_{u_j\le k<u_{j+1}}|(\widehat{f_j}1_{\omega_{k}})\check{\ }|\|_{L^2(\R)}^2)^{1/2}\lesssim (\sum_{j=1}^{J-1}\|{f_j}\|_{L^2(\R)}^2)^{1/2}=\|f\|_{L^2(\R)}.
$$

We next focus on the variational part $\|\cdot\|_{\tilde{V}^r}$ of the norm. For each $f:\R\to \R$ define its Poisson integral $P_tf:=f\ast P_t,\;t>0$, where $\widehat{P}_t(\xi):=e^{-t|\xi|}$. The following is a consequence of the variational result of Lepingle \cite{Le} applied to the Brownian martingale associated with the harmonic function $u(x,t):=(f\ast P_t)(x)$ on the upper half plane:  
\begin{equation}
\label{LePiNgLe}
\|\|P_{2^k}f(x)\|_{V^r_k(\Z)}\|_{L^2_x(\R)}\lesssim \|f\|_{L^2(\R)}.
\end{equation}
We will use this result together with the following corollary for averages. For  each  $\lambda\in\R$ and each $k\in\Z$, define  
$$A_k^{\lambda}f(x):=\int \widehat{f}(\xi)1_{[\lambda-2^{-k},\lambda+2^{-k}]}(\xi)e^{2\pi i\xi x}d\xi.$$ Then ~\eqref{LePiNgLe} and a classical square function argument show that
\begin{equation*}
\label{Varnum2ee}
\|\|A_k^{\lambda}f(x)\|_{V^r_k(\Z)}\|_{L^2_x(\R)}\lesssim \|f\|_{L^2(\R)},
\end{equation*}
with the implicit constant independent of $\lambda.$

For simplicity denote
$$f\ast \check{1}_{\omega_k}(x):=A_{\omega_k}f(x).$$
We proceed by estimating
$$\|\|A_{\omega_k}f(x)\|_{\tilde{V}^r_k(\Z)}\|_{L^2_x(\R)} \le \|\|A_{\omega_k}f_{i(k)}(x)\|_{\tilde{V}^r_k(\Z)}\|_{L^2_x(\R)} +\|\|A_{\omega_k}g_{i(k)}(x)\|_{\tilde{V}^r_k(\Z)}\|_{L^2_x(\R)},$$
where $i(k)$ is such that $k_{i(k)-1}\le k<k_{i(k)},$ and $g_{i}=f_{i-1}-f_{i}$. Note that the functions $g_i$ are pairwise orthogonal. 

Since the sequence $A_{\omega_k}f_{i(k)}$ is constant on each block $k_i\le k<k_{i+1}$ it follows that 
\begin{align*}
\|\|A_{\omega_k}f_{i(k)}(x)\|_{\tilde{V}^r_k(\Z)}\|_{L^2_x(\R)}&=\|\|f_i(x)\|_{\tilde{V}^r_i(\Z)}\|_{L^2_x(\R)}\\&\lesssim \|(\sum_{i}|f_i(x)-P_{2^{k_i}}f(x)|^2)^{1/2}\|_{L^2_x(\R)}+\|\|P_{2^k}f(x)\|_{\tilde{V}^r_k(\Z)}\|_{L^2_x(\R)}.
\end{align*}
Now, since $\widehat{P_t}(0)=1$ and by using the decay of $\widehat{P}_1$ we get that for each $\theta\in\R$ 
\begin{align*}
\sum_{i}|1_{\omega_{k_i}}(\theta)-\widehat{P}_{2^{k_i}}(\theta)|^2&\lesssim \sum_{i:\theta\in\omega_{k_i}}|\theta2^{k_i}|^2 +\sum_{i:\theta\notin\omega_{k_i}}|\widehat{P}_1(\theta 2^{k_i})|^2\\&\lesssim \sum_{i:\theta\in\omega_{k_i}}|\theta2^{k_i}|^2+\sum_{i:\theta\notin\omega_{k_i}}|1+\theta2^{k_i}|^{-2}.
\end{align*}
Note that $\theta\in\omega_{k_i}$ implies that $2^{k_i}|\theta|\le 1$ while $\theta\notin\omega_{k_i}$ implies that $2^{k_{i+1}}|\theta|\ge 1$, and thus we get \begin{equation}
\label{Varnum1ee}
\|\|A_{\omega_k}f_{i(k)}(x)\|_{\tilde{V}^r_k(\Z)}\|_{L^2_x(\R)} \lesssim\|f\|_{L^2(\R)}.
\end{equation}
Denote by $\lambda_i$ the common endpoint of all the intervals $\omega_k$, $k_i\le k< k_{i+1}$.
Finally, ~\eqref{Varnum2ee} and the strong nestedness lead us to
\begin{align*}
\|\|A_{\omega_k}g_{i(k)}(x)\|_{\tilde{V}^r_k(\Z)}\|_{L^2_x(\R)} &\lesssim \|(\sum_i|g_i|^2)^{1/2}\|_{L^2(\R)}+\|(\sum_i\|A_k^{\lambda_i}g_i(x)\|_{\tilde{V}^r(k_i\le k<k_{i+1})}^r)^{1/r}\|_{L^2_x(\R)} \\&\lesssim
\|f\|_{L^2(\R)}+ \|(\sum_i\|A_k^{\lambda_i}g_i(x)\|_{\tilde{V}^r(k_i\le k<k_{i+1})}^2)^{1/2}\|_{L^2_x(\R)}\\&\lesssim\|f\|_{L^2(\R)}.
\end{align*}
This and ~\eqref{Varnum1ee} ends the proof of the lemma.
\end{proof}

\begin{proof}{\bf{of Theorem ~\ref{thm:gh098cf}}}
Denote by $\varphi_{k,l}(z):=(\text{Tr}_{-\lambda_l}(m_{k,l}{1}_{\omega_{k,l}}))\check{\ }(z)$ and by $B$ the best constant for which the following inequality holds for each $f_1,\ldots,f_L\in L^2(\R)$ with  $\supp (\hat{f_l})\subseteq [-2,2]$:
$$\|\sup_{k\ge 0}|\sum_{l=1}^{L}e^{2\pi i\lambda_l z}(f_l\ast\varphi_{k,l})(z)|\|_{L^2_z(\R)}\lesssim B(\sum_{l=1}^{L}\|f_l\|_{L^2(\R)}^2)^{1/2}.$$ It suffices to prove that
$$B\lesssim L^{1/2-1/r}\sup_l\sup_{\|g\|_{L^2(\R)}=1}\|\|(m_{k,l}1_{\omega_{k,l}}\widehat{g})\check{\ }(z)\|_{V^r_k(L)}\|_{L^2_z(\R)}.$$ 
For each $0\le y\le \frac1{100}$ we have by Plancherel's theorem that
$$\|f_l-\text{Tr}_yf_l\|_{L^2(\R)}<\frac12\|f_l\|_{L^2(\R)},$$
and hence we can write 

\begin{align*}
\left\|\sup_{k\ge 0}|\sum_{l=1}^{L}e^{2\pi i\lambda_l z}(f_l\ast\varphi_{k,l})(z)|\right\|_{L^2_z(\R)}&\le \left\|\sup_{k\ge 0}|\sum_{l=1}^{L}e^{2\pi i\lambda_l z}\text{Tr}_y(f_l\ast\varphi_{k,l})(z)|\right\|_{L^2_z(\R)}+\\&+\frac{B}{2}(\sum_{l=1}^{L}\|f_l\|_{L^2(\R)}^2)^{1/2}.
\end{align*}
Thus, by integrating in $y$, it suffices to prove that
$$
\left\|\|\sup_{k\ge 0}|\sum_{l=1}^{L}e^{2\pi i\lambda_l y}e^{2\pi i\lambda_l z}(f_l\ast\varphi_{k,l})(z)|\|_{L^2_y([0,1))}\right\|_{L^2_z(\R)}
$$
$$
\lesssim  L^{1/2-1/r}\sup_l\sup_{\|g\|_{L^2(\R)}=1}\left\|\|(m_{k,l}1_{\omega_{k,l}}\widehat{g})\check{\ }(z)\|_{V^r_k(L)}\right\|_{L^2_z(\R)}(\sum_{l=1}^{L}\|f_l\|_{L^2(\R)}^2)^{1/2}.
$$
We can estimate the first term above by first using Lemma ~\ref{lem:auxq2}
 and then Minkowski's inequality on $l^{r/2}(N)$ (for arbitrary $N$) by
\begin{align*}
L^{1/2-1/r}\|\|e^{2\pi i\lambda_l z}&(f_l\ast\varphi_{k,l})(z)\|_{V^r_k(L)}\|_{L^2_z(\R)}\\&\le  L^{1/2-1/r}\left\|(\sum_{l=1}^{L}\|f_l\ast\varphi_{k,l}(z)\|_{V^r_k(L)}^2)^{1/2}\right\|_{L^2_z(\R)}\\&\le L^{1/2-1/r}\sup_l\sup_{\|g\|_{L^2(\R)}=1}\left\|\|g\ast\varphi_{k,l}(z)\|_{V^r_k(L)}\right\|_{L^2_z(\R)}(\sum_{l=1}^{L}\|f_l\|_{L^2(\R)}^2)^{1/2}
\\&= L^{1/2-1/r}\sup_l\sup_{\|g\|_{L^2(\R)}=1}\left\|\|(m_{k,l}1_{\omega_{k,l}}\widehat{g})\check{\ }(z)\right\|_{V^r_k(L)}\|_{L^2_z(\R)}(\sum_{l=1}^{L}\|f_l\|_{L^2(\R)}^2)^{1/2},
\end{align*}
where the  variational norm in the first term above is understood in the Hilbert space $l^2(L)$. 
\end{proof}

An argument very similar to the above also proves the following version of Theorem ~\ref{thm:gh098cf}:
\begin{theorem}
\label{ttzitzijj6}
Consider a collection $R$ of $L$ disjoint dyadic frequency  intervals $\omega$. For each $\omega\in R$ and each $k\in \Z$ let $m_{k,\omega}:\R\to\C$ be a sequence of multipliers. Define 
$$\Delta_kf(x):=\sum_{\omega\in R}\int_{\omega}m_{k,\omega}(\xi)\widehat{f}(\xi)e^{2\pi i \xi x}d\xi.$$ 
Then for each $r>2$
$$\|\sup_k|\Delta_kf(x)|\|_{L^2_x(\R)}\lesssim L^{1/2-1/r} \sup_{\omega\in R}\sup_{\|g\|_{L^2(\R)}=1}\|\|(m_{k,\omega}1_{\omega}\widehat{g})\check{\ }(z)\|_{V^r_k(L)}\|_{L^2_z(\R)}\|f\|_{L^2(\R)},$$
with the implicit constants depending only on $r$.
\end{theorem}

It turns out that the results of Theorems ~\ref{thm:gh098cf} and  ~\ref{ttzitzijj6} are not general enough for our applications, and so we prove the following more general version. Consider now an arbitrary set $\Lambda=\{\lambda_1,\ldots,\lambda_L\}$ with no further restrictions on it, and for each $k\in \Z$ define $R_k$ as before. We now associate to each $\omega\in \bigcup_k R_k$ a multiplier $m_{\omega}:\R\to\C$ and define
\begin{equation}
\label{Bourgeneral}
\Delta_kf(x):=\sum_{\omega\in R_k}\int_{\omega} m_{\omega}(\xi)\widehat{f}(\xi)e^{2\pi i\xi x}d\xi,
\end{equation}
$$\|m_{\omega}\|_{V_2^{r,*}}:=\sup_{l}\sup_{\lambda_l\in\omega_k\in R_k}\sup_{\|g\|_{L^2(\R)}=1}\|\|(m_{\omega_k}1_{\omega_{k}}\widehat{g})\check{\ }(z)\|_{V^r_k(L)}\|_{L^2_z(\R)}.$$
\begin{theorem}
\label{thm:gh098cfhh}
For each $r>2$ we have the inequality
$$\|\sup_{k}|\Delta_kf(x)|\|_{L^2_x(\R)} \lesssim L^{1/2-1/r}\|m_{\omega}\|_{V_2^{r,*}}\|f\|_2,$$
with the implicit constant depending only on $r$.
\end{theorem}
\begin{proof}
It suffices as before to assume that the index $k$ runs through a finite interval $\{a,a+1,\ldots,b\}$ with $a,b\in\Z$.
We can find a sequence  $a=k_0<k_1<\ldots<k_N=b$ with $N\le L$, such that for each $0\le j\le N-1$, $R_k$ has the same cardinality when $k_j\le k<k_{j+1}$. If $\widehat{f_j}:=(\sum_{\omega\in R_{k_j}}1_{\omega}-\sum_{\omega\in R_{k_{j+1}}}1_{\omega})\widehat{f}$, then the functions $f_j$ are pairwise orthogonal. 
We can now bound $\|\sup_k|\Delta_kf(x)|\|_{L^2_x(\R)}$ by 

\begin{equation}
\label{-1-3-5-7ii}
\|\sup_{j}\sup_{k_j\le k<k_{j+1}}|(\sum_{\omega\in R_{k_{j+1}}}m_{{\omega}(k)}1_{\omega}\sum_{j'>j}\widehat{f}_{j'})\check{\ }(x)|\|_{L^2_x(\R)}+
\end{equation} 
\begin{equation}
\label{-1-3-5-7iii}
+\|\sup_{j}\sup_{k_j\le k<k_{j+1}}|(\sum_{\omega\in R_{k}}m_{\omega}1_{\omega}\widehat{f}_{j})\check{\ }(x)|\|_{L^2_x(\R)}.
\end{equation}
For each $\omega\in R_{k_{j+1}}$ and each $k_j\le k<k_{j+1}$, ${\omega}(k)$ is defined to be the interval in $R_k$ containing $\omega.$ Theorem ~\ref{thm:gh098cf} and scaling invariance show that the term ~\eqref{-1-3-5-7iii} can be bounded by 
\begin{align*}
(\sum_{j}\|\sup_{k_j\le k<k_{j+1}}&|(\sum_{\omega\in R_{k}}m_{\omega}1_{\omega}\widehat{f}_{j})\check{\ }(x)|\|_{L^2_x(\R)}^2)^{1/2}\lesssim\\&\lesssim (\sum_{j}L^{1-2/r}\sup_{l}\sup_{\lambda_l\in\omega_k\in R_k\atop{k_j\le k<k_{j+1}}}\sup_{\|g\|_{L^2(\R)}=1}\|\|(m_{\omega_k}1_{\omega_{k}}\widehat{g})\check{\ }(z)\|_{V^r_k(L)}\|_{L^2_z(\R)}^2\|f_j\|_{L^2(\R)}^2)^{1/2}\\&\lesssim L^{1/2-1/r}\|m_{\omega}\|_{V_2^{r,*}}\|f\|_{L^2(\R)}.
\end{align*}

To estimate the term in ~\eqref{-1-3-5-7ii}, define the maximal operators $$O_j^{*}(h)(x):=\sup_{k_j\le k<k_{j+1}}|(\sum_{\omega\in R_{k_{j+1}}}m_{{\omega}(k)}1_{\omega}\widehat{h})\check{\ }(x)|. $$ 

We will  argue that
$$\|\sup_{1\le j\le L}O_{j}^{*}(\sum_{j\le j'\le L}f_{j'})(x)\|_{L^2_x(\R)}
\lesssim  L^{1/2-1/r}\|m_{\omega}\|_{V_2^{r,*}}(\sum_{j=1}^L\|f_j\|_{L^2(\R)}^2)^{1/2}.$$ It suffices to consider only dyadic values of $L$ so we will assume that $L=2^{M}$, for some $M\ge 0$.
For each $0\le m\le M$, denote by  $A_m$ the best constant for which the following inequality holds for all discrete dyadic intervals $J=(j_1,j_2]:=\{j_1+1,j_1+2,\ldots,j_{2}\}\footnote{Here $j_1=a2^b$ and $j_2=(a+1)2^b$  with  $a,b\in\Z_{+}$}\subseteq \{1,2\ldots,2^M\}$ with $2^m$ elements

$$\|\sup_{j\in J}O_{j}^{*}(\sum_{j\le j'\le j_2}f_{j'})(x)\|_{L^2_x(\R)}\lesssim A_m(\sum_{j\in J}\|f_j\|_{L^2(\R)}^2)^{1/2}.$$ We will use a reasoning  similar to the one  in the proof of the Rademacher-Menshov inequality, to argue that $A_{M}\lesssim B_{M}$, where
$$B_{m}:= 2^{m(1/2-1/r)}\|m_{\omega}\|_{V_2^{r,*}}.$$  We can  write  for each $0\le m\le M-1$ and each discrete dyadic interval $J=(j_1,j_2]\subseteq\{1,2,\ldots,2^{M}\}$ having $2^{m+1}$ elements and midpoint $j_3:=j_1+2^m$ 
$$
\|\sup_{j\in J}O_{j}^{*}(\sum_{j\le j'\le j_2}f_{j'})(x)\|_{L^2_x(\R)}^2
\le
\|\sup_{j_3+1\le j\le j_2}O_{j}^{*}(\sum_{j\le j'\le j_2}f_{j'})(x)\|_{L^2_x(\R)}^2+$$
$$
+\left(\|\sup_{j_1+1\le j\le j_3}O_{j}^{*}(\sum_{j \le j'\le j_3}f_{j'})(x)\|_{L^2_x(\R)}+\|\sup_{j_1+1\le j\le j_3}O_{j}^{*}(\sum_{j_3+1\le j'\le j_2}f_{j'})(x)\|_{L^2_x(\R)}\right)^2.
$$
We then use the definition of $A_m$ for the first two terms above and Theorem ~\ref{ttzitzijj6} for the third one, to bound the sum above by 
\begin{equation*}
 A_{m}^2\sum_{j_3+1\le j'\le j_2}\|f_{j'}\|_{L^2(\R)}^2+ (A_{m}(\sum_{j_1\le j'\le j_3}\|f_{j'}\|_{L^2(\R)}^2)^{1/2}+CB_{m}(\sum_{j_3+1\le j'\le j_2}\|f_{j'}\|_{L^2(\R)}^2)^{1/2})^{2}
\end{equation*}
\begin{equation*}
\le (A_{m}+CB_{m})^2\sum_{j\in J}\|f_{j}\|_{L^2(\R)}^2.
\end{equation*}
We conclude that $A_{{m+1}}\le A_{m}+CB_{m}$ for each $0\le m \le M-1$, which together with the fact that $A_0=0$ proves that $A_{M}\lesssim B_{M}.$
\end{proof}
\begin{remark}
In our later applications of Theorem ~\ref{thm:gh098cfhh} the parameter $r$ will be chosen sufficiently close to $2$, making the dependency on $L$ of the $L^2(\R)$ norm of the weighted  maximal operator  negligible. The fact that the $L^2(\R)$ norm goes to $\infty$ as $L$ gets larger follows from the result in \cite{BKO}, where it is proven that this norm is at least of the order of $(\log L)^{1/4}.$  
\end{remark}

If the multipliers $m_{\omega}$ in the above theorem are chosen to be the constant function 1, we recover (modulo a slightly larger $L$ bound) the result of Bourgain from \cite{Bo1}, via the variational estimates in Lemma ~\ref{lem:auxq4}.  Bourgain's  $L$ bound is $(\log L)^2$ rather then $L^{1/2-1/r}$, however our  slightly larger bound will suffice for our application, since we will take $r$ as close to 2 as we want. We state Bourgain's  result for future reference.

\begin{corollary}
\label{labellhgt64}
Assume we are in the setting of Theorem ~\ref{thm:gh098cfhh} and that $m_{\omega}\equiv 1$ for each $\omega$.
For each $r>2$ we have the inequality
$$\|\sup_{k}|\Delta_kf(x)|\|_{L^2_x(\R)}\lesssim L^{1/2-1/r}\|f\|_{L^2(\R)},$$
with the implicit constant depending only on $r$.
\end{corollary}

An interesting question regards the dependency on $L$ of the $L^q(\R)$ norm of the operator in Theorem ~\ref{thm:gh098cfhh}, for $q\in(1,2)\cup(2,\infty)$. The fact that this norm is large as a function of $L$ is already apparent at a single scale. Due to the equality $\|m\|_{M_q}=\|m\|_{M_{q'}}$ for dual pairs $(q,q')$, it suffices to note the following. 
\begin{proposition}
For each $L\in\N$ and $q\in (2,\infty)$ there is a choice of signs $(\varepsilon_l)_{1\le l\le L}$ such that if $\widehat{f}_L=1_{[0,L]}$ then
$$\|\int \widehat{f}_L(\xi)\sum_{l=0}^{L-1}\varepsilon_l1_{[l,l+1]}(\xi)e^{2\pi i\xi x}d\xi\|_{L^q_x(\R)}\gtrsim L^{|1/2-1/q|}\|f_L\|_{L^q(\R)}.$$
 \end{proposition} 
\begin{proof}
 It  immediately follows that  $$\|f_L\|_{L^q(\R)}\sim L^{1-1/q},$$ $$\|(\sum_{l=0}^{L-1}|\int \widehat{f}_L(\xi)1_{[l,l+1]}(\xi)e^{2\pi i\xi x}d\xi|^2)^{1/2}\|_{L^q_x(\R)}\sim L^{1/2}.$$ Khintchine's inequality ends the proof.
\end{proof}
This shows that the $L^q$ norm of the maximal operator $\Delta^{*}f(x):=\sup_k|\Delta_kf(x)|$, with $\Delta_k$ defined in ~\eqref{Bourgeneral}, satisfies
$$\|\Delta^{*}\|_{L^q(\R) \to L^q(\R)}\gtrsim L^{|1/2-1/q|}\|m_{\omega}\|_{V_2^{r,*}}.$$

Theorem ~\ref{thm:gh098cfhh} will be used in Section ~\ref{sec:PointEst} to control maximal operators. For the proof of the oscillation inequality leading to the dense class results, we will need a  more general version of Theorem ~\ref{thm:gh098cfhh}. 

We will assume that $L$, $\lambda_l$, $R_k$, $m_{\omega}$ and $\Delta_kf$ are  as in Theorem ~\ref{thm:gh098cfhh}. 
Of relevance for the estimates in the next theorem is the following multiplier norm

$$\|m_{\omega}\|_{O_{{\bf U}}\cap V^{r,*}}:=\sup_{l}\sup_{\lambda_l\in\omega_k\in R_k}\sup_{\|g\|_{L^2(\R)}=1}\left\|\|(m_{\omega_k}1_{\omega_{k}}\widehat{g})\check{\ }(z)\|_{O_{{\bf U}} \cap V^r_k(L)}\right\|_{L^2_z(\R)},$$
where ${\bf U}:=u_1<u_2<\ldots<u_J$ is an arbitrary  finite sequence of integers.

\begin{theorem}
\label{Bourgeneralo}
Let ${\bf U}:=u_1<u_2<\ldots<u_J$ be an arbitrary finite sequence of integers.
The following inequality holds for each $r>2$
$$(\sum_{j=1}^{J-1}\|\sup_{u_j\le k <u_{j+1}}|\Delta_kf(x)-\Delta_{u_{j}}f(x)|\|_{L^2_x(\R)}^2)^{1/2}\lesssim J^{\frac{r+4}{4r+4}}L^{1-2/r}\|m_{\omega}\|_{O_{{\bf U}} \cap V^{r,*}}\|f\|_{L^2(\R)},$$
with the implicit constant depending only on $r$. (in particular it does not depend on either $J$ or $u_1,\ldots,u_J$).
\end{theorem}

\begin{remark}
The only relevant thing about the exponents $\frac{r+4}{4r+4}$ and $1-\frac2{r}$ is that the first is  less than $\frac12$ and $\lim_{r\to 2}(1-\frac2r )=0.$
\end{remark}

\begin{proof}
To prove the above theorem we need two inequalities. In the first inequality we  aim for small $L$ dependent bounds and tolerate a trivial $J$ dependent bound. On the other hand, in the second inequality we aim for a $J$ independent bound but we  will tolerate a big $L$ dependent bound. 

Note that for each $j$ we have as a consequence of Theorem ~\ref{thm:gh098cfhh} 
\begin{align*}
\|\sup_{u_j\le k<u_{j+1}}|\Delta_kf(x)-\Delta_{u_{j}}f(x)|\|_{L^2_x(\R)}&\lesssim \|\sup_k|\Delta_kf(x)|\|_{L^2_x(\R)} 
\\&\lesssim L^{1/2-1/r}\|m_{\omega}\|_{O_{{\bf U}} \cap V^{r,*}}\|f\|_{L^2(\R)}.\end{align*}
Thus, we get our first main inequality by doing rough estimates: 
\begin{equation}
\label{eerrthyu:1}
(\sum_{j=1}^{J-1}\|\sup_{u_j\le k <u_{j+1}}|\Delta_kf-\Delta_{u_{j}}f|\|_2^2)^{1/2}
\lesssim J^{1/2}L^{1/2-1/r}\|m_{\omega}\|_{O_{{\bf U}}\cap V^r}\|f\|_2.
\end{equation}

Fix now some $1\le j\le J-1$. For each interval $\omega\in R_{u_j}$ pick some $l\in\{1,\ldots,L\}$ such that $\lambda_l\in\omega$. Denote by $\Lambda(j)$ the set of all these $l$. For each $l\in\Lambda(j)$ and each $u_j\le k<u_{j+1}$, denote as before by $\omega_{k,l}$ the interval in $R_k$ containing $\lambda_l$. Define for each $z\in\R$ 
$$a_{\omega}(z):=(m_{\omega}1_{\omega}\hat{f})\check{\ }(z).$$

Since $|R_{u_j}|\le |R_k|\le |R_{u_j}|+L$ when $u_j\le k<u_{j+1}$, we can evaluate 
\begin{align*}
\sup_{u_j\le k<u_{j+1}}|\Delta_kf(z)-\Delta_{u_{j}}f(z)|^2&=\sup_{u_j\le k<u_{j+1}}|\sum_{\omega\in R_k}a_{\omega}(z)-\sum_{\omega\in R_{u_j}}a_{\omega}(z)|^2\\& \lesssim  \sup_{u_j\le k<u_{j+1}}|\sum_{l\in\Lambda(j)}(a_{\omega_{k,l}}(z)-a_{\omega_{u_j,l}}(z))|^2+L^2\sup_{k,l}|a_{\omega_{k,l}}(z)|^2\\&\le  L\sum_{l\in\Lambda(j)}\sup_{u_j\le k<u_{j+1}}|a_{\omega_{k,l}}(z)-a_{\omega_{u_j,l}}(z)|^2+L^2\sup_{k,l}|a_{\omega_{k,l}}(z)|^2.
\end{align*} 
Note also that  there is a set $V\subset \{1,\ldots J-1\}$ with at most $L$ elements such that $|R_{u_j}|\not=|R_{u_{j+1}}|$ for $j\in V$, and if $j\notin V$ then we can improve on the bound obtained above
$$\sup_{u_j\le k<u_{j+1}}|\Delta_kf(z)-\Delta_{u_{j}}f(z)|^2\lesssim  L\sum_{l\in\Lambda(j)}\sup_{u_j\le k<u_{j+1}}|a_{\omega_{k,l}}(z)-a_{\omega_{u_j,l}}(z)|^2.$$
For each $l\in\Lambda(j)$ define $x_l^{(j)}(z):=\sup_{u_j\le k<u_{j+1}}|(a_{\omega_{k,l}}(z)-a_{\omega_{u_j,l}}(z))|$. By summing over $j$ we get 
\begin{align*}
\sum_{j=1}^{J-1}\sup_{u_j\le k<u_{j+1}}|\Delta_kf(z)-\Delta_{u_{j}}f(z)|^2&\lesssim L\sum_{j=1}^{J-1}\sum_{l\in\Lambda(j)}(x_l^{(j)}(z))^2+L^3\sup_{k,l}|a_{\omega_{k,l}}(z)|^2\\&\le L\sum_{l\in\Lambda}\sum_{j=1}^{J-1}1_{\Lambda(j)}(l)(x_l^{(j)}(z))^2+L^3\sup_{k,l}|a_{\omega_{k,l}}(z)|^2\\&\le L\sum_{l\in\Lambda}\|a_{\omega_{k,l}}(z)\|_{O_{{\bf U}}}^2+L^3\sup_{k,l}|a_{\omega_{k,l}}(z)|^2\\&\lesssim L^3\sup_{l\in\Lambda}\|a_{\omega_{k,l}}(z)\|_{O_{{\bf U}} \cap V^r_k(\Z)}^2.
\end{align*}
This together with integration with respect to $z$ produces the second main inequality

\begin{equation}
\label{eerrthyu:2}
\left(\sum_{j=1}^{J-1}\|\sup_{u_j\le k <u_{j+1}}|\Delta_kf(x)-\Delta_{u_{j}}f(x)|\|_{L^2_x(\R)}^2\right)^{1/2}
\lesssim L^{3/2}\|m_{\omega}\|_{O_{{\bf U}} \cap V^{r,*}}\|f\|_{L^2(\R)}.
\end{equation}

Finally, by interpolating between  ~\eqref{eerrthyu:1} and ~\eqref{eerrthyu:2} it follows that
\begin{align*} 
\left(\sum_{j=1}^{J-1}\|\sup_{u_j\le k <u_{j+1}}|\Delta_kf(x)-\Delta_{u_{j}}f(x)|\|_{L^2_x(\R)}^2\right)^{1/2}&\lesssim (J^{1/2}L^{1/2-1/r})^{\frac{r+4}{2r+2}}(L^{3/2})^{\frac{r-2}{2r+2}}\|m_{\omega}\|_{O_{{\bf U}} \cap V^{r,*}}\|f\|_{L^2(\R)}\\&= J^{\frac{r+4}{4r+4}}L^{1-2/r}\|m_{\omega}\|_{O_{{\bf U}} \cap V^{r,*}}\|f\|_{L^2(\R)}. 
\end{align*}
\end{proof}

We also record the following immediate corollary.

\begin{corollary}
\label{labellhgt64o}
Assume we are in the setting of Theorem ~\ref{Bourgeneralo} and that $m_{\omega}\equiv 1$ for each $\omega$.
For each $r>2$ we have the inequality
$$(\sum_{j=1}^{J-1}\|\sup_{u_j\le k <u_{j+1}}|\Delta_kf(x)-\Delta_{u_{j}}f(x)|\|_{L^2_x(\R)}^2)^{1/2}\lesssim J^{\frac{r+4}{4r+4}}L^{1-2/r}\|f\|_{L^2(\R)},$$
with the implicit constant depending only on $r$.
\end{corollary}

\section{Variational, oscillation and square function  estimates}
\label{sec:8}

In this section we prove a few  auxiliary results of general interest, which combined with Theorem ~\ref{thm:gh098cfhh} will be used later to control the measure of various exceptional sets. The following result is classical.

\begin{proposition}
\label{p.BMO1}
Let $\D_0$ be a finite collection of dyadic intervals included into some interval $\I$, each of which is associated with a function $\phi_I$ satisfying:
\begin{equation}
\label{e.adap1}
\int\phi_I(x)dx=0,
\end{equation}
\begin{equation}
\label{e.adap2}
\phi_I\mbox{\;is C-adapted to\;}I
\end{equation}
If $a_I\in\C$ are such that 
$$(\frac1{|I_0|}\sum_{I\in \D_0\atop{I\subseteq I_0}}|a_I|^2)^{1/2}\le B,$$ for each dyadic $I_0$, then 
$$\|\sum_{I\in \D_0}a_I\phi_I\|_{\BMO(\R)}\lesssim CB,$$
\begin{equation}
\label{aggght}
\|\sum_{I\in \D_0}a_I\phi_I\|_{L^s(\R)}\lesssim CB|\I|^{1/s},\;\;1<s<\infty,
\end{equation}
with the implicit constants depending only on $s$.
\end{proposition}
\begin{proof}
The $\BMO$ estimate follows as in Proposition ~\ref{p.BMO2}, we do not insist on the details here. 
The  estimate ~\eqref{aggght} is an immediate consequence of the first estimate,  John Nirenberg's inequality and the fact that
\begin{equation}
\label{ghycg}
\sum_{I\in \D_0}|a_I\phi_I(x)|\lesssim C\chi_{\I}^3(x),\;\;x\notin 2\I.
\end{equation}
\end{proof}
The next lemma will be used to prove  bounds on the variational norms operators. 
\begin{lemma}
\label{l.BMO2}
Let $\D_0$ be a finite collection of dyadic intervals $I$ each of which is associated with a function $\phi_I$ satisfying ~\eqref{e.adap1} and ~\eqref{e.adap2} for a fixed $C$. Let $\zeta$ be some fixed Schwartz function  with $1_{[-1,1]}\le\widehat{\zeta}\le 1_{[-2,2]}$. Then for each $a_I\in \C$ 
$$\sum_{k\in\Z}\|\sum_{I\in \D_0\atop{|I|\ge 2^k}}a_I\phi_I-\sum_{I\in \D_0}a_I{\phi}_I*\operatorname{Dil}_{2^k}^1\zeta\|_{L^2(\R)}^2\lesssim C^2\sum_{I\in \D_0}|a_I|^2.$$ 
\end{lemma}
\begin{proof}
We start by making a few observations. Define $\Phi_I(x):=\phi_I(x+c(I))$ and note the following consequences of ~\eqref{e.adap1} and ~\eqref{e.adap2}:
\begin{equation*}
\widehat{\Phi_I}(0)=0,\;\;\; \|\frac{d}{d\xi}\widehat{\Phi_I}(\xi)\|_{\infty} \lesssim |I|^{3/2}
\end{equation*}
\begin{equation}
\label{e.adaphot2}
|\widehat{\Phi_I}(\xi)|\lesssim |\xi| |\frac{d}{d\xi}\widehat{\Phi_I}(\xi)| \lesssim C|\xi||I|^{3/2}
\end{equation}
\begin{equation}
\label{e.adaphot3}
\|\frac{d^2}{d\xi^2}\widehat{\Phi_I}(\xi)\|_{L^\infty_\xi(\R)}\lesssim C|I|^{5/2}
\end{equation}
\begin{equation}
\label{e.adaphot4}
|\frac{d^{j}}{d\xi^{j}}\widehat{\Phi_I}(\xi)|\lesssim \frac{1}{|\xi|}\|\frac{d}{dx}(x^j\Phi_I(x))\|_{L^1_x(\R)}\lesssim C\frac{|I|^{j-\frac12}}{|\xi|},\;\;0\le j\le 2.
\end{equation}

Fix $t\ge 0$. The almost orthogonal behavior of the collection $\phi_{I,k}:=\phi_I-\phi_I\ast \text{Dil}_{2^{k}}^{1}\zeta$, with $|I|=2^{k+t}$, is quantified by the following properties:
$$|\langle\phi_{I,k},\phi_{J,k}\rangle|\lesssim \int_{2^{-k}}^{\infty}|\widehat{\Phi_I}(\xi)\widehat{\Phi_J}(\xi)|d\xi\lesssim C^22^{-t}, \;\;\hbox{a consequence of}\; ~\eqref{e.adaphot4},$$
\begin{align*}
|\langle\phi_{I,k},\phi_{J,k}\rangle|&=|\int\widehat{\Phi_I}(\xi)\bar{\widehat{\Phi_J}}(\xi)(1-\widehat{\zeta}(2^k\xi))^2e^{2\pi i(c(J)-c(I))\xi}d\xi|\\&\lesssim \frac{1}{|c(J)-c(I)|^2}\int|\frac{d^2}{d\xi^2}(\widehat{\Phi_I}(\xi)\bar{\widehat{\Phi_J}}(\xi)(1-\widehat{\zeta}(2^k\xi))^2)|d\xi\\&\lesssim 2^{-t}\left(\frac{C|I|}{|c(J)-c(I)|}\right)^2, \;\;\hbox{a consequence of}\; ~\eqref{e.adaphot4},
\end{align*}
and hence for each $|I|=|J|=2^{k+t}$
$$|\langle\phi_{I,k},\phi_{J,k}\rangle|\lesssim\frac{C^22^{-t}}{(1+\frac{|c(J)-c(I)|}{2^{k+t}})^2}.$$
An immediate corollary of this is that
$$\|\sum_{I\in \D_0\atop{|I|=2^{k+t}}}a_I\phi_{I,k}\|_{L^2(\R)}^2\lesssim C^22^{-t}\sum_{I\in \D_0\atop{|I|=2^{k+t}}}|a_I|^2.$$

An application of the triangle inequality first and then Minkowski's inequality gives
\begin{align*}
\sum_{k\in\Z}\|\sum_{I\in \D_0\atop{|I|\ge 2^k}}a_I\phi_I-\F^{-1}(\sum_{I\in \D_0\atop{|I|\ge 2^k}}a_I\widehat{\phi}_I(\xi)\widehat{\zeta}(2^k\xi))\|_{L^2(\R)}^2
&\lesssim \sum_{k\in\Z}(\sum_{t\ge 0}\|\sum_{I\in \D_0\atop{|I|= 2^{k+t}}}a_I\phi_{I,k}\|_{L^2(\R)})^2\\&
\lesssim (\sum_{t\ge 0}(\sum_{k\in \Z}2^{-t}\sum_{I\in \D_0\atop{|I|= 2^{k+t}}}|a_I|^2)^{1/2})^2\\&\lesssim C^2\sum_{I\in \D_0}|a_I|^2.
\end{align*} 

Fix now $t<0$. The almost orthogonal behavior of the collection $\tilde{\phi}_{I,k}:=\phi_I\ast \text{Dil}_{2^{-k}}^{1}\zeta$, with $|I|=2^{k+t}$, follows as before, by now invoking ~\eqref{e.adaphot2} and ~\eqref{e.adaphot3} instead:
$$|\langle\tilde{\phi}_{I,k},\tilde{\phi}_{J,k}\rangle|\lesssim \int_{0}^{2^{1-k}}|\widehat{\Phi_I}(\xi)\widehat{\Phi_J}(\xi)|d\xi\lesssim 2^{-k}\|\widehat{\Phi_I}\|_{L^\infty(\R)}\|\widehat{\Phi_J}\|_{L^\infty(\R)}\lesssim C^22^{t},$$
$$|\langle\tilde{\phi}_{I,k},\tilde{\phi}_{J,k}\rangle|\lesssim\frac{C^22^{t}}{(1+\frac{|c(J)-c(I)|}{2^{k+t}})^2}.$$
We obtain as before
$$\|\F^{-1}(\sum_{I\in \D_0\atop{|I|< 2^k}}a_I\widehat{\phi}_I(\xi)\widehat{\zeta}(2^k\xi))(x)\|_{L^2_x(\R)}^2\lesssim C^2\sum_{I\in \D_0}|a_I|^2.
$$
\end{proof}

\begin{proposition}
\label{p.BMO2}
Let $\D_0$ be a finite collection of dyadic intervals $I$ contained by  some interval $\I$, each of which is associated with a function $\phi_I$ satisfying  ~\eqref{e.adap1} and ~\eqref{e.adap2} for a fixed $C$.
Consider also a sequence ${\bf U}:=(u_j)_{j=-\infty}^{\infty}$ of integers. If $a_I\in\C$ are such that 
\begin{equation}
\label{e:bmoj}
(\frac1{|I_0|}\sum_{I\in \D_0\atop{I\subseteq I_0}}|a_I|^2)^{1/2}\le B
\end{equation}
 for each dyadic $I_0$, then 
$$\|\|\sum_{I\in \D_0\atop{|I|\ge 2^k}}a_I\phi_I(x)\|_{O_{{\bf U}} \cap V^r_k(L)}\|_{\BMO_x(\R)}+\|(\sum_{k\in\Z}|\sum_{I\in \D_0\atop{|I|=2^k}}a_I\phi_I|^2)^{1/2}\|_{\BMO(\R)}\lesssim CB,\;\;r>2$$
$$\|\|\sum_{I\in \D_0\atop{|I|\ge 2^k}}a_I\phi_I\|_{O_{{\bf U}}\cap V^r_k(L)}\|_{L^s(\R)}+\|(\sum_{k\in\Z}|\sum_{I\in \D_0\atop{|I|=2^k}}a_I\phi_I|^2)^{1/2}\|_{L^s(\R)}\lesssim CB|\I|^{1/s},\;\;r>2,\;1<s<\infty,$$ 
with the implicit constants depending only on $s$ and $r$.
\end{proposition}
\begin{proof}
We will only prove the variational estimates, the argument for the oscillation and square function inequalities follows a very similar path. 
It suffices to prove the $\BMO$ bound. Indeed, this together with  John Nirenberg's inequality, trivial estimates of the $V^r$ norm by the $V^1$ and $V^2$ norms and ~\eqref{ghycg} will immediately give the desired $L^s$ estimate.

Consider some arbitrary interval $J$ and define  $\D_1=\{I\in \D_0:|J|>|I|,\; I\subseteq 4J\}$, $\D_2=\{I\in \D_0:|J|>|I|,\;I\cap (4J)^c\not=\emptyset\}$, $\D_3=\{I\in \D_0:|J|\le |I|\}$ and $b_{J}=\|\sum_{I\in \D_3\atop{|I|\ge 2^k}}a_I\phi_I(c(J))\|_{V^r_k(L)}$ . Define also for $i=1,2,3$
$$F_i(x)=\|\sum_{I\in \D_i\atop{|I|\ge 2^k}}a_I\phi_I(x)\|_{V^r_k(L)}.$$
It suffices to prove the following $\BMO$ estimates for each $F_i$ separately, with the implicit constant depending only on $r$ 
\begin{equation}
\label{engtyuij}
\int_{J}|F_i(x)-b_{J}|^2\lesssim C^2B^2|J|.
\end{equation}
We start with estimates for $F_1$ and write
\begin{align}
\nonumber
\int_J|F_1|^2&
\label{svd1}
\lesssim\|F_1\|_{L^2(\R)}^2\\&\lesssim \left\|\|\sum_{I\in \D_1\atop{|I|\ge 2^k}}a_I\phi_I(x)-\sum_{I\in \D_1}a_I{\phi}_I*\operatorname{Dil}_{2^k}^1\zeta(x)	\|_{V^r_k(L)}\right\|_{L^2_x(\R)}^2\\&
\label{svd2}
+\left\|\|\sum_{I\in \D_1}a_I{\phi}_I*\operatorname{Dil}_{2^k}^1\zeta(x)-P_{2^k}(\sum_{I\in \D_1}a_I{\phi}_I)(x)\|_{V^r_k(L)}\right\|_{L^2_x(\R)}^2\\&+
\label{svd3}
\left\|\|P_{2^k}(\sum_{I\in \D_1}a_I{\phi}_I)(x)\|_{V^r_k(L)}\right\|_{L^2_x(\R)}^2.
\end{align}
Using the result of the previous lemma and estimate ~\eqref{e:bmoj} we easily bound the term ~\eqref{svd1} by a universal constant multiple of $C^2B^2|J|$. Then the mean zero of $\zeta-P_1$ and ~\eqref{aggght} show that
$$
~\eqref{svd2}\lesssim \sum_k\left\|\sum_{I\in \D_1}a_I{\phi}_I*\operatorname{Dil}_{2^k}^1(\zeta-P_1)\right\|_{L^2(\R)}^2\lesssim\|\sum_{I\in \D_1}a_I{\phi}_I\|_{L^2(\R)}^2\lesssim C^2B^2|J|,
$$
while ~\eqref{LePiNgLe} and ~\eqref{aggght} imply the same bound for ~\eqref{svd3}. 

The terms corresponding to $F_2$ and $F_3$ are estimated trivially. First, for each $x\in J$  ~\eqref{e.adap2} and ~\eqref{e:bmoj} imply that
$$F_2(x)\lesssim \sum_{I\in\D_2}|a_I\phi_I(x)|\lesssim BC,$$
and hence
$$\int_J|F_2|^2\lesssim  B^2C^2|J|.$$
On the other hand, for each $x\in J$ ~\eqref{e.adap2} and ~\eqref{e:bmoj} imply that
$$|F_3(x)-b_J|\lesssim \sum_{|I|\in\D_3}|a_I||\phi_I(x)-\phi_I(c(J))|\lesssim BC,$$ and so again
$$\int_J|F_3-b_j|^2\lesssim  B^2C^2|J|.$$ An application of triangle's inequality finishes the proof of ~\eqref{engtyuij} and of the proposition.
\end{proof}

Proposition ~\ref{p.BMO2} and the discussion from the end of the  Section ~\ref{sec:trees} implies the following fundamental estimates for a 2-quasitree.

\begin{theorem}
\label{BMOtree}
Let ${\bf U}:=(u_j)_{j=-\infty}^{\infty}$ be an arbitrary sequence  of integers. For each 2-quasitree $\T$  with top $(I_{\T},\xi_{\T})$, each $l,M\ge 0$, $r>2$ and $1<t<\infty$
$$\|\|\sum_{s\in\T\atop{|I_s|< 2^{k}}}\langle f, \varphi_s\rangle \phi_{s,\T}^{(l)}(x,\xi_{\T})\|_{O_{{\bf U}} \cap V^r_k(L)}\|_{\BMO_x(\R)}+
\|(\sum_{k\in\Z}|\sum_{s\in\T\atop{|I|=2^k}}\langle f, \varphi_s\rangle \phi_{s,\T}^{(l)}(x,\xi_{\T})|^2)^{1/2}\|_{\BMO_x(\R)}$$
$$\lesssim 2^{-Ml}\size(\T)$$
and
$$\|\|\sum_{s\in\T\atop{|I_s|< 2^{k}}}\langle f, \varphi_s\rangle \phi_{s,\T}^{(l)}(x,\xi_{\T})\|_{O_{{\bf U}}\cap V^r_k(L)}\|_{L^t_x(\R)}+\|(\sum_{k\in\Z}|\sum_{s\in\T\atop{|I|=2^k}}\langle f, \varphi_s\rangle \phi_{s,\T}^{(l)}(x,\xi_{\T})|^2)^{1/2}\|_{L^t_x(\R)}$$
$$\lesssim 2^{-Ml}\size(\T)|I_{\T}|^{1/t},$$
with the implicit constants depending only on $r$, $t$ and $M$.
\end{theorem}

\section{Pointwise estimates outside  exceptional sets}
\label{sec:PointEst}

Let  $\P$ be a finite set of tiles which can be written as a  disjoint union of trees $\T$ with tops $T$
$$\P=\bigcup_{\T\in\F}\T.$$
To quantify better the contribution to various model sums, coming from individual tiles, we need to reorganize the collection $\F$ in a more suitable way. For each  $\T\in\F$ define its saturation 
$$G(\T):=\{s\in \P:\omega_T\subseteq\omega_s\}.$$
For the purpose of organizing $G(\T)$ as a collection of  disjoint and better spatially localized quasitrees we  define for each $l\ge 0$ and $m\in\Z$ the quasitree
$\T_{l,m}$ to include all tiles $s\in G(\T)$ satisfying  the following requirements:
\begin{itemize}
\item $|I_s\cap (2^lI_T+2^lm|I_T|)|\ge \frac{|I_s|}{2}$
\item either $|I_s\cap (2^lI_T+2^lm|I_T|)|\not= \frac{|I_s|}{2}$ or $|I_s\cap (2^lI_T+2^l(m-1)|I_T|)|\not= \frac{|I_s|}{2}.$
\end{itemize} 
Obviously, for each $l\ge 0$ the collection consisting of $(\T_{l,m})_{m\in\Z}$ forms a partition of $G(\T)$ into quasitrees. The top  of $\T_{l,m}$ is formally assigned to be the pair $(I_{\T_{l,m}},\xi_{\T})$, where $I_{\T_{l,m}}$ is  the interval $2\times(2^lI_T+2^lm|I_T|)$ while $\xi_{\T}$ is the frequency component of the top $(I_{\T},\xi_{\T})$ of the tree $\T$ (considered as a quasitree). 

Let  $\T_{l,m}=\T_{l,m}^{(1)}\cup\T_{l,m}^{(2)}$ be the standard decomposition of $\T_{l,m}$, where both $\T_{l,m}^{(1)}$ and $\T_{l,m}^{(2)}$ are formally assigned the same top  as $\T_{l,m}$. Denote by $\F_{l,m}$, $\F_{l,m}^{(1)}$ and $\F_{l,m}^{(2)}$ the collections of all the quasitrees $\T_{l,m}$, $\T_{l,m}^{(1)}$ and $\T_{l,m}^{(2)}$, respectively.

Consider $\sigma,\gamma>0$, $\beta\ge 1$,  $r>2$ and the complex numbers $a_s,s\in\P.$ Let also $u_1<\ldots<u_{J}$ be an arbitrary finite  sequence of integers. The first result in this section is the crucial estimate behind Theorem ~\ref{ceamaiceahh}.

\begin{theorem}
\label{mainineq}
Assume we are in the settings from above and also that the following additional requirement is satisfied
\begin{equation}
\label{BMNCARseq}
\sup_{s\in\P}\frac{|a_s|}{|I_s|^{1/2}}\le \sigma.
\end{equation}
Define the  exceptional sets
\begin{align*}
E^{(1)}&:=\bigcup_{l\ge 0}\{x:\sum_{\T\in\F}1_{2^lI_T}(x)>\beta2^{2l}\},\\
E^{(2)}&:=\bigcup_{l,m\ge 0}\bigcup_{\T\in\F_{l,m}^{(2)}}\{x:\|\sum_{s\in \T\atop{|I_s|<2^j}}a_s\phi_{s,\T}^{(\alpha(l,m))}(x,\xi_{\T})\|_{V^r_j(\Z	)}>\gamma2^{-l}(|m|+1)^{-2}\},\\
E^{(3)}&:= \bigcup_{l\ge 0}\bigcup_{\T\in\F_{l+1,0}^{(2)}}\{x:\|\sum_{s\in \T\atop{|I_s|<2^j}}a_s\phi_{s,\T}^{(l)}(x,\xi_{\T})\|_{V^r_j(\Z)}>\gamma2^{-l}\},
\end{align*}
where the symbol $\alpha(l,m)$ equals $l$ if $m=0$ and $l+[\log_2|m|]$ if $m\not=0$. 

Then  for each $x\notin E^{(1)}\cup E^{(2)}\cup E^{(3)}$ we have 
\begin{equation}\label{forex}
\|(\sum_{s\in \P\atop{|I_s|< 2^k}}a_s\phi_s(x,\theta))_{k\in\Z}\|_{M_{2,\theta}^*(\R)}\lesssim \beta^{1/2-1/r}(\gamma+\sigma),
\end{equation}
with the implicit constants depending only on $r$.
\end{theorem}
\begin{proof}
For each $l\ge 0$  and each $x\in \R$ define inductively
\begin{align*}
\F_{0,x} &:=\{\T\in\F,x\in I_T\} \\
\F_{l,x} &:=\{\T\in\F,x\in 2^lI_T\setminus 2^{l-1}I_T\},\;\;l\ge 1\\
\P_{0,x} &:=\bigcup_{\T\in \F_{0,x}}G(\T)\\
\P_{l,x}&:=\bigcup_{\T\in \F_{l,x}}G(\T)\setminus\bigcup_{l'<l}\P_{l',x},\;\;l\ge 1\\
\Xi_{x,l}&:=\{c(\omega_T):\;\T\in \F_{l,x}\}.
\end{align*}
Note that for each $x\in\R$, $\{\P_{l,x}\}_{l\ge 0}$ forms a partition of $\P$. Since $x\notin E^{(1)}$, it also follows that $\sharp \Xi_{x,l}\le \beta 2^{2l}$.

Fix  $x\not\in E^{(1)}\cup E^{(2)}\cup E^{(3)}$  and focus on estimates for the left-hand side of \eqref{forex}.
For each $k\in\Z$ and $l\ge 0$ 
let $\Omega_{k,l}$ be the collection of dyadic frequency intervals of 
length $2^{-k}$  which contain an element of $\Xi_{x,l}$.
Let $\tilde{\Omega}_{k,l}$ be the collection of all (dyadic) siblings of
intervals in $\Omega_{k-1,l}$ that are not themselves in $\Omega_{k,l}$.
Observe that both $\bigcup_{k'} \tilde{\Omega}_{k',l}$ and $\Omega_{k,l}\cup \bigcup_{k'\le k}\tilde{\Omega}_{k',l}$ are  collections of pairwise disjoint intervals which cover $\{\;\omega_{s,2}:s\in \P_{l,x}\}$. Moreover we can write 
\begin{align*}
\sum_{s\in\P_{l,x} :|I_s|< 2^k} a_s\phi_{s}(x,\theta) &=\sum_{\omega \in \Omega_{k,l}} 1_{\omega}(\theta)
\sum_{s\in\P_{l,x}\atop{|I_s|< 2^k,\;\omega\cap \omega_{s,2}\not=\emptyset}}a_s \phi_s(x,\theta) \\
&\quad +
\sum_{k'\le k} 
\sum_{\omega \in \tilde{\Omega}_{k',l}}
1_{\omega}(\theta) \sum_{s\in\P_{l,x}\atop{|I_s|< 2^{k},\;\omega \cap \omega_{s,2}\not=\emptyset}}a_s \phi_s(x,\theta).
\end{align*}
Indeed, if $1_\omega(\theta)\phi_s(x,\theta)\not\equiv 0$ for some $\omega\in \Omega_{k,l}\cup \bigcup_{k'\le k}\tilde{\Omega}_{k',l}$ and $s\in \P_{l,x}$, then this implies that $\omega\cap\omega_{s,2}\not=\emptyset$. Moreover, when $\omega \in  \Omega_{k,l}$ this latter restriction together with  $|I_s|<2^{k}$ is equivalent with just asking that $\omega\subseteq \omega_{s,2}$. Similarly, when $\omega \in\bigcup_{k'\le k}\tilde{\Omega}_{k',l}$ then $\omega_{s,2}\subsetneq \omega$ is impossible, which in turn makes the requirement  $|I_s|<2^{k}$ superfluous. Indeed  $\omega_{s,2}\subsetneq \omega$ would imply that $\omega_{s}\subseteq \omega$, contradicting the fact that  $\omega_s$ contains an element from $\Xi_{x,l}$ while   $\omega$ does not. Hence we can rewrite 
\begin{align}
\sum_{s\in\P_{l,x} :|I_s|< 2^k} a_s \phi_{s}(x,\theta) \label{mktrees1}
&=\sum_{\omega \in \Omega_{k,l}} 1_{\omega}(\theta)
\sum_{s\in\P_{l,x}\atop{\omega\subseteq \omega_{s,2}}}a_s \phi_s(x,\theta)\\
\label{mktrees2}
&\quad +
\sum_{k'\le k} 
\sum_{\omega \in \tilde{\Omega}_{k',l}}
1_{\omega}(\theta) \sum_{s\in\P_{l,x}\atop{\omega\subseteq \omega_{s,2}}}
a_s \phi_s(x,\theta).
\end{align}

The  multiplier in  ~\eqref{mktrees2} can be written more conveniently as  
$$(1-\sum_{\tilde{\omega}\in \Omega_{k,l}}1_{\tilde{\omega}})\left(\sum_{k'} 
\sum_{\omega \in \tilde{\Omega}_{k',l}}
1_{\omega}(\theta) \sum_{s\in\P_{l,x},\;\omega \subseteq\omega_{s,2}}
a_s \phi_s(x,\theta)\right)=(1-\sum_{\tilde{\omega}\in \Omega_{k,l}}1_{\tilde{\omega}})\sum_{s\in\P_{l,x}} a_s \phi_{s}(x,\theta),$$
given the fact that $(\bigcup_{I\in\Omega_{k,l}} I)^{c}=\bigcup_{k'\le k}\bigcup_{I\in\tilde{\Omega}_{k',l}}I$ and $(\bigcup_{k'\le k}\tilde{\Omega}_{k',l})\bigcap(\bigcup_{k'> k}\tilde{\Omega}_{k',l})=\emptyset,$ modulo the endpoints of intervals.
The above multiplier operator is the composition of two operators. The first one is the identity minus  an operator for which Corollary ~\ref{labellhgt64} provides good bounds. The second one is associated with the multiplier $\sum_{s\in\P_{l,x}} a_s \phi_{s}(x,\theta)$ and hence its $L^2$ norm will equal
\begin{equation}
\label{pieceCarl}
\|\sum_{s\in\P_{l,x}} a_s \phi_{s}(x,\theta)\|_{L^{\infty}_\theta(\R)}.
\end{equation}
We will start by estimating ~\eqref{pieceCarl}, and note that this will implicitly provide a  proof of inequality ~\eqref{discterecarleson} (and thus of Carleson-Hunt's Theorem), along the lines of the argument in the next section. We leave the details to the interested reader.

Fix a $\theta$\footnote{It suffices to assume $\theta$ is not a dyadic point} and note that the main contribution to ~\eqref{pieceCarl} comes from a single tree. More precisely, let $s_{\theta}$ be a maximal element with respect to the ordering of tiles in the collection
$$\X=\{s\in\P_{l,x}:\theta\in\omega_{s,2}\}.$$If $\T\in\F_{l,x}$ (with top $T$) denotes one of the  trees such that $s_{\theta}\in G(\T)$, then  nestedness implies that $\X\subseteq G(\T).$
We also recognize as a consequence  of the definition of $\P_{l,x}$ that 
\begin{align}
\nonumber
|\sum_{s\in\P_{l,x}} a_s \phi_{s}(x,\theta)|&=
|\sum_{s\in G(\T)\cap\P_{l,x}} a_s \phi_{s}(x,\theta)|\\&
\label{dfg6208hgghta}
\le|\sum_{s\in \T_{l+1,0}^{(1)}\cap\P_{l,x}} a_s \phi_{s}(x,\theta)|
\\&
\label{dfg6208hgghtalater}
\quad+|\sum_{s\in \T_{l+1,0}^{(2)}\cap\P_{l,x}} a_s \phi_{s,\T_{l+1,0}^{(2)}}^{(l)}(x,\theta)|
\\&
\label{dfg6208hgghta2}
\quad+\sum_{s\in \P:|I_s|\le |I_{T}|\atop{\dist(x,I_s)\ge 2^{l-2}|I_{T}|}}| a_s \phi_{s}(x,\theta)|,
\end{align}
where $\T_{l+1,0}$ is the quasitree obtained from $\T$ by using the procedure in the beginning of the section, while $\T_{l+1,0}=\T^{(1)}_{l+1,0}\cup\T^{(2)}_{l+1,0}$ is the standard decomposition of $\T_{l+1,0}$.
The term ~\eqref{dfg6208hgghta2} is an error term and it is bounded crudely by $\sigma2^{-Ml},$ by using the triangle inequality, ~\eqref{e.m_s-size} and ~\eqref{BMNCARseq}.

We next focus on ~\eqref{dfg6208hgghta}. Note that  at  most one scale in  $\T^{(1)}_{l+1,0}$ contributes to the  summation ~\eqref{dfg6208hgghta}. Thus crude estimates relying on ~\eqref{BMNCARseq} prove that
\begin{align*}
|\sum_{s\in \T_{l+1,0}^{(1)}\cap\P_{l,x}} a_s \phi_{s}(x,\theta)|&\le \sup_{j\in\Z}\sum_{s\in \P:|I_s|=2^{j}\atop{x\notin 2^{l-1}I_s}} |a_s \phi_{s}(x,\theta)|\\&\lesssim \sigma2^{-Ml}.
\end{align*}

Before we evaluate the  sum corresponding to the 2-quasitree, we make two useful remarks. The first one concerns the fact that   there exists $n_{\T}$ depending on $x$ and $l$ such that 
$$\T^{(2)}_{l+1,0}\cap\P_{l,x}=\{s\in\T^{(2)}_{l+1,0}:2^{n_{\T}}\le |I_s|\}.$$
The lower bound on the scale is an immediate consequence of the definition of $\P_{l,x}.$
The second observation states that if $\phi_{s,\T_{l+1,0}^{(2)}}^{(l)}(x,\theta)\not=0$ for some $s\in \T_{l+1,0}^{(2)}$, then $|\theta-\xi_{\T_{l+1,0}}|\le |\omega_s|$.

 We then invoke inequalities ~\eqref{e.m_s-smoothgrr} and ~\eqref{BMNCARseq} to estimate
\begin{align*}
|\sum_{s\in\T^{(2)}_{l+1,0}\atop{|I_s|\ge 2^{n_{\T}}}}a_s \phi_{s,\T_{l+1,0}^{(2)}}^{(l)}(x,\theta)|&=|\sum_{s\in\T^{(2)}_{l+1,0}\atop{2^{n_{\T}}\le |I_s|\le|\theta-\xi_{\T_{l+1,0}}|^{-1}}}a_s \phi_{s,\T_{l+1,0}^{(2)}}^{(l)}(x,\theta)|\\&\le \sum_{s\in\T^{(2)}_{l+1,0}\atop{|I_s|\le|\theta-\xi_{\T_{l+1,0}}|^{-1}}}|a_s||\phi_{s,\T_{l+1,0}^{(2)}}^{(l)}(x,\theta)-\phi_{s,\T_{l+1,0}^{(2)}}^{(l)}(x,\xi_{\T_{l+1,0}})|\\
&\quad +|\sum_{s\in\T^{(2)}_{l+1,0}\atop{2^{n_{\T}}\le |I_s|\le|\theta-\xi_{\T_{l+1,0}}|^{-1}}}a_s \phi_{s,\T_{l+1,0}^{(2)}}^{(l)}(x,\xi_{\T_{l+1,0}})|\\&\lesssim \sigma2^{-Ml}+\|\sum_{s\in\T^{(2)}_{l+1,0}\atop{|I_s|<2^j}}a_s \phi_{s,\T_{l+1,0}^{(2)}}^{(l)}(x,\xi_{\T_{l+1,0}})\|_{V^r_j(\Z)}\\&\lesssim \sigma2^{-Ml}+\gamma 2^{-l},
\end{align*}
where in the last inequality we rely on the observation that  $x\notin E^{(3)}$.

We thus end up having the following estimate for ~\eqref{pieceCarl}
\begin{equation}
\label{laterneeded112}
\|\sum_{s\in\P_{l,x}} a_s \phi_{s}(x,\theta)\|_{L^{\infty}_\theta(\R)}\lesssim \sigma2^{-Ml}+\gamma2^{-l}.
\end{equation}
Finally, triangle inequality in $l$, an application of Corollary ~\ref{labellhgt64}  and the fact that $|\Xi_{x,l}|\le \beta 2^{2l}$ conclude that
\begin{equation}
\label{g55rdcchanged}
\left\|\left(\sum_{l\ge 0}(1-\sum_{\tilde{\omega}\in \Omega_{k,l}}1_{\tilde{\omega}})\sum_{s\in\P_{l,x}} a_s \phi_{s}(x,\theta)\right)_{k \in \Z}\right\|_{M_{2,\theta}^{*}(\R)}\lesssim {\beta}^{1/2-1/r}(\sigma+\gamma).
\end{equation}

We will next turn our attention to the term ~\eqref{mktrees1}. 
The multiplier in ~\eqref{mktrees1} is of the form $\sum_{\omega\in \Omega_{k,l}}1_{\omega}(\theta)m_{\omega}(x,\theta)$, where we define
$$m_{\omega}(x,\theta):=\sum_{s\in\P_{l,x}\atop{\omega_\subseteq \omega_{s,2}}}a_s \phi_s(x,\theta).$$

To estimate the norm $M_2^{*}$ of this sequence of multipliers  we will use  Theorem ~\ref{thm:gh098cfhh} with  $\Omega_{k,l}$ as the collection $R_k.$ Fix $l$ and consider a collection of nested intervals $\omega_k\in\Omega_{k,l}$, $k\in\Z$. For the remaining part of the proof we will be concerned with obtaining pointwise estimates in $x$ for the quantity
\begin{equation}
\label{jinjinjin}
\left\|\|(\widehat{g}1_{\omega_k}\sum_{s\in\P_{l,x}\atop{\omega_k\subseteq \omega_{s,2}}}a_s \phi_s(x,\cdot))\check{\ }(z)\|_{V^r_k(L)}\right\|_{L^2_z(\R)},
\end{equation}
which are uniform over all functions $g$ with $\|g\|_{L^2(\R)}= 1$,
where the inverse Fourier transform of the inner most expression is taken with respect to the variable $\theta$.

Fix $g$.
We observe that  the collection
$$\B:=\{\omega_{s,2}:s\in\P_{l,x},\; \omega_{k}\subseteq\omega_{s,2}\;\hbox{for some}\; k\in\Z\}$$  consists of nested intervals. Since this collection  is finite, it contains  a smallest element, corresponding to some $s_0\in \P_{l,x}$, such that $\omega_{s_0,2}\subseteq \omega_{s,2}$ whenever $\omega_{s,2}\in\B$. Now $s_0\in G(\T)$ for some $\T\in \F_{l,x}$ and hence all the tiles contributing to the term  ~\eqref{jinjinjin} are in $G(\T)$. 
 For each  $m\in\Z$, we denote by $\T_{l,m}$ the quasitree obtained from $\T$ by the procedure described in the beginning of the section.  Choose some arbitrary  $c\in \bigcap_{k}\omega_k$. This kind of choice for  $c$ will make possible the estimation of the two error terms below by rather trivial methods.

The next adjustment has to do with the fact that  $m_{\omega}$ is not constant. We will write it as the sum of a main (constant) term and two error terms $m_{\omega}=m_{\omega}^{(1)}+m_{\omega}^{(2)}+m_{\omega}^{(3)}$, with
\begin{align*}
m_{\omega_k}^{(1)}(x,\theta)&:=\sum_{s\in G(\T)\cap\P_{l,x}\atop{\omega_k\subset \omega_{s,2}}}a_s \phi_s(x,\xi_{\T})\\
m_{\omega_k}^{(2)}(x,\theta)&:=\sum_{s\in G(\T)\cap\P_{l,x}\atop{\omega_k\subset \omega_{s,2}}}a_s (\phi_s(x,\theta)-\phi_s(x,c)),\\ 
m_{\omega_k}^{(3)}(x,\theta)&:=\sum_{s\in G(\T)\cap\P_{l,x}\atop{\omega_k\subset \omega_{s,2}}}a_s(\phi_s(x,c)-\phi_s(x,\xi_{\T})).
\end{align*}
In dealing with the first error term we get the following sequence of inequalities, uniformly in  $x$
\begin{align}
\nonumber
\|\|(\widehat{g}1_{\omega_k}m_{\omega_k}^{(2)}(x,\cdot))\check{\ }(z)\|_{V^r_k(L)}\|_{L^2_z(\R)}&\le \|(\sum_k|(\widehat{g}1_{\omega_k}m_{\omega_k}^{(2)}(x,\cdot))\check{\ }(z)|^2)^{1/2}\|_{L^2_z(\R)}\\\nonumber&\le \sup_{\theta}(\sum_{k:|\omega_k|\ge |\theta-c|}|m_{\omega_k}^{(2)}(x,\theta)|^2)^{1/2}
\\\nonumber&\lesssim \sup_{\theta}(\sum_{k:|\omega_k|\ge |\theta-c|}(\sigma\sum_{s\in:G(\T)\cap\P_{l,x}\atop{|\omega|\ge |\omega_k|}}{|I_s|}^{1/2}|\phi_{s}(x,\theta)-\phi_{s}(x,c)|)^2)^{1/2}
\\\nonumber&\lesssim \sup_{\theta}(\sum_{k:|\omega_k|\ge |\theta-c|}(\sigma\sum_{s\in:G(\T)\cap\P_{l,x}\atop{|\omega|\ge |\omega_k|}}|c-\theta||I_s|\chi_{I_s}^{M}(x))^2)^{1/2}\\\nonumber&\lesssim \sup_{\theta}(\sum_{k:|\omega_k|\ge |\theta-c|}(\sigma|c-\theta||\omega_k|^{-1}2^{-Ml})^2)^{1/2}\\
\label{giro1qas1}
&\lesssim \sigma2^{-Ml}.
\end{align}
The passage from the first to the second line above is insured by the trivial inequality $\|\cdot\|_{V^r}\lesssim 2\|\cdot\|_{l^2}$, while the passage from the fourth line to the fifth relies on the estimate ~\eqref{xbfgr1456bb} on the  $\theta$ derivative of $\phi_{s}(x,\theta)$. The passage from the fifth line to the sixth line relies on the fact that $s\in\P_{l,x}$ implies  $x\notin 2^{l-1}I_s.$

To estimate the second error term we invoke  Lemma ~\ref{lem:auxq3}, Lemma ~\ref{lem:auxq4} and the fact that $\|\cdot\|_{V^r}\lesssim \|\cdot\|_{l^1}$
\begin{align}
\nonumber
\|\|(\widehat{g}1_{\omega_k}m_{\omega_k}^{(3)}(x,\cdot))\check{\ }(z)\|_{V^r_k(L)}\|_{L^2_z(\R)}&\lesssim\|m_{\omega_k}^{(3)}(x,0)\|_{V^r_k(L)} \|\|(\widehat{g}1_{\omega_k})\check{\ }(z)\|_{V^r_k(L)}\|_{L^2_z(\R)}
\\\nonumber&\lesssim\sum_{s\in G(\T)\cap\P_{l,x}\atop{\omega_k\subset \omega_{s,2}}}|a_s(\phi_s(x,c)-\phi_s(x,\xi_{\T}))|\\ \nonumber&\lesssim|c-\xi_{\T}| \sum_{2^n\le |c-\xi_{\T}|^{-1}}\sum_{|I_s|=2^n\atop{x\notin 2^{l-1}I_s}}\sigma|I_s|\chi_{I_s}^M(x)\\
\label{giro1qas2}&\lesssim \sigma2^{-Ml}.
\end{align}

The last task is to get estimates for the main term. We decompose each $m_{\omega_k}^{(1)}$ as 
$$m_{\omega_k}^{(1)}=\sum_{m\in\Z}m_{\omega_k,m}^{(1,1)}+\sum_{m\in\Z}m_{\omega_k,m}^{(1,2)},$$
where 
$$m_{\omega_k,m}^{(1,i)}(x,\theta)=\sum_{s\in\T_{l,m}^{(i)}\cap\P_{l,x}\atop{\omega_k\subseteq\omega_{s,2}}}a_s\phi_s(x,\xi_{\T}).$$ Then we estimate

\begin{align}
\label{-o-y6}
\|\|(\widehat{g}1_{\omega_k}m_{\omega_k}^{(1)}(x,\cdot))\check{\ }(z)\|_{V^r_k(L)}\|_{L^2_z(\R)}&\lesssim \sum_{m\in\Z} \|\|(\widehat{g}1_{\omega_k}m_{\omega_k,m}^{(1,1)}(x,\cdot))\check{\ }(z)\|_{V^r_k(L)}\|_{L^2_z(\R)}+\\&
\label{-o-y67}
+\sum_{m\in\Z} \|\|(\widehat{g}1_{\omega_k}m_{\omega_k,m}^{(1,2)}(x,\cdot))\check{\ }(z)\|_{V^r_k(L)}\|_{L^2_z(\R)}.
\end{align}

In analyzing the term ~\eqref{-o-y6} we note that for each $\T$, $l$ and $m$ the collection 
$$\U=\bigcup_{k}\{s\in \T_{l,m}^{(1)}\cap\P_{l,x}:\omega_k\subseteq \omega_{s,2}\}$$ contains at most one scale. This is because  the collections $\{\omega_{s,2}:s\in \U\}$ and  $\{\omega_{s,1}:s\in \U\}$ are nested. Also, if $l\ge 1$ and $s\in \U$ then $x\notin 2^{\alpha(l,m)-1}I_s$. Thus, for each $m\in\Z$ and each $l\ge 1$, Lemma ~\ref{lem:auxq4} gives 
\begin{equation}\label{nnggffccxsd2}
\|\|(\widehat{g}1_{\omega_k}m_{\omega_k,m}^{(1,1)}(x,\cdot))\check{\ }(z)\|_{V^r_k(L)}\|_{L^2_z(\R)}\le \sup_{j\in\Z}\sum_{s\in\P:|I_s|=2^j\atop{x\notin 2^{\alpha(l,m)-1}I_s}}|a_s||\phi_s(x,\xi_{\T})|\lesssim_M \sigma2^{-M\alpha(l,m)},
\end{equation}
and by a similar argument, the same works for $l=0$, too.

Next, we consider the term ~\eqref{-o-y67}. 
We first  acknowledge the fact that for each $k,l,m$ there exist $n_3\le n_4$ such that
$$\{s\in \T_{l,m}^{(2)}\cap\P_{l,x}:\omega_k\subseteq \omega_{s,2}\}=\{s\in \T_{l,m}^{(2)}:2^{n_3}\le |I_s|<2^{n_4}\}.$$ The number $n_3$ is independent of $k$ and appears as a restriction due to the fact that at level $l$ we only consider tiles that have not been selected at previous stages. The  restriction $|I_s|<2^{n_4}$ replaces the restriction  $\omega_k\subseteq \omega_{s,2}$ and $n_4$ is increasing as a function of $k$. This observation together with  Lemma ~\ref{lem:auxq3}, Lemma ~\ref{lem:auxq4}  and the fact that $x\notin E^{(2)}$ implies

\begin{align}
\nonumber
\|\|(\widehat{g}1_{\omega_k}m_{\omega_k,m}^{(1,2)}(x,\cdot))\check{\ }(z)\|_{V^r_k(L)}\|_{L^2_z(\R)}&\le \|\sum_{s\in\T_{l,m}^{(2)}\cap\P_{l,x}\atop{\omega_k\subset\omega_{s,2}}}a_s\phi_s(x,\xi_{\T})\|_{V^r_k(L)} \|\|(\widehat{g}1_{\omega_k})\check{\ }(z)\|_{V^r_k(L)}\|_{L^2_z(\R)}\\\nonumber&\lesssim \|\sum_{s\in\T_{l,m}^{(2)}\atop{|I_s|<2^j}}a_s\phi_{s,\T_{l,m}^{(2)}}^{(\alpha(l,m))}(x,\xi_{\T})\|_{V^r_j(L)}\\
\label{nnggffccxsd1}
&\lesssim  \gamma 2^{-Ml}(|m|+1)^{-2}.
\end{align}

Thus, summation over $m$ in inequalities  ~\eqref{nnggffccxsd2} and ~\eqref{nnggffccxsd1} leads to 
\begin{equation}
\label{giro1qas3}
\|\|(\widehat{g}1_{\omega_k}m_{\omega_k}^{(1)}(x,\cdot))\check{\ }(z)\|_{V^r_k(L)}\|_{L^2_z(\R)}\lesssim (\sigma+\gamma)2^{-l}.
\end{equation}

A final application of the triangle inequality with respect to  $l$ in ~\eqref{giro1qas1}, ~\eqref{giro1qas2} and ~\eqref{giro1qas3}, together with Theorem ~\ref{thm:gh098cfhh} and the fact that $|\Xi_{x,l}|\le \beta 2^{2l}$ conclude to 

\begin{equation}
\label{g55rdciiii}
\|( \sum_{l\ge 0} 
\sum_{\omega \in {\Omega}_{k,l}}
1_{\omega}(\theta) \sum_{s\in\P_{l,x}\atop{\omega \subseteq\omega_{s,2}}}a_s\phi_s(x,\theta))_{k \in \Z} \|_{M_{2,\theta}^{*}(\R)}\lesssim {\beta}^{1/2-1/r}(\sigma+\gamma).
\end{equation}

By putting together the estimates from ~\eqref{g55rdcchanged} and ~\eqref{g55rdciiii} the conclusion of our theorem follows.
\end{proof}

We continue with the variant of Theorem ~\ref{mainineq} that will prove useful in the proof of the oscillation inequality in Theorem ~\ref{ceamaiceahho1}. To this end, let ${\bf U}:=(u_j)_{j=1}^{J}$  be a finite sequence of integers and recall the oscillation-variational norm $\|\cdot\|_{O_{{\bf U}} \cap V^r}$ introduced in ~\eqref{defoscvar}.

\begin{theorem}
\label{mainineqo1}
Assume we are in the settings preceding Theorem ~\ref{mainineq} and also that the following additional requirement is satisfied
\begin{equation*}
\label{BMNCARseqo1}
\sup_{s\in\P}\frac{|a_s|}{|I_s|^{1/2}}\le \sigma.
\end{equation*}
Define the  exceptional sets
\begin{align*}
E^{(1)}&=\bigcup_{l\ge 0}\{x:\sum_{\T\in\F}1_{2^lI_T}(x)>\beta2^{2l}\},\\
E^{(2)}&=\bigcup_{l,m\ge 0}\bigcup_{\T\in\F_{l,m}^{(2)}}\{x:\|\sum_{s\in \T\atop{|I_s|<2^j}}a_s\phi_{s,\T}^{(\alpha(l,m))}(x,\xi_{\T})\|_{O_{{\bf U}} \cap V^r_j(\Z)}>\gamma2^{-l}(|m|+1)^{-2}\},\\
E^{(3)}&= \bigcup_{l\ge 0}\bigcup_{\T\in\F_{l+1,0}^{(2)}}\{x:\|\sum_{s\in \T\atop{|I_s|<2^j}}a_s\phi_{s,\T}^{(l)}(x,\xi_{\T})\|_{O_{{\bf U}} \cap V^r_j(\Z)}>\gamma2^{-l}\},
\end{align*}
where the symbol $\alpha(l,m)$ equals $l$ if $m=0$ and $l+[\log_2|m|]$ if $m\not=0$.
Then  for each $x\notin E^{(1)}\cup E^{(2)}\cup E^{(3)}$ and each $g$ with $\|g\|_{L^2(\R)}=1$ we have the uniform pointwise estimate
$$
(\sum_{j=1}^{J-1}\|\sup_{\;u_j\le
k<u_{j+1}}|\F^{-1}_{\theta}\{\sum_{s\in \S\atop{2^{u_j}\le
|I_s|<2^k}}a_s\phi_s(x,\theta)\widehat{g}(\theta)\}(z)|\|_{L^2_z(\R)}^2)^{1/2}\lesssim J^{\frac{r+4}{4r+4}}\beta^{1-2/r}(\gamma+\sigma),$$
with the implicit constants depending only on $r$.
\end{theorem}

\begin{proof}
The proof follows closely the lines of the proof of Theorem ~\ref{mainineq}. 
Fix  $x\notin E^{(1)}\cup E^{(2)}\cup E^{(3)}$ and fix $g$ with $\|g\|_{L^2(\R)}=1$.
We will use the notation introduced in the beginning of the proof of Theorem ~\ref{mainineq} and the representation 
\begin{align}
\label{mktrees1o}
\sum_{s\in\P_{l,x} :|I_s|< 2^k} a_s \phi_{s}(x,\theta) &=\sum_{\omega \in \Omega_{k,l}} 1_{\omega}(\theta)
\sum_{s\in\P_{l,x}\atop{\omega\subseteq \omega_{s,2}}}a_s \phi_s(x,\theta)\\
\label{mktrees2o}
&\quad +
(1-\sum_{\tilde{\omega}\in \Omega_{k,l}}1_{\tilde{\omega}})\sum_{s\in\P_{l,x}} a_s \phi_{s}(x,\theta).
\end{align}
Then, by using the triangle inequality in $l$, Corollary ~\ref{labellhgt64o},  inequality ~\eqref{laterneeded112} and the fact that $x\notin E^{(1)}$ we get the following estimate for contribution to the term ~\eqref{mktrees2o}: if $\|g\|_{L^2(\R)} = 1$, then
\begin{align*}
&\left(\sum_{j=1}^{J-1} 
\left\|\sup_{\;u_j\le k<u_{j+1}} |\sum_{l\ge 0}\F^{-1}_{\theta}
\{[(1-\sum_{\tilde{\omega}\in \Omega_{k,l}}1_{\tilde{\omega}})-
(1-\sum_{\tilde{\omega}\in \Omega_{u_j,l}}1_{\tilde{\omega}})]
 \sum_{s\in\P_{l,x}} a_s \phi_{s}(x,\theta)\widehat{g}(\theta)\}(z)|
\right\|_{L^2_z(\R)}^2\right)^{1/2}
\\&\lesssim J^{\frac{r+4}{4r+4}}{\beta}^{1-2/r}\sum_{l\ge 0}\|\F^{-1}_{\theta}\{\sum_{s\in\P_{l,x}} a_s \phi_{s}(x,\theta)\widehat{g}(\theta)\}(z)\|_{L^2_z(\R)}\\\label{carestpartdoi}&\lesssim J^{\frac{r+4}{4r+4}}{\beta}^{1-2/r}(\sigma+\gamma).
\end{align*}

The multiplier in ~\eqref{mktrees1o} is of the form $\sum_{\omega\in \Omega_{k,l}}1_{\omega}(\theta)m_{\omega}(x,\theta)$, where we define as before
$$m_{\omega}(x,\theta)=\sum_{s\in\P_{l,x}\atop{\omega_\subseteq \omega_{s,2}}}a_s \phi_s(x,\theta).$$

To estimate 
$$\left(\sum_{j=1}^{J-1}\left\|\sup_{\;u_j\le
k<u_{j+1}}|\F^{-1}_{\theta}\{\sum_{\omega\in \Omega_{k,l}}1_{\omega}(\theta)m_{\omega}(x,\theta)\widehat{g}(\theta)-\sum_{\omega\in \Omega_{u_j,l}}1_{\omega}(\theta)m_{\omega}(x,\theta)\widehat{g}(\theta)\}(z)|\right\|_{L^2_z(\R)}^2\right)^{1/2}$$
we will use  Theorem ~\ref{Bourgeneralo} with  $\Omega_{k,l}$ as the collection $R_k.$ Fix $l\ge 0$ and consider a collection of nested intervals $\omega_k\in\Omega_{k,l}$, $k\in\Z$. For the remaining part of the proof  we will be concerned with obtaining pointwise estimates in $x$ for the quantity
\begin{equation}
\label{jinjinjin1}
\left\|\|(\widehat{g}1_{\omega_k}\sum_{s\in\P_{l,x}\atop{\omega_k\subseteq \omega_{s,2}}}a_s \phi_s(x,\cdot))\check{\ }(z)\|_{O_{{\bf U}} \cap V^r_k(L)}\right\|_{L^2_z(\R)},
\end{equation}
where the inverse Fourier transform of the inner most expression is taken with respect to the variable $\theta$. Also, given the estimates for the $V^r$ norm from the proof of the previous theorem, all that is left is getting the corresponding oscillation estimates. Split as before $m_{\omega}=m_{\omega}^{(1)}+m_{\omega}^{(2)}+m_{\omega}^{(3)}$. 

The same type of estimates as in ~\eqref{giro1qas1} lead to the following estimate for the error term associated with the multiplier $m_{\omega}^{(2)}$

\begin{equation}
\label{giro1qas1o}
\left\|\|(\widehat{g}1_{\omega_k}m_{\omega_k}^{(2)}(x,\cdot))\check{\ }(z)\|_{O_{\bf U}(k)}\right\|_{L^2_z(\R)}
\lesssim \sigma2^{-Ml}.
\end{equation}

To estimate the second error term associated with the multiplier $m_{\omega}^{(3)}$ we proceed like in ~\eqref {giro1qas2}. By invoking the second part of  Lemma ~\ref{lem:auxq3}, Lemma ~\ref{lem:auxq4} and the fact that $\|\cdot\|_{O_{{\bf U}}}\lesssim \|\cdot\|_{l^1}$ we get 
\begin{equation}
\label{giro1qas2o}
\left\|\|(\widehat{g}1_{\omega_k}m_{\omega_k}^{(3)}(x,\cdot))\check{\ }(z)\|_{O_{\bf U}}\right\|_{L^2_z(\R)}
\lesssim \sigma2^{-Ml}.
\end{equation}

The last task is to get estimates for the main term. We split as before  
$$m_{\omega_k}^{(1)}=\sum_{m\in\Z}m_{\omega_k,m}^{(1,1)}+\sum_{m\in\Z}m_{\omega_k,m}^{(1,2)},$$
and  estimate

\begin{align}
 \label{-o-y6o}
\left\|\|(\widehat{g}1_{\omega_k}m_{\omega_k}^{(1)}(x,\cdot))\check{\ }(z)\|_{O_{{\bf U}}}\right\|_{L^2_z(\R)}&\lesssim \sum_{m\in\Z}\left\|\|(\widehat{g}1_{\omega_k}m_{\omega_k,m}^{(1,1)}(x,\cdot))\check{\ }(z)\|_{O_{{\bf U}}}\right\|_{L^2_z(\R)}+\\&
\label{-o-y670}
+\sum_{m\in\Z}\left\|\|(\widehat{g}1_{\omega_k}m_{\omega_k,m}^{(1,2)}(x,\cdot))\check{\ }(z)\|_{O_{{\bf U}}}\right\|_{L^2_z(\R)}.
\end{align}

The same discussion as the one regarding the derivation of inequalities ~\eqref{nnggffccxsd2} and ~\eqref{nnggffccxsd1} shows that 

\begin{equation}
\label{nnggffccxsd2o}
\left\|\|(\widehat{g}1_{\omega_k}m_{\omega_k,m}^{(1,1)}(x,\cdot))\check{\ }(z)\|_{O_{\bf U}}\right\|_{L^2_z(\R)}\lesssim_M \sigma2^{-M\alpha(l,m)},
\end{equation}

\begin{equation}
\label{nnggffccxsd1o}
\left\|\|(\widehat{g}1_{\omega_k}m_{\omega_k,m}^{(1,2)}(x,\cdot))\check{\ }(z)\|_{O_{{\bf U}}}\right\|_{L^2_z(\R)}\lesssim  \gamma 2^{-Ml}(|m|+1)^{-2}.
\end{equation}

Thus, summation over $m$ in inequalities  ~\eqref{nnggffccxsd2o} and ~\eqref{nnggffccxsd1o} leads to 
\begin{equation}
\label{giro1qas3o}
\left\|\|(\widehat{g}1_{\omega_k}m_{\omega_k}^{(1)}(x,\cdot))\check{\ }(z)\|_{O_{{\bf U}}}\right\|_{L^2_z(\R)}\lesssim (\sigma+\gamma)2^{-l}.
\end{equation}

A final application of the triangle inequality with respect to  $l$ in ~\eqref{giro1qas1o}, ~\eqref{giro1qas2o} and ~\eqref{giro1qas3o}, together with Theorem ~\ref{Bourgeneralo} and the fact that $|\Xi_{x,l}|\le \beta 2^{2l}$ concludes to 

\begin{equation*}
(\sum_{j=1}^{J-1}\left\|\sup_{\;u_j\le
k<u_{j+1}}|\sum_{l\ge 0}\F^{-1}_{\theta}\{\sum_{\omega\in \Omega_{k,l}}1_{\omega}(\theta)m_{\omega}(x,\theta)\widehat{g}(\theta)-\sum_{\omega\in \Omega_{u_j,l}}1_{\omega}(\theta)m_{\omega}(x,\theta)\widehat{g}(\theta)\}|\right\|_{L^2_x(\R)}^2)^{1/2}
\end{equation*}
\begin{equation}
\label{g55rdciiiio}
\lesssim J^{\frac{r+4}{4r+4}}{\beta}^{1-2/r}(\sigma+\gamma).
\end{equation}

By putting together the estimates from ~\eqref{carestpartdoi} and ~\eqref{g55rdciiiio}, the conclusion of our theorem follows.
\end{proof}

We close this section with a square function estimate for the Carleson-Hunt operator that will play the decisive role in the proof of Theorem ~\ref{ceamaiceahho2}. The proof does not contain any serious new ideas, other than the ones used in the proof of Theorem ~\ref{mainineq} to estimate the $L^{\infty}$ norm of the Carleson-Hunt operator.

\begin{theorem}
\label{mainineqo2}
Assume we are in the settings preceding Theorem ~\ref{mainineq} and also that the following additional requirement is satisfied
\begin{equation*}
\label{BMNCARseqo2}
\sup_{s\in\P}\frac{|a_s|}{|I_s|^{1/2}}\le \sigma.
\end{equation*}
Define the  exceptional set
$$E:= \bigcup_{l\ge 0}\bigcup_{\T\in\F_{l+1,0}^{(2)}}\{x:(\sum_{j\in\Z}|\sum_{s\in \T\atop{|I_s|=2^j}}a_s\phi_{s,\T}^{(l)}(x,\xi_{\T})|^2)^{1/2}>\gamma2^{-l}\}.$$
Then  for each $x\notin E$ we have 
\begin{equation}\label{forex-2}
\left\|(\sum_{k\in\Z}|\sum_{s\in \P\atop{|I_s|=2^k}}a_s\phi_s(x,\theta)|^2)^{1/2}\right\|_{M_{2,\theta}(\R)}\lesssim \sigma+\gamma.
\end{equation}
\end{theorem}
\begin{proof}

We will assume again the notation introduced in the beginning of the proof of Theorem ~\ref{mainineq}.
Fix $x\notin E$ and  $\theta\in\R$. Note that the main contribution to \eqref{forex-2} comes from a single tree. More precisely, let $s_{\theta}$ be a maximal element in the collection
$$\X=\{s\in\P_{l,x}:\theta\in\omega_{s,2}\}.$$If $\T\in\F_{l,x}$ (with top $T$) denotes one of the  trees such that $s_{\theta}\in G(\T)$, then  nestedness implies that $\X\subseteq G(\T).$
We also recognize as a consequence  of the definition of $\P_{l,x}$ that 
\begin{align}
\nonumber
(\sum_{k\in\Z}|\sum_{s\in\P_{l,x}\atop{|I_s|=2^k}} a_s \phi_{s}(x,\theta)|^2)^{1/2}&=
(\sum_{k\in\Z}|\sum_{s\in G(\T)\cap\P_{l,x}\atop{|I_s|=2^k}} a_s \phi_{s}(x,\theta)|^2)^{1/2}\\&
\label{dfg6208hgghtao2}
\lesssim(\sum_{k\in\Z}|\sum_{s\in \T_{l+1,0}^{(1)}\cap\P_{l,x}\atop{|I_s|=2^k}} a_s \phi_{s}(x,\theta)|^2)^{1/2}
\\&
\label{dfg6208hgghtao2jherf}
\quad +(\sum_{k\in\Z}|\sum_{s\in \T_{l+1,0}^{(2)}\cap\P_{l,x}\atop{|I_s|=2^k}} a_s \phi_{s,\T_{l+1,0}^{(2)}}^{(l)}(x,\theta)|^2)^{1/2}
\\&
\label{dfg6208hgghta2o2}
\quad +(\sum_{k\in\Z}(\sum_{s\in \P:|I_s|\le |I_{T}|\atop{\dist(x,I_s)\ge 2^{l-2}|I_{T}|\atop{|I_s|=2^k}}}| a_s \phi_{s}(x,\theta)|)^2)^{1/2},
\end{align}

The term ~\eqref{dfg6208hgghta2o2} is bounded in the same manner as the term ~\eqref{dfg6208hgghta2}  by $\sigma2^{-Ml}$ (uniformly in $\theta$).

We next focus on ~\eqref{dfg6208hgghtao2}. Note that  at  most one scale in  $\T^{(1)}_{l+1,0}$ contributes to the  summation ~\eqref{dfg6208hgghtao2}. Thus estimates like the ones for ~\eqref{dfg6208hgghta} prove that
$$
(\sum_{k\in\Z}|\sum_{s\in \T_{l+1,0}^{(1)}\cap\P_{l,x}\atop{|I_s|=2^k}} a_s \phi_{s}(x,\theta)|^2)^{1/2}\le \sup_{j\in\Z}\sum_{s\in \P:|I_s|=2^{j}\atop{x\notin 2^{l-1}I_s}} |a_s \phi_{s}(x,\theta)|\lesssim \sigma2^{-Ml}.
$$

Before we evaluate the  sum corresponding to the 2-quasitree, we recall from the proof of Theorem ~\ref{mainineq} that  there exists $n_{\T}$ depending on $x$ and $l$ such that 
$$\T^{(2)}_{l+1,0}\cap\P_{l,x}=\{s\in\T^{(2)}_{l+1,0}:2^{n_{\T}}\le |I_s|\}.$$
Also, recall that if $\phi_{s,\T_{l+1,0}^{(2)}}^{(l)}(x,\theta)\not=0$ for some $s\in \T_{l+1,0}^{(2)}$, then $|\theta-\xi_{\T_{l+1,0}}|\le |\omega_s|$.

 We then estimate
$$
(\sum_{k\in\Z}|\sum_{s\in\T^{(2)}_{l+1,0}\cap\P_{l,x}\atop{|I_s|=2^k}}a_s \phi_{s,\T_{l+1,0}^{(2)}}^{(l)}(x,\theta)|^2)^{1/2}=(\sum_{2^{n_{\T}}\le 2^k\le\atop{\le |\theta-\xi_{\T_{l+1,0}}|^{-1}}}|\sum_{s\in\T^{(2)}_{l+1,0}\atop{|I_s|=2^k}}a_s \phi_{s,\T_{l+1,0}^{(2)}}^{(l)}(x,\theta)|^2)^{1/2}$$
$$\le (\sum_{ 2^k\le|\theta-\xi_{\T_{l+1,0}}|^{-1}}(\sum_{s\in\T^{(2)}_{l+1,0}\atop{|I_s|=2^k}}|a_s||\phi_{s,\T_{l+1,0}^{(2)}}^{(l)}(x,\theta)-\phi_{s,\T_{l+1,0}^{(2)}}^{(l)}(x,\xi_{\T_{l+1,0}})|)^2)^{1/2}$$
$$
+(\sum_{k\in\Z}|\sum_{s\in\T^{(2)}_{l+1,0}\atop{|I_s|=2^k}}a_s \phi_{s,\T_{l+1,0}^{(2)}}^{(l)}(x,\xi_{\T_{l+1,0}})|^2)^{1/2}\lesssim \sigma2^{-Ml}+\gamma 2^{-l},$$
where in the last inequality we have used the fact that $x\notin E$.

We thus end up having the following estimate 
\begin{equation}
\label{laterneeded112square}
\|(\sum_{k\in\Z}|\sum_{s\in \P_{l,x}\atop{|I_s|=2^k}}a_s\phi_s(x,\theta)|^2)^{1/2}\|_{M_{2,\theta}(\R)}
\lesssim \sigma2^{-Ml}+\gamma2^{-l}.
\end{equation}
Finally, the triangle inequality in $l$ concludes that
\begin{align*}
\|(\sum_{k\in\Z}|\sum_{s\in \P\atop{|I_s|=2^k}}a_s\phi_s(x,\theta)|^2)^{1/2}\|_{M_{2,\theta}(\R)}&\le\sum_{l\ge 0}\|(\sum_{k\in\Z}|\sum_{s\in \P_{l,x}\atop{|I_s|=2^k}}a_s\phi_s(x,\theta)|^2)^{1/2}\|_{M_{2,\theta}(\R)}\\&\lesssim \sigma+\gamma.
\end{align*}
\end{proof}

\section{Proof of Theorems ~\ref{ceamaiceahh}, ~\ref{ceamaiceahho1} and ~\ref{ceamaiceahho2}}
\label{section:last}

We will present the proof  of  Theorem ~\ref{ceamaiceahh} in detail and then indicate the modifications that have to be made in the argument to get Theorems ~\ref{ceamaiceahho1} and ~\ref{ceamaiceahho2}.
Let ${\bf U}=(u_j)_{j=1}^{J}$ be an arbitrary finite sequence of integers.
For each collection of tiles $\S'\subseteq \S$ define the following operators, relevant for the three theorems we aim to prove:
\begin{align*}
T_{\S'}f(x)&:=\|(\sum_{s\in \S'\atop{|I_s|< 2^k}}\langle f,\varphi_s\rangle\phi_s(x,\theta))_{k\in\Z}\|_{M_{2,\theta}^*(\R)}. \\
O_{\S'}f(x)&:=\sup_{\|g\|_{L^2(\R)}=1}(\sum_{j=1}^{J-1}\|\sup_{\;u_j\le
k<u_{j+1}}|\F^{-1}_{\theta}\{\sum_{s\in \S'\atop{2^{u_j}\le
|I_s|<2^k}}\langle f,\varphi_s\rangle\phi_s(x,\theta)\widehat{g}(\theta)\}(z)|\|_{L^2_z(\R)}^2)^{1/2}\\
Q_{\S'}f(x)&:=\|(\sum_{k\in\Z}|\sum_{s\in \S'\atop{|I_s|=2^k}}\langle
f,\varphi_s\rangle\phi_s(x,\theta)|^2)^{1/2}\|_{M_{2,\theta}(\R)}.
\end{align*}
Let $V$ denote any of these operators. Define  $c(V,p)$ to equal 1 if $V$ is either $T$ or $Q$, and $J^{1/2-\delta(p)}$ when $V=O$. For each $1<p<\infty$, the index $\delta(p)$ is some number in  $(0,\frac12)$ whose value will become implicit later in the argument, without however being computed explicitly.  

There is a common part in the argument  for all  three operators above, and we will describe it in the following.
Note that for each $\S'$ the operator $V_{\S'}$ is sublinear as a function of $f$. Also, for each $f$ and $x$ the mapping  $\S'\to V_{\S'}f(x)$ is sublinear as a function of the tile set $\S'$. We will prove in the following that
\begin{equation}
\label{en.6}
m\{x:V_{\S}1_{F}(x)\gtrsim \lambda\}\lesssim c(V,p)\frac{|F|}{\lambda^p},
\end{equation}
for each $F\subseteq \R$ of finite measure, each $\lambda>0$ and each $1<p<\infty.$  Then, by invoking the Marcinkiewicz interpolation theorem and restricted weak type interpolation we get for each $1<p<\infty$ that
$$\|V_{\S}f\|_{p}\lesssim d(V,p)\|f\|_p,$$
where $d(V,p)$  equals 1 if $V$ is either $T$ or $Q$, and $J^{1/2-\epsilon(p)}$ when $V=O$, for some appropriate $\epsilon(p)\in (0,\frac12)$ whose value will become implicit later.

Fix $F$ and $\lambda$. 
We first prove ~\eqref{en.6} in the case $\lambda\le 1$. 
Define the  first exceptional set
$$E:=\{x:M_{p}1_{F}(x)\ge \lambda\}$$
and note that $|E|\lesssim \frac{|F|}{\lambda^p}.$
Split $\S=\S_1\cup\S_2$ where
\begin{align*}
\S_1&:=\{s\in\S:I_s\cap E^{c}\not=\emptyset\}\\
\S_2&:=\{s\in\S:I_s\cap E^{c}=\emptyset\}.
\end{align*}
Decompose $E=\bigcup_{i}E_i$ as a disjoint union of intervals $E_i$ and define $E':=\bigcup_i 2E_i.$
Let us  first show that 
\begin{equation}
m\{x\in (E')^c: V_{\S_2}1_F(x)\gtrsim \lambda\}\lesssim \frac{|F|}{\lambda^p}.
\end{equation}
Split $\S_2=\bigcup_{\kappa\ge 1}\S_2^{\kappa}$ where
$$\S_2^{\kappa}:=\{s\in\S_2:2^{\kappa-1}I_s\cap E^{c}=\emptyset,\;2^{\kappa}I_s\cap E^{c}\not=\emptyset\}.$$
Note further that if $s\in\S_2^{\kappa}$ then 
$$
\frac{|\langle 1_F,\varphi_s\rangle|}{\sqrt{|I_s|}}\lesssim \inf_{x\in I_s}\M_p1_F(x)\lesssim 2^{\kappa}\inf_{x\in 2^{\kappa}I_s}\M_p1_F(x)\lesssim\lambda2^{\kappa}.
$$
We next partition $\S_2^\kappa := \bigcup_{i,l \in \Z} \S_2^{\kappa,i,l}$, where
$\S_2^{\kappa,i,l}:=\{s\in\S_2^{\kappa}: I_s\subseteq E_i,\;|I_s|=2^l\}$ and  observe that for each $x\in (E')^c$
$$
V_{\S_2^{\kappa}}1_F(x)\le\sum_{i}\sum_{l}V_{\S_2^{\kappa,i,l}}1_F(x)\le\sum_{i}\sum_{l} \|\sum_{s\in \S_2^{\kappa,i,l}}\langle 1_F,\varphi_s\rangle\phi_s(x,\theta)\|_{L^{\infty}_\theta(\R)}$$
$$\lesssim \sum_{i}\sum_{l}\sum_{s\in\S\atop{2^\kappa I\subseteq E_i\atop{|I|=2^l}}}\lambda2^{\kappa}\chi_{I}^2(x)\lesssim\sum_{i}\sum_{l}\sum_{s\in\S\atop{2^{\kappa}I\subseteq E_i\atop{|I|=2^l}}}\lambda2^{\kappa}\chi_{E_i}^{2}(x)\left(\frac{2^l}{|E_i|}\right)^M\lesssim \sum_i\lambda2^{-M\kappa}\sum_{i}\chi_{E_i}^{2}(x).
$$
We next apply  the Fefferman-Stein inequality \cite{FS}
$$
\|V_{\S_2^{\kappa}}1_{F}\|_{L^p((E')^c)}\lesssim \lambda2^{-M{\kappa}}\|\sum_{i}\chi_{E_i}^{2}\|_{L^p(\R)}\lesssim \lambda2^{-M\kappa}|E|^{1/p}.
$$
Now we can write 
$$
m\{x\in (E')^c: V_{\S_2}1_F(x)\gtrsim \lambda\}\lesssim \lambda^{-p}(\sum_{\kappa\ge 1}\|V_{\S_2^{\kappa}f}\|_{L^p((E')^c)})^p\lesssim |E|\lesssim \frac{|F|}{\lambda^p}.
$$
The rest of the proof in the case $\lambda\le 1$ is devoted to arguing that 
\begin{equation}
\label{en.11}
m\{x\in \R: V_{\S_1}1_F(x)\gtrsim \lambda^{1-\epsilon}\}\lesssim c(V,p)\frac{|F|}{\lambda^p},
\end{equation}
for each $\epsilon>0$, with the implicit constant depending only on $\epsilon$ and $p$. By combining this with the previous estimates, we get  that for each $1<p<\infty$, $0<\lambda\le 1$, $\epsilon>0$ and each $F\subset \R$ of finite measure, there is an exceptional set of measure $\lesssim  \lambda^{-p}|F|$ such that for $x$ outside this set 
$$V_{\bf S}1_{F}(x)\lesssim c(V,p)\lambda^{1-\epsilon}.$$
Finally, this will easily imply ~\eqref{en.6} for $\lambda\le 1$, since the range of $p$ is open. The proofs of ~\eqref{en.11} and ~\eqref{en.6} in the case $\lambda>1$ for the operators $T_{\S}$, $O_{\S}$ and $Q_{\S}$ are very similar, the only difference appears in the choice of the exceptional set. We start by giving the full details for the operator $T_{\S}$, and then briefly indicate the modifications needed for the other two operators.

\subsection{The estimates for $T_{\S}$.}

We start by proving ~\eqref{en.11}.
Proposition ~\ref{g11a00uugg55} guarantees that $\size(\S_1)\lesssim \lambda$, where the size is understood here with respect to the  function $1_F$. Define $\Delta:=[-\log_2(\size(\S_1))].$ Use the result of Proposition ~\ref{Besselsineq} to split $\S_1$ as a disjoint union $\S_1=\bigcup_{n\ge \Delta}\P_n,$ where $\size(\P_n)\le 2^{-n}$ and each $\P_n$ consists of a family $\F_{\P_n}$ of trees  satisfying
\begin{equation}
\label{e:intp1}
\sum_{\T\in\F_{\P_n}}|I_T|\lesssim 2^{2n}|F|.
\end{equation}  

Let $\epsilon>0$ be an arbitrary  positive number.
For each $n\ge \Delta$ define $\sigma:=2^{-n}$, $\beta:=2^{3n}\lambda^{p}$, $\gamma:=2^{-n/2}\lambda^{1/2-\epsilon}$. Define $a_s:=\langle 1_F,\varphi_s\rangle$ for each $s\in\P_n$ and note that the collection $\P_n$ together with the coefficients $(a_s)_{s\in\P_n}$ satisfy the requirements of Theorem ~\ref{mainineq}. Let $\F_{\P_n,l,m}^{(2)}$ be the collection of all the 2-quasitrees $\T_{l,m}^{(2)}$ obtained from all the trees $\T\in\F_{\P_n}$ by the procedure described in the beginning of the previous section. Define the corresponding exceptional sets
\begin{align*}
E^{(1)}_{n}&:=\bigcup_{l\ge 0}\{x:\sum_{\T\in\F_{\P_n}}1_{2^lI_T}(x)>\beta2^{2l}\},\\
E^{(2)}_{n}&:=\bigcup_{l,m\ge 0}\bigcup_{\T\in\F_{\P_n,l,m}^{(2)}}\{x:\|\sum_{s\in \T\atop{|I_s|<2^j}}a_s\phi_{s,\T}^{(\alpha(l,m))}(x,\xi_{\T})\|_{V^r_j(\Z)}>\gamma2^{-l}(|m|+1)^{-2}\},\\
E^{(3)}_{n}&:= \bigcup_{l\ge 0}\bigcup_{\T\in\F_{\P_n,l+1,0}^{(2)}}\{x:\|\sum_{s\in \T\atop{|I_s|<2^j}}a_s\phi_{s,\T}^{(l)}(x,\xi_{\T})\|_{V^r_j(\Z)}>\gamma2^{-l}\}.
\end{align*}
By  ~\eqref{e:intp1} and the fact that $\lambda\le 1$ we get
$$|E^{(1)}_{n}|\lesssim 2^{-n}\lambda^{-p}|F|.$$
By Theorem ~\ref{BMOtree} and the fact that $\lambda\le 1$, for each $1<s<\infty$ we get
$$|E^{(2)}_{n}|\lesssim \gamma^{-s}\sigma^{s-2}|F|\lesssim 2^{-n(s/2-2)}\lambda^{-s(1/2-\epsilon)}|F|, $$
$$|E^{(3)}_{n}|\lesssim \gamma^{-s}\sigma^{s-2}|F|\lesssim 2^{-n(s/2-2)}\lambda^{-s(1/2-\epsilon)}|F|.$$
Define 
$$E^{*}:=\bigcup_{n\ge \Delta}(E^{(1)}_{n}\cup E^{(2)}_{n}\cup E^{(3)}_{n}).$$ Note that since $\Delta\gtrsim \log_2(\lambda^{-1})$, we have $|E^{*}|\lesssim \lambda^{-p}|F|,$ an estimate which can be seen by using a sufficiently large $s$.

For each $x\notin E^{*}$, Theorem ~\ref{mainineq} guarantees that
$$
\|( \sum_{s\in \S_1\atop{|I_s|< 2^k}}\langle1_{F},\varphi_s\rangle\phi_s(x,\theta) )_{k\in\Z} \|_{M_{2,\theta}^*(\R)} \le \sum_{n\ge \Delta}\|(\sum_{s\in \P_n\atop{|I_s|< 2^k}}\langle1_{F},\varphi_s\rangle\phi_s(x,\theta))_{k\in\Z}\|_{M_{2,\theta}^*(\R)}\lesssim$$
$$\lesssim \sum_{n\ge \Delta}n[2^{(3(r/2-1)-1)n}\lambda^{p(r/2-1)}+2^{(3(r/2-1)-1/2)n}\lambda^{p(r/2-1)+1/2-\epsilon}]\lesssim \lambda^{1-2\epsilon},
$$
if $r$ is chosen sufficiently close to 2, depending on $p$ and $\epsilon$. 
This ends the proof of ~\eqref{en.11}, and hence the proof of  ~\eqref{en.6} in the case $\lambda\le 1$.

We next focus on proving ~\eqref{en.6} in the case $\lambda>1$. In the remaining part of the discussion the size will be  understood with respect to the function $\lambda^{-1}1_F$. Proposition ~\ref{g11a00uugg55} implies that $\size(\S)\lesssim \lambda^{-1}$. Define $\Delta:=[-\log_2(\size(\S))].$ Split  $\S$ as before, as a disjoint union $\S=\bigcup_{n\ge \Delta}\P_n,$ where $\size(\P_n)\le 2^{-n}$ and each $\P_n$ consists of a family $\F_{\P_n}$ of trees  satisfying 
\begin{equation}
\label{e:intpikh}
\sum_{\T\in\F_{\P_n}}|I_T|\lesssim 2^{2n}\lambda^{-2}|F|.
\end{equation}

For each $n\ge \Delta$ define $\sigma:=2^{-n}$, $\beta:=2^{(p+1)n}$ and $\gamma:=2^{-n/2}$.
Define also $a_s:=\langle \lambda^{-1}1_F,\varphi_s\rangle$ for each $s\in\P_n$ and note that the collection $\P_n$ together with  the coefficients $(a_s)_{s\in\P_n}$ satisfy the requirements of Theorem ~\ref{mainineq}. Let $\F_{\P_n,l,m}^{(2)}$ the collection of all the 2-quasitrees $\T_{l,m}^{(2)}$ obtained from all the trees $\T\in\F_{\P_n}$ by the procedure described in the beginning of the previous section. Define the corresponding exceptional sets 
\begin{align*}
E^{(1)}_{n}&:=\bigcup_{l\ge 0}\{x:\sum_{\T\in\F_{\P_n}}1_{2^lI_T}(x)>\beta2^{2l}\},\\
E^{(2)}_{n}&:=\bigcup_{l,m\ge 0}\bigcup_{\T\in\F_{\P_n,l,m}^{(2)}}\{x:\|\sum_{s\in \T\atop{|I_s|<2^j}}a_s\phi_{s,\T}^{(\alpha(l,m))}(x,\xi_{\T})\|_{V^r_j(\Z)}>\gamma2^{-l}(|m|+1)^{-2}\},\\
E^{(3)}_{n}&:= \bigcup_{l\ge 0}\bigcup_{\T\in\F_{\P_n,l+1,0}^{(2)}}\{x:\|\sum_{s\in \T\atop{|I_s|<2^j}}a_s\phi_{s,\T}^{(l)}(x,\xi_{\T})\|_{V^r_j(\Z)}>\gamma2^{-l}\}.
\end{align*}
By  ~\eqref{e:intpikh} and the fact that $\lambda\ge 1$ we get
$$|E^{(1)}_{n}|\lesssim 2^{-(p-1)n}\lambda^{-2}|F|.$$
By Theorem ~\ref{BMOtree} and the fact that $\lambda\ge 1$, for each $1<s<\infty$ we get
$$|E^{(2)}_{n}|\lesssim \gamma^{-s}\sigma^{s-2}\lambda^{-2}|F|\lesssim 2^{-n(s/2-2)}\lambda^{-2}|F|, $$
$$|E^{(3)}_{n}|\lesssim \gamma^{-s}\sigma^{s-2}\lambda^{-2}|F|\lesssim 2^{-n(s/2-2)}\lambda^{-2}|F|.$$
Define 
$$E^{*}:=\bigcup_{n\ge \Delta}(E^{(1)}_{n}\cup E^{(2)}_{n}\cup E^{(3)}_{n}).$$ Note that since $\Delta\gtrsim \log_2(\lambda)$, we have $|E^{*}|\lesssim \lambda^{-p}|F|,$ an estimate which can be seen by using a sufficiently large $s$.

For each $x\notin E^{*}$, Theorem ~\ref{mainineq} guarantees that
$$
\|( \sum_{s\in \S\atop{|I_s|< 2^k}}\langle \lambda^{-1}1_{F},\varphi_s\rangle\phi_s(x,\theta))_{k \in \Z}\|_{M_{2,\theta}^*(\R)}\le \sum_{n\ge \Delta}\|(\sum_{s\in \P_n\atop{|I_s|< 2^k}}\langle \lambda^{-1}1_{F},\varphi_s\rangle\phi_s(x,\theta))_{k\in\Z}\|_{M_{2,\theta}^*(\R)}$$
$$
\lesssim \sum_{n\ge \Delta}n2^{(p+1)(r/2-1)n}(2^{-n}+2^{-n/2})\lesssim 1,
$$
if $r$ is chosen sufficiently close to 2, depending only on $p$. This ends the proof of ~\eqref{en.6}
in the case $\lambda>1$.

\subsection{The estimates for $O_{\bf S}$} 
To prove ~\eqref{en.11} and ~\eqref{en.6} in the case $\lambda>1$ for $O_{{\bf S}}$, we use the same values for $a_s$,  $\sigma$, $\beta$ and $\gamma$ as in the case of $T$ and work with the exceptional sets 
\begin{align*}
E^{(1)}_{n}&:=\bigcup_{l\ge 0}\{x:\sum_{\T\in\F_{\P_n}}1_{2^lI_T}(x)>\beta2^{2l}\},\\
E^{(2)}_{n}&:=\bigcup_{l,m\ge 0}\bigcup_{\T\in\F_{\P_n,l,m}^{(2)}}\{x:\|\sum_{s\in \T\atop{|I_s|<2^j}}a_s\phi_{s,\T}^{(\alpha(l,m))}(x,\xi_{\T})\|_{O_{{\bf U}} \cap V^r_j(\Z)}>\gamma2^{-l}(|m|+1)^{-2}\},\\
E^{(3)}_{n}&:= \bigcup_{l\ge 0}\bigcup_{\T\in\F_{\P_n,l+1,0}^{(2)}}\{x:\|\sum_{s\in \T\atop{|I_s|<2^j}}a_s\phi_{s,\T}^{(l)}(x,\xi_{\T})\|_{O_{{\bf U}} \cap V^r_j(\Z)}>\gamma2^{-l}\}.
\end{align*}

\subsection{The estimates for $Q_{{\bf S}}$} 
To prove ~\eqref{en.11} and ~\eqref{en.6} in the case $\lambda>1$ for the operator $Q_{{\bf S}}$, we use the same values for $a_s$,  $\sigma$, $\beta$ and $\gamma$ as in the case of $T$ and define the exceptional set
$$E:= \bigcup_{l\ge 0}\bigcup_{\T\in\F_{l+1,0}^{(2)}}\{x:(\sum_{j\in\Z}|\sum_{s\in \T\atop{|I_s|=2^j}}a_s\phi_{s,\T}^{(l)}(x,\xi_{\T})|^2)^{1/2}>\gamma2^{-l}\}.$$


\end{document}